\documentclass[12pt,letterpaper,leqno]{amsart}
\usepackage{microtype}
\usepackage{amssymb}
\usepackage{amsfonts}
\usepackage{geometry}
\usepackage{hyperref}

\usepackage{mathrsfs}
\usepackage{enumitem}
\usepackage{booktabs,array}
\usepackage{tikz}
\usetikzlibrary{arrows.meta,calc}

\renewcommand{\Im}{\operatorname{Im}}

\newcommand{\supp}{\operatorname{supp}}

\newtheorem{theorem}{Theorem}
\theoremstyle{plain}

\newtheorem{conjecture}[theorem]{Conjecture}

\newtheorem{lemma}[theorem]{Lemma}

\newtheorem{proposition}[theorem]{Proposition}
\newtheorem{remark}[theorem]{Remark}

\numberwithin{equation}{section}
\numberwithin{theorem}{section}
\IfFileExists{tcilatex.tex}{\input{tcilatex}}{}
\numberwithin{figure}{section}
\numberwithin{table}{section}
\newcommand{\Mom}{\mathbf P}
\newcommand{\R}{\mathbb R}\newcommand{\C}{\mathbb C}
\newcommand{\cH}{\mathcal H}\newcommand{\cR}{\mathcal R}
\newcommand{\cD}{\mathcal D}
\newcommand{\cA}{\mathcal A}
\newcommand{\eps}{\varepsilon}\newcommand{\dd}{\,d}
\newcommand{\LE}{\operatorname{LE}}\newcommand{\LG}{\operatorname{LG}}
\newcommand{\OO}{\mathsf O}\newcommand{\Sim}{\operatorname{Sim}}

\DeclareMathOperator{\dist}{dist}

\begin{document}
\title[Radiation and persistent shape dynamics]
  {Radiation and persistent shape dynamics for the five-dimensional energy-critical nonlinear wave equation}

\author[X. Chen]{Xuwen Chen}
\address{Department of Mathematics, University of Rochester, Rochester, NY, USA}
\email{xuwenmath@gmail.com}

\author[J. Holmer]{Justin Holmer}
\address{Department of Mathematics, Brown University, Providence, RI, USA}
\email{\detokenize{justin_holmer@brown.edu}}

\author[S. Shen]{Shunlin Shen}
\address{School of Mathematical Sciences, University of Science and Technology of China, Hefei, 230026, China}
\email{slshen@ustc.edu.cn}


\date{}

\begin{abstract}
We construct a nonradial solution of the five-dimensional focusing energy-critical wave equation that remains bounded in the energy space for all positive times. After subtracting a nonzero outgoing free wave, the solution is asymptotic to a single sign-changing core whose subsequential limits lie in one conformal class and represent a continuum of distinct similarity orbits. Paired incoming radiation packets generate quadratic impulses that move the core along its conformal orbit. This disproves the fixed-profile soliton resolution conjecture and suggests a reformulation in which each soliton has a fixed velocity and a profile that may vary within a fixed conformal class. For forward-global solutions bounded in the energy space, we prove without symmetry assumptions that a decomposition along a single sequence of times tending to infinity into one fully conformally nondegenerate core and a free wave extends to all late times, with the same conformal class, velocity and free wave.
\end{abstract}
\maketitle
\tableofcontents

\section{Introduction}\label{master:intro}
We construct a counterexample to fixed-profile soliton resolution for the
five-dimensional focusing energy-critical wave equation
\begin{equation}\label{master:nlw}
u_{tt}-\Delta u=|u|^{\frac{4}{3}}u,
 \qquad (t,x)\in\R\times\R^5,
\end{equation}
in the energy space
\[
 \cH=\dot H^1(\R^5)\times L^2(\R^5),\qquad
 \|(f,g)\|_{\cH}^2=\|\nabla f\|_2^2+\|g\|_2^2.
\]
We write $\vec u=(u,u_t)$ and denote the free wave group by $U_0(t)$.
By an energy solution on $I$ we mean a real solution of Duhamel's formula with
\[
 \vec u\in C(I;\cH),\qquad
 u\in L_{\mathrm{loc}}^{7/3}(I;L^{14/3}(\R^5)).
\]
The equation conserves the energy
\[
E(u,u_t)=\frac{1}{2}\int_{\R^5}(|\nabla u|^2+|u_t|^2)\,dx
             -\frac{3}{10}\int_{\R^5}|u|^{\frac{10}{3}}\,dx.
\]
Our solution remains close to a sign-changing stationary state $Q_*$
for all $t\ge0$. After we subtract a free wave, the core has two distinct
subsequential limits in the same conformal class. Their distance remains
positive modulo translations, dilations and rotations. In fact, a continuum
of similarity classes occurs along different sequences of times.

The soliton resolution problem asks whether every forward-global solution bounded in $\cH$
separates into traveling waves and free radiation as $t\to\infty$.
In \cite[Conjecture~5.5]{DKM15}, Duyckaerts, Kenig and
Merle fix the stationary profiles and their Lorentz velocities before
choosing an arbitrary sequence of times.
The spatial parameters may depend on the sequence, and the parameters
of distinct profiles must be asymptotically orthogonal. In dimensions
three through five, Duyckaerts, Jia, Kenig and Merle \cite{DJKM17}
prove a decomposition along a suitable sequence for solutions bounded
in $\cH$. We ask whether, in the nonradial problem, one finite list of
profiles suffices for every sequence.

The stationary equation
\begin{equation}\label{intro:stationary}
 -\Delta Q=|Q|^{\frac{4}{3}}Q
\end{equation}
has a larger symmetry group than the spatial similarity group of
\eqref{master:nlw}. Its conformal orbit may therefore contain stationary
states with different shapes modulo similarities. We use the symmetric
single-node state in the isoparametric setting of Fern\'andez and
Petean \cite{FP20}, studied by Neves and Pistoia \cite{NP26}.
For this state $Q$, a mode decomposition and the spectral bound of
Appendix~\ref{app:single-node} give
\[
 L_Q=-\Delta-\frac{7}{3}|Q|^{\frac{4}{3}},\qquad
 \ker_{\dot H^1}L_Q=T_Q\bigl(\operatorname{Conf}(\mathbb S^5)Q\bigr),
 \qquad \dim\ker_{\dot H^1}L_Q=15
\]
(Lemma~\ref{np:nondegenerate}). We obtain the spectral bound by ODE
comparison and exact polynomial inequalities. We choose a conformal
representative at which incoming radiation produces quadratic impulses
in every kernel direction. One of these directions is transverse to
the similarity orbit, and we use it to prescribe the changes of shape.

We denote the spatial similarity group by $\Sim$, including reflections.
For $\lambda>0$, $a\in\R^5$ and $R\in O(5)$, its action on energy pairs is
\[
 \mathcal T_{\lambda,a,R}(f,g)=
 \left(\lambda^{-\frac{3}{2}}f\left(\frac{R^{-1}(x-a)}{\lambda}\right),
       \lambda^{-\frac{5}{2}}g\left(\frac{R^{-1}(x-a)}{\lambda}\right)\right).
\]
Since $\mathcal T_{\lambda,a,R}$ is an isometry of $\cH$, we may
measure the separation of two cores by the distance to a similarity orbit.

\begin{theorem}[Main Theorem]\label{master:main}
There exist a sign-changing stationary solution $Q_*$, a $C^2$ local
chart $q_\xi$ in its conformal orbit, $\xi\in\R^{15}$, with $q_0=Q_*$, and a unit
vector $e\in\R^{15}$ with the following property. For every sufficiently
small $d>0$, there exist a nonradial forward-global energy solution $u$ of
\eqref{master:nlw}, a $C^1$ curve $\xi:[0,\infty)\to\R^{15}$, and a free
wave $v_L$ such that
\begin{gather}
 \sup_{t\ge0}\left(\|u(t)-Q_*\|_{\dot H^1}+\|u_t(t)\|_2\right)\lesssim d,
 \qquad \|\xi\|_\infty\lesssim d,\label{master:bounds}\\
 \|\vec u(t)-(q_{\xi(t)},0)-\vec v_L(t)\|_{\cH}\longrightarrow0.
 \label{master:radiation}
\end{gather}
Along two sequences tending to infinity, the stationary core converges
in $\dot H^1$ to $Q_*$ and $q_{de}$, respectively. These limits satisfy
\begin{equation}\label{master:anchors}
 \dist_{\dot H^1}(q_{de},\Sim\cdot Q_*)\gtrsim d.
\end{equation}
The omega-limit set
\begin{equation}\label{master:omega}
 \begin{aligned}
 \Omega
 &=\bigcap_{T>0}\overline{\{q_{\xi(t)}\mid t\ge T\}}^{\dot H^1}\\
 &=\{Q\in\dot H^1\mid q_{\xi(t_n)}\to Q\text{ in }\dot H^1
       \text{ for some }t_n\to\infty\}
 \end{aligned}
\end{equation}
is compact and connected, lies in the conformal orbit of $Q_*$,
and contains a continuum of distinct similarity classes.

There exists a sequence $t_n\to\infty$ on which no fixed finite family
of nonzero energy pairs represents $\vec u(t_n)$, up to one free wave
and an $o_{\cH}(1)$ error, by a uniformly bounded number of pairwise
orthogonal spatial transforms. The active entries may vary with $n$.
In particular, the solution violates the fixed-profile conclusion of
\cite[Conjecture~5.5]{DKM15}.
\end{theorem}

\begin{remark}[A continuum of similarity classes]\label{master:main-remark}
For the single solution in Theorem~\ref{master:main}, the set of core
limits $\Omega$ in \eqref{master:omega} satisfies
\begin{equation}\label{master:omega-interval}
 [0,\dist_{\dot H^1}(q_{de},\Sim\cdot Q_*)]
 \subseteq
 \{\dist_{\dot H^1}(Q,\Sim\cdot Q_*)\mid Q\in\Omega\}.
\end{equation}
Since similarities act isometrically, the function
\[
 Q\longmapsto\dist_{\dot H^1}(Q,\Sim\cdot Q_*)
\]
is $1$-Lipschitz and constant on each similarity orbit. The connectedness
of $\Omega$ and the limits $Q_*,q_{de}\in\Omega$ imply
\eqref{master:omega-interval}. The interval has length $\gtrsim d$ by
\eqref{master:anchors}. For each $r$ in this interval, we choose $Q_r\in\Omega$ so that
\begin{equation}\label{master:class-separation}
 \begin{gathered}
 \dist_{\dot H^1}(Q_r,\Sim\cdot Q_*)=r,\\
 \dist_{\dot H^1}(Q_r,\Sim\cdot Q_s)\ge|r-s|.
 \end{gathered}
\end{equation}
The $Q_r$ therefore belong to distinct similarity orbits, even if we
allow reflections. In particular, no countable union of such orbits
contains $\Omega$.

In Section~\ref{i5:i5:obstruction}, we use uniqueness of the radiation
and Lemma~\ref{lem:profile-collapse} to exclude every fixed finite family
of nonzero energy pairs. The selection may depend on time, provided
that the number of active entries stays bounded and their parameters
are pairwise orthogonal.
\end{remark}

Every limit in \eqref{master:omega} belongs to the conformal orbit of
$Q_*$. This suggests fixing the conformal class of each soliton while
allowing its representative to depend on time. More precisely, let
$G=\operatorname{Conf}(\mathbb S^5)$ act on stationary states by
stereographic projection, with critical weight
$|\det Dg^{-1}|^{\frac{3}{10}}$. We set $\mathcal K_Q=O(6)Q$. Its
translations and positive dilations exhaust $GQ$. Indeed, given $g\in G$,
we choose $O\in O(6)$ with $O(g^{-1}(\infty))=\infty$. Then $gO^{-1}$
fixes $\infty$ and acts by a Euclidean similarity, whose orthogonal factor
we absorb into $O(6)$. For $\ell\in\R^5$
with $|\ell|<1$, we write
\[
 A_\ell=(I-\ell\otimes\ell)^{-\frac{1}{2}},\qquad
 \phi_q(y)=q(A_\ell y),\qquad
 \mathbf Q_\ell(q)=(\phi_q,-\ell\cdot\nabla\phi_q).
\]
Thus $\mathbf Q_\ell(q)$ is the initial pair of the traveling wave
$\phi_q(x-\ell t)$. We put $\mathcal T_{\lambda,a}=\mathcal T_{\lambda,a,I}$
and denote the family with class $GQ$ and velocity $\ell$ by
\[
 \mathcal M_\ell(Q)=
 \{\mathcal T_{\lambda,a}\mathbf Q_\ell(q):
            \lambda>0,\ a\in\R^5,\ q\in\mathcal K_Q\}.
\]

\begin{conjecture}[Soliton resolution by conformal classes]\label{conj:classes}
Let $u$ be a forward-global strong solution of \eqref{master:nlw}
bounded in $\cH$. There exist an integer $J\ge0$, conformal classes
$GQ_1,\ldots,GQ_J$ of nonzero stationary states, velocities
$|\ell_j|<1$ and a finite-energy free wave $v$, all fixed, such that
\[
 \vec u(t)=\vec v(t)+\sum_{j=1}^{J}
     \mathcal T_{\lambda_j(t),a_j(t)}\mathbf Q_{\ell_j}(q_j(t))
     +o_{\cH}(1)\qquad(t\to\infty).
\]
Here $q_j(t)\in\mathcal K_{Q_j}$, and the scales and centers of distinct
summands satisfy \eqref{obstruction:orthogonality}.
\end{conjecture}

We prove the conjecture in the nonradial case when a single
nondegenerate core occurs along one sequence of times. We use the full
conformal nondegeneracy condition of \cite{DKM16}. In particular,
$T_Q(GQ)$ includes the infinitesimal symmetries obtained by conjugating
translations with Kelvin inversion.

\begin{theorem}[One nondegenerate core]\label{thm:sequence}
Let $u$ be a forward-global strong solution of \eqref{master:nlw}
bounded in $\cH$. Let $Q\ne0$ solve \eqref{intro:stationary} and satisfy
\[
 \ker_{\dot H^1} L_Q=T_Q(GQ),\qquad
 L_Q=-\Delta-\frac{7}{3}|Q|^{\frac{4}{3}} .
\]
Suppose that for some $|\ell|<1$, a finite-energy free wave $v$ and
a sequence $t_n\to\infty$,
\[
 \dist_{\cH}\bigl(\vec u(t_n)-\vec v(t_n),\mathcal M_\ell(Q)\bigr)
 \longrightarrow0 .
\]
Then $U_0(-t)\vec u(t)\rightharpoonup\vec v(0)$ in $\cH$.
There exist $T>0$ and $C^1$ functions $\lambda(t)>0$, $a(t)\in\R^5$,
$q(t)\in\mathcal K_Q$ on $[T,\infty)$ such that, with the same $v$,
\begin{gather*}
 \|\vec u(t)-\vec v(t)
       -\mathcal T_{\lambda(t),a(t)}\mathbf Q_\ell(q(t))\|_{\cH}
       \longrightarrow0,\\
 \lambda'(t)\longrightarrow0,\qquad a'(t)\longrightarrow\ell,\qquad
 \lambda(t)\|q'(t)\|_{\dot H^1\cap L^2}\longrightarrow0 .
\end{gather*}
In particular, $\lambda(t)=o(t)$ and $a(t)=\ell t+o(t)$.
\end{theorem}

We prove Theorem~\ref{thm:sequence} in Section~\ref{sec:one-core}.
The state $Q_*$ of Theorem~\ref{master:main} satisfies its nondegeneracy
hypothesis by Lemma~\ref{np:nondegenerate} and conformal covariance.
The even times in \eqref{global:stationary-limits} verify the
sequence hypothesis with $(Q,\ell,v)=(Q_*,0,v_L)$. For any fixed
stationary state $Q'\ne0$, the triangle inequality along the even
and odd times implies
\[
 \limsup_{t\to\infty}
 \dist_{\cH}\bigl(\vec u(t)-\vec v_L(t),\Sim\cdot(Q',0)\bigr)
 \ge\frac{1}{2}\dist_{\dot H^1}(q_{de},\Sim\cdot Q_*)\gtrsim d .
\]
Thus, even under full conformal nondegeneracy, the representative
$q(t)$ in Theorem~\ref{thm:sequence} must in general depend on time.

\subsection{Related results and the role of the stationary family}\label{intro:related-results}
Bahouri and G\'erard \cite{BG99} established the wave profile decomposition
in dimension three, and Bulut \cite{Bulut10} extended it to all dimensions
at least three. Kenig and Merle \cite{KM08} developed the
concentration-compactness and rigidity method for the focusing equation.
Duyckaerts, Kenig and Merle proved radial all-time soliton resolution
in dimension three \cite{DKM13} and then in all odd dimensions at
least three \cite{DKM23}. Jendrej and Lawrie \cite{JL23} proved the
radial result in every dimension at least four. In the nonradial
problem, Collot, Duyckaerts, Kenig and Merle \cite{CDKM26} consider
compactly supported data with ground-state traveling waves in a
sequential decomposition. They obtain all-time resolution in the
one-wave case and in the two-wave case with distinct velocities.
For three or more waves with distinct velocities, they obtain an
alternative between all-time resolution and infinitely many collisions.

The choice of profiles can also depend on the time sequence in
geometric evolution equations. Topping
\cite{Topping04} constructed finite-time singularities of harmonic map
heat flow with nonunique limiting bubbles for certain target manifolds.
Jendrej, Lawrie and Schlag \cite{JLS25} distinguish convergence to the
family of bubble configurations from the existence of a choice of
harmonic maps independent of the time sequence, and cite Topping's
example in discussing the latter question. The related discussion in
Collot, Duyckaerts, Kenig and Merle \cite{CDKM26} refers to Topping's
earlier rigidity paper \cite{Topping97}.

For wave maps, Engelstein and Mendelson \cite{EM22} use the regularity
and bubbling theory of Sterbenz and Tataru \cite{ST10a,ST10} to construct
solutions into a compact smooth target with a continuum of distinct
harmonic-map bubbles along different time sequences. For the scalar
wave equation \eqref{master:nlw}, we distinguish the limiting cores
modulo spatial similarities.
For harmonic map heat flow, under an additional symmetry assumption,
Gustafson, Nakanishi and Tsai \cite{GNT10} construct global degree-two
solutions whose scale oscillates indefinitely while the solution
approaches the harmonic-map family.

Ding \cite{Ding86} constructed sign-changing stationary solutions.
Del Pino, Musso, Pacard and Pistoia \cite{DPMPP11} constructed families
of nodal states, and Musso and Wei \cite{MW} studied their nondegeneracy. Duyckaerts, Jia, Kenig and Merle
\cite{DJKM17} discuss these nodal constructions among the possible
stationary profiles in subsequential resolution.

Duyckaerts, Kenig and Merle
\cite[(1.7)--(1.8), Sections~3.1 and~3.3]{DKM16} incorporated the full
conformal orbit of the stationary equation into their nondegeneracy
condition and modulation analysis, including the directions obtained
by conjugating translations with Kelvin inversion. Their rigidity results assume nondegeneracy of the stationary state.
Once we verify this hypothesis, their Proposition~4.1 also gives a
qualitative argument for nonzero radiation. In Appendix~\ref{app:energy},
we obtain a lower bound from the incoming packet energies and the
energy expansion.

For the stationary input, we construct a symmetric single-node state in
the isoparametric setting of Fern\'andez and Petean \cite{FP20}.
Neves and Pistoia \cite{NP26} prove nondegeneracy in the invariant
space in the odd dimensions they consider, and in the full space in
sufficiently large odd dimensions. In dimension five, we verify full nondegeneracy directly in
Lemma~\ref{np:nondegenerate} and Appendix~\ref{app:single-node},
independently of the Musso--Wei theorem.

Krieger, Nakanishi and Schlag \cite{KNS15} construct a center-stable
manifold near the ground state. Martel and Merle \cite{MM16} construct
five-dimensional multi-solitons. Their more recent construction
\cite{MM25} allows arbitrary directions with distinct subluminal velocities
and proves inelasticity under a non-cancellation assumption.
Yuan \cite{Yuan21} constructs excited multi-solitons with distinct
collinear velocities from nondegenerate stationary states satisfying
suitable decay assumptions. Conformal transformations of the nodal family
of del Pino, Musso, Pacard and Pistoia provide faster-decaying examples,
which serve as fixed shapes in the asymptotic decomposition. In dimension four,
Yuan \cite{Yuan22} also observes formally that the long-range
interaction between separating ground-state solitons can give rise to
divergent modulation terms.

Krieger, Schlag and Tataru \cite{KST09} construct radial type-II blowup
solutions for the three-dimensional critical wave equation with a
family of prescribed concentration rates. In the five-dimensional
radial wave equation, Jendrej \cite{J17} prescribes the regular part
of a solution at blowup and uses it as a source for the concentration
of the ground state. For equivariant wave maps, Jendrej, Lawrie and
Rodriguez \cite{JLR22} prescribe the radiation and determine the
resulting blowup rate. In these constructions the evolving parameter
is the scale. Donninger and Krieger \cite{DK13} construct global
solutions of the three-dimensional wave equation whose ground-state
scale varies as $t\to\infty$. Pillai \cite{Pillai21} constructs global
radial solutions with admissible time-dependent scales, including
oscillatory scales, together with one free wave and a decaying correction.
These core trajectories remain in one similarity orbit. Kadar \cite{Kadar26}
constructs solutions near a soliton with prescribed, sufficiently decaying
radiation at null infinity.

In our construction, the quadratic interaction with radiation moves the
core in a conformal direction transverse to the similarity orbit. The
stabilizer of the nodal state and properness of the similarity action
give the separation \eqref{master:anchors}. Since the core trajectory is
continuous and precompact, its limit set then satisfies
\eqref{master:omega-interval}. For $\vec v_L(t)=U_0(t)W_+$,
Appendix~\ref{app:energy} gives
\begin{equation*}
 \|W_+\|_{\mathcal H}\asymp d^3,\qquad
 E(u,u_t)-E(Q_*,0)\asymp d^6.
\end{equation*}

\subsection{Outline of the proofs}
We prove Theorem~\ref{master:main} in four steps.

\begin{enumerate}[label=(\arabic*),ref=(\arabic*),leftmargin=*,itemsep=0.7ex]
\item\label{outline:decomposition}\textbf{Construction of the core.}
We seek a forward-global solution satisfying \eqref{master:bounds} and
\[
 \vec u(t)=\underbrace{(q_{\xi(t)},0)}_{\text{stationary core}}
 +\underbrace{\vec v_L(t)}_{\text{free radiation}}+o_{\cH}(1),
 \qquad \vec v_L(t)=U_0(t)W_+.
\]
The map $t\mapsto q_{\xi(t)}$ is continuous in $\dot H^1$, and
\begin{equation}\label{outline:core-compactness}
 \overline{\{q_{\xi(t)}\mid t\ge0\}}^{\dot H^1}
 \quad\text{is compact in }\dot H^1\setminus\{0\}.
\end{equation}
Each $q_{\xi(t)}$ solves \eqref{intro:stationary}. On finite intervals,
we obtain uniform bounds and strongly convergent initial data
(Sections~\ref{part:linear}--\ref{part:endpoint}). Local stability and
the source estimates then give the global solution and $W_+$
(Section~\ref{sec:global-radiation}).

\item\label{outline:limits}\textbf{Separation of the two limits.}
For $T_j=4^jd^{-5}$, we require
\[
 \begin{gathered}
 q_{\xi(T_{2k})}\longrightarrow Q_*,\qquad
 q_{\xi(T_{2k+1})}\longrightarrow q_{de}
 \quad\text{in }\dot H^1\quad(k\to\infty),\\
 \dist_{\dot H^1}(q_{de},\Sim\cdot Q_*)\gtrsim d>0.
 \end{gathered}
\]
In Section~\ref{part:background}, we choose the transverse direction
$e$ and prove the separation. We take
\[
 \eta_j=d^3 2^{-j},\qquad
 \sum_{j\ge0}\eta_j\lesssim d^3,\qquad T_j\eta_j^2=d.
\]
If we place the two normalized impulses $(x,y)=\pm(e,-e)$ at times
$2T_j$ and $3T_j$, respectively, we obtain
\[
 \underbrace{\eta_j^2(x+y)}_{\text{velocity increment}}=0,
 \qquad
 \underbrace{T_j\eta_j^2(2x+y)}_{\text{displacement}}=\pm de.
\]
We realize these displacements by force estimates and endpoint selection
(Sections~\ref{part:catalogue}--\ref{part:endpoint}). The endpoints pass
to the global solution in \eqref{global:stationary-limits}
(Section~\ref{sec:global-passage}).

\item\label{outline:radiation}\textbf{Uniqueness of the radiation.}
We consider a competing decomposition with fixed profiles,
\[
 \begin{gathered}
 \vec u(t)=U_0(t)\widetilde W_+
 +\sum_{\ell=1}^{J(t)}\mathcal T_{\ell,t}F_{\nu_\ell(t)}+o_{\cH}(1),\\
 F_1,\ldots,F_M\in\cH\setminus\{0\},\qquad
 J(t)\le J_*,\qquad \nu_\ell(t)\in\{1,\ldots,M\},
 \end{gathered}
\]
where $\mathcal T_{\ell,t}$ is a spatial similarity whose parameters
satisfy \eqref{obstruction:orthogonality}. We subtract
\eqref{master:radiation} and pair with $U_0(t)(W_+-\widetilde W_+)$.
By \eqref{outline:core-compactness} and the uniform separation estimate
\eqref{short:radiation-separation}, we obtain
\[
 \|W_+-\widetilde W_+\|_{\cH}^2=o(1),\qquad
 \widetilde W_+=W_+.
\]
Section~\ref{sec:radiation-comparison} gives the same conclusion
along any sequence tending to infinity.

\item\label{outline:obstruction}\textbf{Exclusion of a fixed family.}
The limit set $\Omega$ in \eqref{master:omega} satisfies
\[
 \{Q_*,q_{de}\}\subset\Omega\subset\dot H^1\setminus\{0\}.
\]
Continuity and \eqref{outline:core-compactness} give compactness and
connectedness of $\Omega$. Hence \eqref{master:omega-interval} gives
a continuum of similarity classes.

For each $q\in\Omega$, we choose $t_n\to\infty$ such that
$q_{\xi(t_n)}\to q$ in $\dot H^1$. Along a subsequence, the number $J$
and indices $\nu_\ell$ of the active profiles are fixed. We write
$\mathcal T_{\ell,n}=\mathcal T_{\ell,t_n}$. After subtracting
$\vec v_L(t_n)$, we apply Lemma~\ref{lem:profile-collapse} to obtain
\[
 \sum_{\ell=1}^{J}\mathcal T_{\ell,n}F_{\nu_\ell}
 \longrightarrow(q,0)\ne0
 \quad\Longrightarrow\quad
 J=1,\quad (q,0)\in\Sim\cdot F_{\nu_1}.
\]
With a fixed list of active profiles, the even and odd times would
give $(Q_*,0),(q_{de},0)\in\Sim\cdot F_1$, contrary to
\eqref{master:anchors}. Even if the indices vary, we would still have
\[
 \{(q,0)\mid q\in\Omega\}
 \subseteq\bigcup_{m=1}^{M}\Sim\cdot F_m.
\]
This contradicts \eqref{master:omega-interval}. We give the details in
Sections~\ref{sec:connected-limits}--\ref{sec:profile-collapse}.
\end{enumerate}

\begin{figure}[!htbp]
\centering
\begin{tikzpicture}[
  font=\small,
  proofbox/.style={draw,rounded corners=2pt,
    line width=0.45pt,align=center,inner sep=5pt,text width=14.3cm},
  branch/.style={proofbox,text width=6.5cm,minimum height=3.6cm,anchor=north},
  route/.style={-{Stealth[length=2mm]},line width=0.6pt,draw}
]
\node[proofbox] (decomposition) at (0,0) {
  \textbf{The stationary core $(q_{\xi(t)},0)$}\quad Step~\ref{outline:decomposition}\\[2pt]
  $\vec u=(q_\xi,0)+\vec v_L+o_{\cH}(1)$, with $\vec v_L=U_0(t)W_+$\\[2pt]
  Uniform finite-interval construction and global passage in
  Sections~\ref{part:linear}--\ref{sec:global-radiation}
};
\node[branch] (limits) at ([xshift=-3.8cm,yshift=-5mm]decomposition.south) {
  \textbf{Limits $Q_*$ and $q_{de}$}\\
  Step~\ref{outline:limits}, Sections~\ref{part:background},
  \ref{part:endpoint}--\ref{sec:global-radiation}\\[3pt]
  $q_{\xi(T_{2k})}\to Q_*$,\quad $q_{\xi(T_{2k+1})}\to q_{de}$\\[3pt]
  $\dist_{\dot H^1}(q_{de},\Sim\cdot Q_*)>0$\\[3pt]
  $\displaystyle\sum_j\eta_j\lesssim d^3,\qquad T_j\eta_j^2=d$
};
\node[branch] (radiation) at ([xshift=3.8cm,yshift=-5mm]decomposition.south) {
  \textbf{The unique free state $W_+$}\\
  Step~\ref{outline:radiation}, Section~\ref{sec:radiation-comparison}\\[3pt]
  Any competing fixed-profile resolution\\[3pt]
  $\Longrightarrow\ \widetilde W_+=W_+$\\[3pt]
  We subtract $U_0(t)W_+$.
};
\node[proofbox,anchor=north] (obstruction) at ([xshift=-3.8cm,yshift=-5mm]radiation.south) {
  \textbf{Profile collapse and the contradiction}\quad
  Step~\ref{outline:obstruction}, Sections~\ref{sec:connected-limits}--\ref{sec:profile-collapse}\\[3pt]
  The two separated limits exclude a single similarity orbit.\\[3pt]
  Connectedness of $\Omega$ and \eqref{master:omega-interval}
  give a continuum of similarity classes.\\[3pt]
  A fixed finite family would give
  $\{(q,0):q\in\Omega\}\subseteq\bigcup_{m=1}^M\Sim\cdot F_m$.
};
\draw[route] (decomposition.south) -- ++(0,-0.25) -| (limits.north);
\draw[route] (decomposition.south) -- ++(0,-0.25) -| (radiation.north);
\draw[route] (limits.south) -- ++(0,-0.25) -| (obstruction.north);
\draw[route] (radiation.south) -- ++(0,-0.25) -| (obstruction.north);
\end{tikzpicture}
\caption{Four goals for the fixed-profile obstruction in Theorem~\ref{master:main}.
The limit set $\Omega$ is defined by \eqref{master:omega} in $\dot H^1$.
Its energy-space counterpart is $\{(q,0)\mid q\in\Omega\}$.}
\label{fig:proof-outline}
\end{figure}
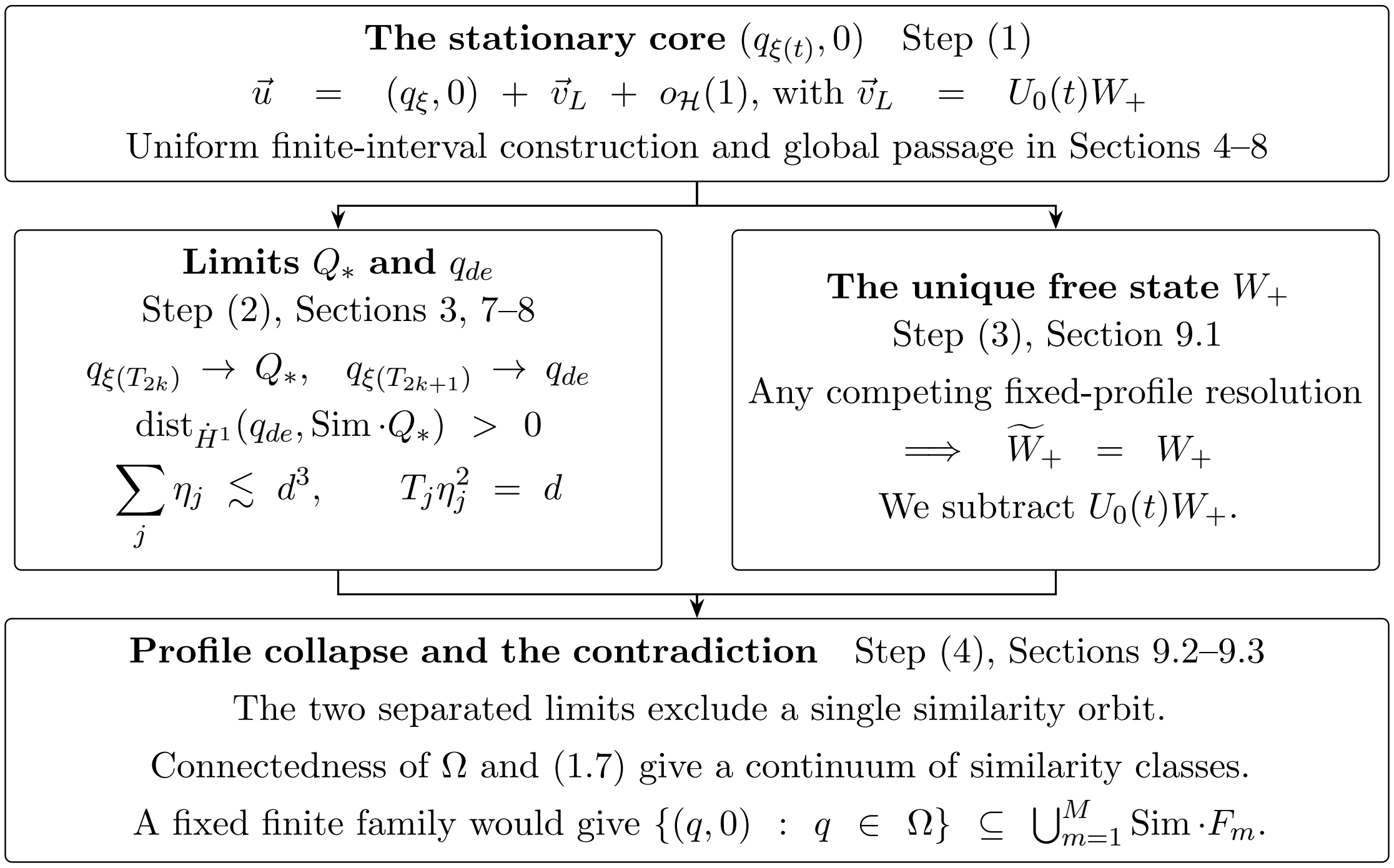

Figure~\ref{fig:proof-outline} records the four steps for
Theorem~\ref{master:main}. Table~\ref{tab:proof-sections} lists the sections used in the
proofs of Theorems~\ref{master:main} and~\ref{thm:sequence}.
\begin{table}[!htbp]
\caption{Sections used in the proofs of the two main theorems.}
\label{tab:proof-sections}
\centering
{\renewcommand{\arraystretch}{1.1}
\begin{tabular}{@{}>{\raggedright\arraybackslash}p{0.14\textwidth}>{\raggedright\arraybackslash}p{0.82\textwidth}@{}}
\toprule
\textbf{Sections} & \textbf{Role in the proof}\\
\midrule
\ref{master:proof} & We deduce Theorem~\ref{master:main} from the main propositions.\\[2pt]
\ref{part:background} & We prove orbit separation and controllability for step~\ref{outline:limits}.\\[2pt]
\ref{part:linear}--\ref{part:nonlinear} & We develop the linear estimates, radiation controls and nonlinear correction.\\[2pt]
\ref{part:endpoint} & We construct finite-interval solutions with the prescribed endpoints.\\[2pt]
\ref{sec:global-radiation} & We obtain the global decomposition and two limits, completing steps~\ref{outline:decomposition} and~\ref{outline:limits}.\\[2pt]
\ref{sec:radiation-comparison} & We identify the radiation in step~\ref{outline:radiation}.\\[2pt]
\ref{sec:connected-limits}--\ref{sec:profile-collapse} & We prove the continuum of classes and profile collapse, completing step~\ref{outline:obstruction}.\\[2pt]
\ref{sec:one-core} & We extend the single-core decomposition to all large times
in Theorem~\ref{thm:sequence}.\\
\bottomrule
\end{tabular}}
\end{table}

\medskip
\noindent\emph{One nondegenerate core.}
To prove Theorem~\ref{thm:sequence} in Section~\ref{sec:one-core}, we
start from the decomposition along $t_n\to\infty$. By the radiation
comparison in Step~\ref{outline:radiation}, we identify $v$ with the
canonical free wave (Section~\ref{sc:radiation}). We subtract $\vec v(t)$
and choose the scale, center and conformal representative so that the
error is orthogonal to the tangent space of $\mathcal M_\ell(Q)$.
In Section~\ref{sc:geometry}, we use full conformal nondegeneracy to
identify this tangent space with the kernel of the Hessian of
$E+\ell\cdot\Mom$, where $\Mom$ is the momentum. In
Section~\ref{sc:exterior}, we detect the growing modes by their energy
on an exterior cone and prove the coercivity estimate \eqref{sc:coerc}.
In Section~\ref{sc:dynamics}, we integrate the damped mode equations
in rescaled time. With $\eta$ as in \eqref{sc:tails}, we obtain
\[
 \begin{gathered}
 \mathrm d\tau=\frac{\mathrm dt}{\lambda(t)},\qquad
 d(t)=\dist_{\cH}\bigl(\vec u(t)-\vec v(t),\mathcal M_\ell(Q)\bigr),\\
 \sup_{t\ge t_n} d(t)\lesssim d(t_n)+\eta(t_n)\longrightarrow0 .
 \end{gathered}
\]
Thus the same radiation $v$, velocity $\ell$ and conformal class $GQ$
describe the solution throughout the intervals $[t_n,t_{n+1}]$.

\section{Main propositions and proof of Theorem~\ref*{master:main}}\label{master:proof}

We deduce Theorem~\ref{master:main} from three propositions and the
obstruction in Proposition~\ref{prop:abstract-obstruction}. Throughout
the construction, we write
\begin{equation*}
f(s)=|s|^{\frac{4}{3}}s,\qquad D=\sqrt{-\Delta},\qquad G_0=D^{-2},
 \qquad \rho=(1+|x|^2)^{\frac{1}{2}},
\end{equation*}
and
\begin{equation*}
X_s=L^2(\rho^{-2s}dx),\qquad Y_s=L^2(\rho^{2s}dx),
 \qquad S(J)=L_t^{\frac{7}{3}}L_x^{\frac{14}{3}}(J\times\R^5).
\end{equation*}
We take the columns of $K$ to be an $L^2$-orthonormal basis of the full
linearized kernel. The compactly supported tests $\Phi$ of
Section~\ref{presentation:linear-setting} satisfy $\Phi K=I$ and
annihilate the negative modes. We fix the background and the radiation
family before choosing $d$ and the number $N$ of intervals. Unless we
indicate otherwise, the implicit constants are independent of $d$ and $N$.

We prove the following geometric statement in Section~\ref{part:background}.

\begin{proposition}[Stationary geometry]\label{presentation:geometry}
There are a sign-changing stationary solution $Q_*$, a local stationary
chart
\begin{equation*}
q_\xi=Q_*+K\xi+R(\xi),\qquad \xi\in\R^{15},\qquad
K^*R=0,\quad R(0)=DR(0)=0,
\end{equation*}
and a unit vector $e\in\R^{15}$ with the following properties.
The chart lies in the conformal orbit of $Q_*$. The columns of $K$ span
the full $\dot H^1$ kernel of
\begin{equation*}
 L=-\Delta-f'(Q_*).
\end{equation*}
The chart satisfies
\begin{equation*}
 q_\xi\in C^2_\xi H^1,\qquad D^2q_\xi\in C^2_\xi Y_4,
\end{equation*}
and, uniformly for $|\alpha|\le2$,
\begin{equation*}
 |\partial_\xi^\alpha q_\xi|\lesssim\rho^{-3},\qquad
 |\nabla\partial_\xi^\alpha q_\xi|\lesssim\rho^{-4}.
\end{equation*}
For all sufficiently small $d>0$,
\begin{equation*}
 \operatorname{dist}_{\dot H^1}(q_{de},\Sim\cdot Q_*)\gtrsim d.
\end{equation*}
Moreover, there are finitely many lines $b_\nu+\R\omega_\nu$ such that
\begin{equation*}
 j_\nu=\frac{1}{4}\int_\R
       f''(Q_*(b_\nu+s\omega_\nu))K(b_\nu+s\omega_\nu)\,ds,
 \qquad 0\in\operatorname{int}\operatorname{conv}\{j_\nu\}.
\end{equation*}
\end{proposition}

The correction $\alpha_N$ below removes the oscillations of the
incoming waves.

\begin{proposition}[Construction on finite time intervals]
\label{presentation:finite}
We fix the background and chart of Proposition~\ref{presentation:geometry}.
For sufficiently small $d>0$, we set
\begin{equation*}
 \begin{gathered}
 \varepsilon=d^3,\qquad \eta_j=\varepsilon2^{-j},\qquad T_j=4^jd^{-5},\\
 \theta_j=\begin{cases}0,&j\text{ even},\\de,&j\text{ odd}.
 \end{cases}
 \end{gathered}
\end{equation*}
For every $N\ge1$, there is an exact critical energy solution
$u_N=q_{\xi_N}+w_N$ on $[0,T_N]$, with $\xi_N\in C^1$ and $\Phi w_N=0$,
such that, with all norms below taken on $[0,T_N]$,
\begin{gather}
 \|\xi_N\|_\infty\lesssim d,\label{presentation:finite-chart}\\
 \|\nabla w_N\|_{L_t^\infty L_x^2}
 +\|\partial_tu_N\|_{L_t^\infty L_x^2}
 +\|w_N\|_{S}+\|w_N\|_{L_t^2X_2}
 +\|\xi_N'\|_{L_t^2}+\|\xi_N'\|_\infty
 \lesssim \varepsilon.\label{presentation:finite-bounds}
\end{gather}
There are functions $\alpha_N\in C_b^1(\R;\R^{15})$ for which
$Y_N=\xi_N-\alpha_N$ satisfies
\begin{equation*}
Y_N(0)=0,\qquad Y_N(T_j)=\theta_j\quad(1\le j\le N).
\end{equation*}
After passage to a subsequence as $N\to\infty$, the physical initial
data $\vec u_N(0)$ converge strongly in $\mathcal H$, and $\alpha_N$
converge in $C_b^1(\R)$ to a function $\alpha\in L^2(0,\infty)$.
All constants and the smallness threshold for $d$ are independent of $N$.
\end{proposition}

We prove Proposition~\ref{presentation:finite} in
Section~\ref{part:endpoint}, and the following scattering criterion in
Section~\ref{sec:scattering-limits}.

\begin{proposition}[Scattering criterion]\label{presentation:scattering}
Let $u=q_\xi+w$ be a forward-global energy solution, where $\xi$ is a
$C^1$ curve in a fixed compact subset of the chart in
Proposition~\ref{presentation:geometry}. If
\begin{equation*}
 \begin{aligned}
 A:={}&\sup_{t\ge0}\bigl(\|\nabla w(t)\|_2+\|u_t(t)\|_2\bigr)\\
 &+\|w\|_{S(0,\infty)}+\|w\|_{L_t^2X_2(0,\infty)}
       +\|\xi'\|_{L^2(0,\infty)}<\infty,
 \end{aligned}
\end{equation*}
then there is a unique free wave $v_L$ such that
\begin{equation*}
 \|(w(t),u_t(t))-(v_L(t),v_{L,t}(t))\|_{\mathcal H}\longrightarrow0.
\end{equation*}
Moreover,
\begin{equation*}
 \|\vec v_L(0)\|_{\mathcal H}\lesssim A+A^2+A^{\frac{7}{3}}.
\end{equation*}
\end{proposition}

\begin{proof}[Proof of Theorem~\ref{master:main}]
By Proposition~\ref{presentation:finite} and
Section~\ref{sec:global-passage}, we obtain a forward-global solution with
\begin{equation*}
 \|\xi\|_\infty\lesssim d,\qquad A\lesssim d^3,\qquad
 \xi(T_j)-\alpha(T_j)=\theta_j\quad(j\ge1).
\end{equation*}
By the chart bounds and Proposition~\ref{presentation:scattering},
we have \eqref{master:bounds}--\eqref{master:radiation}. Since
$\alpha\in C_b^1\cap L^2$, we have $\alpha(t)\to0$, and hence
\begin{equation}\label{global:stationary-limits}
 q_{\xi(T_{2k})}\longrightarrow Q_*,\qquad
 q_{\xi(T_{2k+1})}\longrightarrow q_{de}
 \quad\text{in }\dot H^1.
\end{equation}
These limits are separated by Proposition~\ref{presentation:geometry}.
If $Q_*$ were radial about a point, its stabilizer would contain a
conjugate of $SO(5)$. This would give
\[
 \dim(GQ_*)\le 21-\dim SO(5)=11<15=\dim(GQ_*).
\]
By \eqref{master:radiation} and free dispersion,
$\vec u(T_{2k})\rightharpoonup(Q_*,0)$ in $\cH$. Since the data radial about a given point form a weakly closed
subspace, $u$ is nonradial. The core trajectory is continuous and
precompact, so Proposition~\ref{prop:abstract-obstruction} completes the proof.
\end{proof}

\section{Stationary geometry}\label{part:background}

\subsection{The full conformal kernel and chart}
We use the profile and the spectral bound of Appendix~\ref{app:single-node}
to identify the conformal kernel. We write $y=(y',y'')\in S^5\subset\R^3\times\R^3$ and take the
regular profile $w$ satisfying
\begin{equation}\label{np:profile-r}
 \begin{gathered}
 w''+2\cot r\,w'+\frac{15}{16}(|w|^{\frac{4}{3}}w-w)=0,
 \qquad w'(0)=w'(\pi)=0,\\
 w(\pi-r)=-w(r),\qquad w>0\ \text{on }[0,\frac{\pi}{2}),
 \qquad w'<0\ \text{on }(0,\pi).
 \end{gathered}
\end{equation}
For the profile of Lemma~\ref{np:certificate}, we set
\begin{equation}\label{np:spherical-state}
 \begin{gathered}
 v(y)=w\bigl(\arccos(|y'|^2-|y''|^2)\bigr),\qquad
 \Omega(x)=\frac{2}{1+|x|^2},\\
 Q(x)=\left(\frac{15}{4}\right)^{\frac{3}{4}}
       \Omega(x)^{\frac{3}{2}}v(\sigma(x)),
 \end{gathered}
\end{equation}
where $\sigma:\R^5\to S^5$ is inverse stereographic projection.
The regular endpoint expansions of $w$ remove the apparent singularities
in this composition. The spherical equation and its Euclidean normalization are
\begin{equation*}
 -\Delta_{S^5}v+\frac{15}{4}v=\frac{15}{4}|v|^{\frac{4}{3}}v,
 \qquad -\Delta Q=|Q|^{\frac{4}{3}}Q.
\end{equation*}

\begin{lemma}[Nondegeneracy of the five-dimensional single-node state]
\label{np:nondegenerate}
The full $\dot H^1(\R^5)$ kernel of $-\Delta-f'(Q)$ is
fifteen-dimensional and equals the tangent space to the conformal orbit
of $Q$. Its elements belong to $L^2(\R^5)$.
\end{lemma}
\begin{proof}
We use the coordinates
\begin{equation*}
 y=\left(\cos\frac{r}{2}\,\theta,\sin\frac{r}{2}\,\omega\right),
 \qquad \theta,\omega\in S^2.
\end{equation*}
We expand a spherical kernel element in product spherical harmonics of
degrees $(j,\ell)$. The radial coefficients satisfy
\begin{equation}\label{np:angular-equation}
 h''+2\cot r\,h'+\left(\mathfrak V(w)
 -\frac{j(j+1)}{4\cos^2(\frac{r}{2})}
 -\frac{\ell(\ell+1)}{4\sin^2(\frac{r}{2})}\right)h=0,
 \qquad \mathfrak V(s)=\frac{35}{16}|s|^{\frac{4}{3}}-\frac{15}{16}.
\end{equation}
By completeness and elliptic regularity, we determine the kernel from
\eqref{np:angular-equation} with the finite-energy branches
\begin{equation*}
 h(r)=O(r^\ell)\quad(r\to0),\qquad
 h(r)=O((\pi-r)^j)\quad(r\to\pi),
\end{equation*}
since the other branches, $r^{-\ell-1}$ and $(\pi-r)^{-j-1}$,
have infinite gradient energy.

The radial sector has trivial kernel by Lemma~\ref{np:radial-kernel}.
By conformal covariance, the conformal gradients and mixed rotations
generate kernel elements in the sectors $(1,0)$, $(0,1)$ and $(1,1)$.
We compute
\begin{equation}\label{np:conformal-modes}
 \begin{aligned}
 \frac{3}{2}y_i'v-\nabla_{S^5}y_i'\cdot\nabla_{S^5}v
 &=\left(\frac{3}{2}w\cos\frac{r}{2}+2w'\sin\frac{r}{2}\right)\theta_i,\\
 \frac{3}{2}y_a''v-\nabla_{S^5}y_a''\cdot\nabla_{S^5}v
 &=\left(\frac{3}{2}w\sin\frac{r}{2}-2w'\cos\frac{r}{2}\right)\omega_a,\\
 \frac{1}{2}(y_i'\partial_{y_a''}-y_a''\partial_{y_i'})v
 &=w'(r)\theta_i\omega_a,\qquad 1\le i,a\le3.
 \end{aligned}
\end{equation}

The regular solution space in each sector is at most one-dimensional.
For the mixed sectors, we use $w'<0$ and the ground-state identity
\begin{equation*}
 \int_0^\pi\sin^2r
 \left(h'^2-\mathfrak V(w)h^2+\frac{2h^2}{\sin^2r}\right)dr
 =\int_0^\pi\sin^2r\,(w')^2
 \left|\left(\frac{h}{w'}\right)'\right|^2dr.
\end{equation*}
The regular endpoint conditions make the boundary terms vanish.
Increasing either positive degree adds a strictly positive term, so
no mixed sector other than $(1,1)$ contributes.
For a one-sided degree $\ell\ge2$,
we would have, by Lemmas~\ref{np:principal-level} and~\ref{np:certificate},
\begin{equation*}
 \ell(\ell+1)=\Lambda_*,\qquad
 4\cdot5<\Lambda_*<5\cdot6,
\end{equation*}
which is impossible. We treat the other one-sided family by reflection.
Thus the three sectors in \eqref{np:conformal-modes} exhaust the kernel.
Their multiplicities give
\[
 \dim\ker L_Q=3+3+3\cdot3=15.
\]

Finally, the conformal identity
\begin{equation*}
 \int_{\R^5}\!\left|\nabla\bigl(\Omega^{\frac{3}{2}}
                  \psi\circ\sigma\bigr)\right|^2dx
 =\int_{S^5}\!\left(|\nabla\psi|^2+\frac{15}{4}|\psi|^2\right)dS
\end{equation*}
extends to an isomorphism $H^1(S^5)\to\dot H^1(\R^5)$.
The projection pole has zero $H^1$ capacity, and the potential term
is conformally covariant, so the weak kernels correspond. By elliptic
regularity, the spherical kernel elements are bounded, and we have
\begin{equation*}
 K_i=O(\rho^{-3})\in L^2(\R^5).
\end{equation*}
\end{proof}

By elliptic bootstrap, $v\in C^{4,\alpha}(S^5)$ for
$\alpha<\frac{1}{3}$, since $f\in C^{2,\frac{1}{3}}$.
At any conformal representative $Q_*$, we choose an $L^2$-orthonormal
kernel basis $K$. As in \cite[Section~3.3]{DKM16}, we apply the inverse
function theorem on a conformal slice, with coordinate $K^*(q-Q_*)$,
to obtain a $C^3$ chart
\begin{equation}\label{i3:eq:chart}
q_\xi=Q_*+K\xi+R(\xi),\qquad
K^*R=0,\quad R(0)=DR(0)=0.
\end{equation}
Differentiating the conformal maps and their weights gives
\begin{equation}\label{proof:conformal-decay}
 |\partial_\xi^\alpha q_\xi|\lesssim\rho^{-3},\qquad
 |\nabla\partial_\xi^\alpha q_\xi|\lesssim\rho^{-4},
 \qquad |\alpha|\le3.
\end{equation}
By the chain rule, we also have
\begin{equation}\label{proof:chart-densities}
 D^2\partial_\xi^\alpha q_\xi
 =\partial_\xi^\alpha f(q_\xi)=O(\rho^{-7})\in Y_4,
 \qquad |\alpha|\le2.
\end{equation}
Dominated convergence gives $D^2q_\xi\in C^2_\xi Y_4$.
Since $f''$ is $\frac{1}{3}$-H\"older continuous, the same identity gives
\begin{equation}\label{proof:chart-holder}
 \|D^2R_{ij}(\xi)-D^2R_{ij}(\zeta)\|_{Y_4}
 \lesssim|\xi-\zeta|^{\frac{1}{3}}.
\end{equation}
The same decay estimates give
\begin{equation*}
 V=f'(Q_*)=O(\rho^{-4}),\qquad \nabla V=O(\rho^{-5}),\qquad
 f''(Q_*)K=O(\rho^{-4}),
\end{equation*}
with continuity across the nodal set.
Differentiating $D^2q_\xi=f(q_\xi)$ twice at $\xi=0$, we obtain
\begin{equation}\label{short:chart-identities}
LR_{ij}(0)=w_iK_j,\quad w_i=f''(Q_*)K_i,\quad
 \int w_iK_jK_k=0,\quad D^2R_{ij}(0)\in Y_4.
\end{equation}
We obtain the cancellation by pairing the first identity with $K_k$.
The equation and decay place all functions in $D(L)$.

\subsection{Balancing the impulse cone}

For a line $b+\mathbb R\omega$, with $\omega\in S^4$ and $b\perp\omega$,
we define the line transform by
\begin{equation*}
 Xv(b,\omega)=\int v(b+s\omega)ds.
\end{equation*}
For a continuous function $v$ with $|v|\lesssim\rho^{-4}$, the line transform is continuous and
\begin{equation*}
 |Xv(b,\omega)|\lesssim\langle b\rangle^{-3}.
\end{equation*}
For $R\in SO(6)$, we define $Q_R$ by replacing $v$ with $v\circ R$
in \eqref{np:spherical-state}.
With the projection pole at $y_6=1$, we have
\begin{equation*}
 (G_0f'(Q_R))(0)=\frac{35}{32\pi^2}
 \int_{S^5}\frac{|v(Ry)|^{\frac{4}{3}}}{(1-y_6^2)^{\frac{3}{2}}}\,dS(y).
\end{equation*}
The weight is integrable, so this is a $C^1$ function of $R$.
We choose $Q_*=Q_R$ at a maximum. The profile and the weight are even
under $y\mapsto-y$, since $v(-y)=v(y)$. The conformal gradients of the
profile are odd, so the corresponding first variations vanish. The
rotational variations vanish at the maximum as well. By Lemma~\ref{np:nondegenerate}, we obtain
\begin{equation}\label{short:cone-balance}
 (G_0w_i)(0)=0,\qquad 1\le i\le15.
\end{equation}
Conformal covariance preserves
nondegeneracy and all chart properties at this representative.

For $|g|\lesssim\rho^{-4}$, we obtain, by polar coordinates and
$D^{-1}D^{-1}=G_0$,
\begin{equation*}
 \begin{aligned}
 D^{-1}g(x)&=\frac{1}{4\pi^3}\int_{S^4}
 Xg(x-(x\cdot\omega)\omega,\omega)\,d\omega,\\
 (G_0g)(0)&=\frac{1}{2\pi^3}\int_{\R^5}
 \frac{D^{-1}g(x)}{|x|^4}\,dx.
 \end{aligned}
\end{equation*}
Both integrals converge absolutely, since $|D^{-1}g|\lesssim\rho^{-3}$.
For $w_c=f''(Q_*)Kc$ with $Xw_c\ge0$, we have, by
\eqref{short:cone-balance},
\begin{equation*}
 D^{-1}w_c\ge0,\qquad
 \int_{\R^5}\frac{D^{-1}w_c(x)}{|x|^4}\,dx=0.
\end{equation*}
It follows that $D^{-1}w_c=0$, and hence $w_c=0$. On the dense set
$\{Q_*\ne0\}$, we then have $Kc=0$. Continuity and linear independence
of the columns of $K$ give $c=0$.
Applying this argument to $\pm c$, we conclude that $Xw_c$ takes both
signs whenever $c\ne0$.

By compactness of the dual unit sphere, we choose finitely many ray columns
\begin{equation*}
 j_\nu=\frac{1}{4}(Xw_i(b_\nu,\omega_\nu))_{i=1}^{15}
\end{equation*}
satisfying
\begin{equation}\label{i3:cone:margin}
 \delta\overline B_{\mathbb R^{15}}
 \subset\operatorname{conv}\{j_1,\ldots,j_M\},\qquad
 \delta=\min_{|c|=1}\max_\nu c\cdot j_\nu>0.
\end{equation}
Consequently,
\begin{equation*}
 \max_\nu|\widetilde j_\nu-j_\nu|<\frac{\delta}{2}
 \quad\Longrightarrow\quad
 (\frac{\delta}{2})\overline B_{\R^{15}}
 \subset\operatorname{conv}\{\widetilde j_1,\ldots,\widetilde j_M\}.
\end{equation*}

\subsection{A direction transverse to similarities}\label{i5:i5:geometry}

By \eqref{np:spherical-state} and our choice of $Q_*$, a conjugate of
$SO(3)\times SO(3)$ fixes $Q_*$. Its orbits on $S^5$ have dimension at
most four. Since $\Sim$ is the stabilizer of $\infty$ in
$\operatorname{Conf}(S^5)$, we have
\begin{equation}\label{proof:similarity-dimension}
 \dim(\Sim\cdot Q_*)\le16-(6-4)=14<15=\dim\ker L.
\end{equation}
If $\mathcal V=T_{Q_*}(\Sim\cdot Q_*)$, then
\eqref{proof:similarity-dimension} allows us to choose $e$ with
\[
 e\in\R^{15},\qquad |e|=1,\qquad Ke\perp_{\dot H^1}\mathcal V.
\]

By compact smooth approximation, escaping similarities send $Q_*$
weakly to zero. We therefore apply the inverse function theorem modulo
isotropy and write the local orbit as
\begin{equation*}
 Q_*+v+O_{\dot H^1}(\|v\|_{\dot H^1}^2),\qquad v\in\mathcal V.
\end{equation*}
For orbit points within $O(d)$ of $q_{de}$, we have
$\|v\|_{\dot H^1}=O(d)$. Projecting onto $\mathcal V^\perp$, we obtain
\begin{equation}\label{proof:profile-separation}
 \dist_{\dot H^1}(q_{de},\Sim\cdot Q_*)
 \ge d\|Ke\|_{\dot H^1}-Cd^2\gtrsim d.
\end{equation}

This separation estimate and \eqref{i3:cone:margin} complete the proof of
Proposition~\ref{presentation:geometry}.

\section{Linear estimates for the projected equation}\label{part:linear}

We control the zero modes by compact tests and the negative modes
by mixed initial and terminal conditions. We then derive the decay
estimates for earlier forcing from the same Green operator.

\subsection{Setting}\label{presentation:linear-setting}

We consider $L=D^2-V$ on $L^2$ with domain $H^2$, where $V$ is real, bounded,
and $|V|\lesssim\rho^{-4}$. We assume its full $\dot H^1$ kernel has an
$L^2$-orthonormal basis
\begin{equation*}
 K=(K_1,\ldots,K_m),
\end{equation*}
with
\begin{equation*}
 K=O(\rho^{-3}),\qquad \nabla K=O(\rho^{-4}),\qquad
 \nu=D^2K=VK\in Y_4.
\end{equation*}
Let $\Psi_1,\ldots,\Psi_{n_-}$ be an orthonormal basis of the negative
spectral subspace. Its elements are exponentially localized, and we write
\begin{equation*}
 L\Psi_l=-\kappa_l^2\Psi_l,\qquad \kappa_l>0,\qquad
 \Pi_- h=\sum_{l=1}^{n_-}\langle\Psi_l,h\rangle\Psi_l .
\end{equation*}
We choose compact smooth tests satisfying
\begin{equation*}
 \Phi K=I,\qquad \Phi\Psi_l=0,
\end{equation*}
and put
\begin{equation*}
P=I-K\Phi,\qquad P_0=I-KK^*.
\end{equation*}
The map $P$ is bounded on energy positions. We use $P_0$ only on $L^2$ variables.
We call energy pairs annihilated by $\Phi$ and all $\Psi_l^*$ complementary.

We extend the sources by zero after $T$ and consider
\begin{equation}\label{short:green-system}
x_t=y-G_0\beta,\qquad y_t=-PLx+P(F_2+F_1),\qquad\Phi x=\Phi y=0,
\end{equation}
where
\begin{equation*}
 \beta\in L^2Y_4,\qquad F_2\in L^2Y_1,\qquad F_1\in L^1L^2,
 \qquad \Phi G_0\beta=K^*F_2=0.
\end{equation*}
We prescribe $(I-\Pi_-)(x(0),y(0))=h_c$, where $h_c\in\mathcal H$
is complementary and $\Pi_-$ acts componentwise. For the negative
coordinates, we impose
\begin{equation*}
a_l=\frac{\Psi_l^*y+\kappa_l\Psi_l^*x}{2\kappa_l},\quad
b_l=\frac{\kappa_l\Psi_l^*x-\Psi_l^*y}{2\kappa_l},\qquad b_l(0)=a_l(T)=0.
\end{equation*}
For $S_J(x)=\|x\|_{S(J)}$, we use
\begin{equation*}
 \|(x,y)\|_{\mathcal R(J)}=\|(\nabla x,y)\|_{L^\infty L^2(J)}
 +S_J(x)+\LE_J(x)+\LG_J(x)+\OO_J(y),
\end{equation*}
where
\begin{equation*}
 \LE=\|x\|_{L^2X_2},\qquad \LG=\|\nabla x\|_{L^2X_1},\qquad
 \OO=\|G_0y\|_{L^2X_4}.
\end{equation*}
\begin{theorem}[Uniform Green estimate]\label{i1:thm:main}
There is a unique energy solution with these data and conditions, and
\begin{equation}\label{short:green-bound}
 \|(x,y)\|_{\mathcal R([0,\infty))}\lesssim
 \|h_c\|_{\mathcal H}+\|\beta\|_{L^2Y_4}+\|F_2\|_{L^2Y_1}+\|F_1\|_{L^1L^2}
 =:\mathfrak D.
\end{equation}
The Green map is linear and its constant is independent of $T$.
\end{theorem}

\subsection{Free estimates and kernels}\label{linear:free}

For the free system
\begin{equation*}
 u_t=y-G_0\beta,\qquad y_t=\Delta u+F_2+F_1,
\end{equation*}
we prove
\begin{equation}\label{short:free-estimate}
 \|(u,y)\|_{\mathcal R}+\|y\|_{L^2X_1}+\|u\|_{L^2X_{\frac{5}{3}}}
 \lesssim\|(u_0,y_0)\|_{\mathcal H}
       +\|\beta\|_{L^2Y_2}+\|F_2\|_{L^2Y_1}+\|F_1\|_{L^1L^2}.
\end{equation}
The same right-hand side bounds $\|u\|_{L^qL^r}$ for
\begin{equation}\label{proof:energy-strichartz-pairs}
 \frac{1}{q}+\frac{5}{r}=\frac{3}{2},\qquad
 4\le r\le\frac{14}{3}.
\end{equation}
We set
\begin{equation*}
 C_0(t)=\cos(tD),\qquad S_0(t)=\frac{\sin(tD)}{D}.
\end{equation*}
For $v_{tt}-\Delta v=F_2+F_1$, we use the local energy estimates
\cite[Lemmas~3.1--3.2]{HWY12}, with
$\delta=\frac{1}{4}$ and $\delta'=\frac{1}{12}$, to obtain
\begin{equation}\label{proof:free-local-energy}
 \|\nabla_{t,x}v\|_{L^\infty L^2\cap L^2X_{\frac{2}{3}}}
 +\|v\|_{L^2X_{\frac{5}{3}}}
 \lesssim\|\vec v(0)\|_{\mathcal H}
       +\|F_2\|_{L^2Y_1}+\|F_1\|_{L^1L^2}.
\end{equation}
For the response to $\beta$ alone, we set
\begin{equation*}
 z(t)=\int_0^tS_0(t-s)G_0\beta(s)\,ds,\qquad
 (u,y)=(-z_t,-\Delta z).
\end{equation*}
We apply \eqref{proof:free-local-energy} to $\partial_jz$ and use
the weighted Hardy and Stein--Weiss inequalities \cite{SW58},
\begin{equation*}
 \|u\|_{X_{\frac{5}{3}}}\lesssim\|\nabla u\|_{X_{\frac{2}{3}}},\qquad
 \|\nabla G_0\beta\|_{Y_1}\lesssim\|\beta\|_{Y_2}.
\end{equation*}
This proves the energy, position, gradient and velocity bounds.
For the homogeneous Newton observation, we have
\begin{equation*}
 \widehat{G_0e^{-itD}g}(\lambda,x)=c\lambda^2\int_{S^4}
e^{i\lambda x\cdot\omega}\widehat g(\lambda\omega)d\omega
\end{equation*}
for $\lambda>0$. By Plancherel, Cauchy--Schwarz on $S^4$, and $\rho^{-4}\in L^2$,
we estimate
\begin{align*}
 \|G_0e^{-itD}g\|_{L^2X_4}^2
 &\lesssim\int_0^\infty\lambda^4
       \int_{S^4}|\widehat g(\lambda\omega)|^2d\omega\,d\lambda\\
 &\lesssim\|g\|_2^2.
\end{align*}

For the spacetime norms, we combine the dual homogeneous estimates
with the Strichartz estimate \cite[Corollary~1.3]{KT} (see also
\cite{Tao}). Since $q\ge\frac{7}{3}>2$, the Christ--Kiselev lemma
\cite{CK01} gives the retarded bound. We estimate $F_1$ by
Minkowski's inequality.
For $t>0$, the Newton kernels in dimension five are
\begin{equation}\label{short:newton-kernels}
 \begin{split}
 [S_0(t)G_0](r)&=\frac{t}{8\pi^2r^3}\mathbf1_{r>t},\\
 [C_0(t)G_0](r)&=\frac{1}{8\pi^2r^3}\mathbf1_{r>t}
                   -\frac{1}{8\pi^2r^2}\delta_{t=r},\\
 [C_0(t)G_0^2](r)&=\frac{r^{-1}-t^2r^{-3}}{16\pi^2}\mathbf1_{r>t}.
 \end{split}
\end{equation}
We obtain the first kernel by the Laplace transform. For the second
and third kernels, we use
\[
 C_0(t)G_0=\partial_t[S_0(t)G_0],\qquad
 C_0(t)G_0^2=G_0^2-\int_0^tS_0(s)G_0\,ds.
\]
For the retarded Newton observation, we include the cone measure and use
\begin{equation*}
 \int_0^\infty[S_0(t)G_0](r)\,dt=\frac{1}{16\pi^2r},\qquad
 \|[C_0(\cdot)G_0](r)\|_{\mathcal M_t}\lesssim r^{-2}.
\end{equation*}
By Young's inequality in time, we reduce the source terms to
\begin{equation*}
 I_3:Y_1\to X_2,\qquad I_4:Y_2\to X_2\subset X_4.
\end{equation*}
By density and truncation, we obtain \eqref{short:free-estimate}, also
with Laplace damping. By the weighted wave estimate
\cite[Theorem~1.5]{FW11}, with $a=1$ and $b=3$, we also have
\begin{equation*}
 \|S_0(t)g\|_{L^2X_{\frac{3}{2}}}
 \le\||x|^{-\frac{3}{2}}S_0(t)g\|_{L^2_{t,x}}
 \lesssim\|g\|_2.
\end{equation*}
On every time interval, we obtain by Sobolev embedding and interpolation
\begin{equation}\label{proof:critical-L4-interpolation}
 \|z\|_{L^4L^4}\lesssim
 S(z)^{\frac{7}{12}}\|\nabla z\|_{L^\infty L^2}^{\frac{5}{12}}.
\end{equation}
By duality with the homogeneous Newton, weighted $S_0$, and
Strichartz estimates, we also have
\begin{equation}\label{proof:free-position-duality}
 \begin{aligned}
 &\left\|\int_0^t\bigl[-C_0(t-s)G_0\beta(s)
                  +S_0(t-s)(F_2+F_1)(s)\bigr]ds\right\|_{L_t^\infty L_x^2}\\
 &\qquad\lesssim\|\beta\|_{L^2Y_4}+\|F_2\|_{L^2Y_{\frac{3}{2}}}
                  +\|F_1\|_{L^{\frac{7}{4}}L^{\frac{14}{11}}}.
 \end{aligned}
\end{equation}

We fix $\psi\in C_c^\infty$ with $m(\psi)=1$, where $m(g)=\int g$,
and write $\kappa_0(t)=S_0(t)G_0\psi$. For $g\in Y_4$, we have
\begin{equation*}
 |\widehat{g-m(g)\psi}(\zeta)|
 \lesssim|\zeta|\|g\|_{Y_4},\qquad |\zeta|\le1.
\end{equation*}
By Plancherel, we obtain
\begin{equation}\label{proof:mean-zero-newton}
 \sup_t\|S_0(t)G_0(g-m(g)\psi)\|_2
 \le\|D^{-3}(g-m(g)\psi)\|_2\lesssim\|g\|_{Y_4}.
\end{equation}

For $\operatorname{Im}z>0$, let $G_z=(D^2-z^2)^{-1}$ and
$H_z=(G_z-G_0)/z^2$. Their kernels extend to $\operatorname{Im}z\ge0$ as
\begin{equation}\label{short:resolvent-kernels}
G_z(r)=\frac{e^{izr}(1-izr)}{8\pi^2r^3},\qquad
H_z(r)=\frac{h(zr)}{8\pi^2r},\quad
h(w)=\frac{e^{iw}(1-iw)-1}{w^2},\quad h(0)=\frac{1}{2}.
\end{equation}
We use
\begin{equation*}
 h(w)=\int_0^1t e^{itw}\,dt,\qquad h(w)+wh'(w)=e^{iw}-h(w).
\end{equation*}
Thus $h,h',h''$ and $h+wh'$ are bounded in the closed upper half-plane.
In particular,
\begin{equation}\label{proof:kernel-derivatives}
 \begin{gathered}
 \partial_zG_z(r)=\frac{ze^{izr}}{8\pi^2r},\\
 |H_z|\lesssim r^{-1},\qquad |H'_z|\lesssim1,\qquad
 |H_z''|\lesssim r,\qquad
 |zH_z|\lesssim r^{-2},\qquad |(zH_z)'|\lesssim r^{-1}.
 \end{gathered}
\end{equation}
We use Stein--Weiss \cite{SW58} in the forms
\begin{equation*}
I_2:Y_1\to X_1,\quad I_3:Y_1\to X_2,\quad
I_3:Y_2\to X_1,\quad I_4:Y_2\to X_2.
\end{equation*}
The two weights are below $\frac{5}{2}$ and sum to the order of integration.
Since $|V|\lesssim\rho^{-4}$, we also have the bounded map
\begin{equation*}
 V:X_2\longrightarrow Y_2.
\end{equation*}
At high frequency, we use the weighted free resolvent bound
of Feng, Soffer, Wu and Yao \cite{FSWY} in the form
\begin{equation}\label{short:high-resolvent}
 \|G_z\|_{Y_2\to X_2}+\|\partial_zG_z\|_{Y_2\to X_2}
 \lesssim |z|^{-1},\qquad |z|\ge1,\quad\Im z\ge0.
\end{equation}
To verify both bounds, take Schwartz functions $f_1,f_2$ and use the
unitary spatial Fourier transform. For $r>0$ set
$a_j(r,\omega)=r^2\widehat f_j(r\omega)$ and extend $a_j$ by zero to $r<0$.
Hardy's and Rellich's inequalities in the Fourier variable, followed by
the one-dimensional Sobolev product estimate, give
\[
 \begin{gathered}
 \|a_j\|_{H^2(\R;L^2(S^4))}\lesssim\|f_j\|_{Y_2},\\
 B(r)=\int_{S^4}a_1(r,\omega)\overline{a_2(r,\omega)}\,d\omega,\qquad
 \|B\|_{H^2(\R)}\lesssim\|f_1\|_{Y_2}\|f_2\|_{Y_2}.
 \end{gathered}
\]
For $C_B(\zeta)=\int_{\R}B(r)/(r-\zeta)\,dr$ and $k=0,1$, the Fourier
representation of the Cauchy kernel gives
$|C_B^{(k)}(\zeta)|\lesssim\|\widehat{B^{(k)}}\|_1\lesssim\|B\|_{H^2}$
in both open half-planes, with the same bounds for their boundary values.
Polar coordinates and partial fractions give
\[
 2z\int_{\R^5}(G_z f_1)\overline{f_2}\,dx=C_B(z)-C_B(-z).
\]
Differentiating this identity and using $X_2=Y_2^*$ proves
\eqref{short:high-resolvent} by density.
By the positive-energy exclusion in \cite[Lemma~1.3(1)]{Miz}, based on
\cite{IJ03,IS06}, we have, for $-V\in L^{\frac{5}{2}}$,
\begin{equation*}
 \lambda>0,\quad \zeta=G_\lambda V\zeta\in L^{\frac{10}{3},2}
 \quad\Longrightarrow\quad\zeta=0.
\end{equation*}
We treat negative frequencies by conjugation,
\begin{equation}\label{proof:resolvent-conjugation}
 G_{-\overline z}=\overline{G_z},\qquad \Im z\ge0.
\end{equation}

\subsection{Construction and observation identities}

The negative coordinates satisfy
\begin{equation}\label{short:negative-green}
 \begin{aligned}
 a_l(t)&=-\frac{1}{2\kappa_l}\int_t^Te^{-\kappa_l(s-t)}
   \langle\Psi_l,(F_2+F_1-\kappa_lG_0\beta)(s)\rangle\,ds,\\
 b_l(t)&=-\frac{1}{2\kappa_l}\int_0^te^{-\kappa_l(t-s)}
   \langle\Psi_l,(F_2+F_1+\kappa_lG_0\beta)(s)\rangle\,ds.
 \end{aligned}
\end{equation}
We have $|a_l(0)|\lesssim\mathfrak D$ and choose
\begin{equation*}
 (x_0,y_0)=h_c+
 \left(\sum a_l(0)\Psi_l,\sum\kappa_la_l(0)\Psi_l\right).
\end{equation*}
Since $-PL=\Delta+V+K\Phi L$ gives a bounded perturbation of the
free energy generator, Duhamel's formula gives a unique solution
on each finite interval.
Since $V:\dot H^1\to\dot H^{-1}$ is compact and $\Phi K=I$,
we apply the argument of \cite[proof of Proposition~3.6]{DKM16}
on the negative complement to obtain
\begin{equation*}
 \Phi x=0,\quad\Psi_l^*x=0\ \text{for all }l
 \quad\Longrightarrow\quad
 \langle Lx,x\rangle\gtrsim\|\nabla x\|_2^2.
\end{equation*}

After $T$, we subtract the decaying negative modes and use the conserved energy
\begin{equation*}
 \frac{1}{2}\langle Lx,x\rangle+\frac{1}{2}\|P_0y\|_2^2
 \gtrsim\|\nabla x\|_2^2+\|y\|_2^2,
\end{equation*}
where $y=PP_0y$. By the Hardy and Hardy--Littlewood--Sobolev inequalities,
we have
\begin{equation*}
 \|x\|_{X_2}\lesssim\|\nabla x\|_2,\qquad
 \|G_0y\|_{X_4}\le\|\rho^{-4}\|_{\frac{5}{2}}\|G_0y\|_{10}
 \lesssim\|y\|_2.
\end{equation*}
Thus the causal extension is tempered and its damped observations belong to $L^2$.

We denote by $(u_f,y_f)$ the causal free response with the same full
initial data and unprojected sources, and set $z_f=G_0y_f$. We extend
these functions by zero to $t<0$ and test the velocity equation against
$K$. In distributions, we have
\begin{equation}\label{short:charge}
 \begin{gathered}
 q=K^*y_0\delta_0+K^*F_1(t)\,dt,\qquad
 \|q\|_{\rm TV}\lesssim\mathfrak D,\\
 \partial_t(\nu^*z_f)=-\nu^*u_f+q.
 \end{gathered}
\end{equation}
Here $\delta_0$ denotes the unit atom at zero. Taking the Fourier
transform of \eqref{short:charge} and testing the position equation
against $\Phi$, we obtain
\begin{equation*}
 \nu^*\widehat u_f=iz\nu^*\widehat z_f+\widehat q,\qquad
 \Phi D^2\widehat z_f=-iz\Phi\widehat u_f.
\end{equation*}
We write
\begin{equation*}
 X=\widehat x,\qquad Z=G_0\widehat y
\end{equation*}
and subtract the free equations to obtain
\begin{align}\label{short:resolvent-system}
X&=\widehat u_f+G_zVX+H_z\nu\gamma,\qquad
Z=\widehat z_f-izH_zVX-\frac{i}{z}\mathscr L_z\gamma,\notag\\
 \mathscr L_z&=(H_z-H_0)\nu,\qquad
 \gamma=\Phi LX-\Phi(\widehat F_2+\widehat F_1).
\end{align}
Here we use the identity
\begin{equation*}
 G_0H_z\nu=\frac{(H_z-H_0)\nu}{z^2}
\end{equation*}
on the difference of two energy solutions. We introduce two auxiliary waves,
\begin{equation*}
 U(t)=\int_{[0,t]}S_0(t-s)K\,dq(s),\qquad
 W(t)=\int_{[0,t]}G_0C_0(t-s)K\,dq(s).
\end{equation*}
Including the initial atom in $dq$, we have
\begin{equation*}
 \widehat U=H_z\nu\widehat q,\qquad
 \widehat W=-\frac{i}{z}\mathscr L_z\widehat q.
\end{equation*}
By \eqref{short:free-estimate} and Minkowski, we obtain
\begin{equation*}
 \|U\|_{L^2X_2}+\|W\|_{L^2X_4}
 \lesssim\|q\|_{\rm TV}\lesssim\mathfrak D.
\end{equation*}
We set
\begin{equation*}
 \mathbb X=(X_2\oplus X_4)^2,\qquad
 \mathcal Z_f=(u_f,z_f,U,W).
\end{equation*}

\subsection{A bounded Lipschitz multiplier}

At zero, $G_0V:X_2\to X_2$ is compact. Its potential tail has norm
$O(R^{-1})$ into $Y_1$. For a compactly supported potential, we use
local Rellich compactness and the Newton $r^{-3}$ tail.
If $x=PG_0Vx$, then $x\in\dot H^1\cap\ker\Phi$.
We apply $L$ and pair with $K$ to obtain
\begin{gather*}
 Lx=-\nu\Phi G_0Vx,\qquad (K^*\nu)\Phi G_0Vx=0,\\
 K^*\nu=(\langle\nabla K_i,\nabla K_j\rangle)>0.
\end{gather*}
Thus $x\in\ker L\cap\ker\Phi=\{0\}$. By the Fredholm alternative,
$I-PG_0V$ is invertible on $X_2$.

We set
\begin{equation*}
 A_z=\nu^*H_z\nu,\qquad T_z=PH_z\nu.
\end{equation*}
Since $\nu\in Y_4$, we have $\rho\nu\in L^1$ and, by
\eqref{proof:kernel-derivatives},
\begin{equation*}
 \|A_z''\|\lesssim\iint|x-y||\nu(x)||\nu(y)|\,dx\,dy<\infty.
\end{equation*}
Since $A_0=\nu^*G_0^2\nu=K^*K=I$, we have
\begin{equation*}
 \frac{A_z-I}{z}=\int_0^1A'_{\theta z}\,d\theta.
\end{equation*}
The following maps are bounded and Lipschitz near zero,
\begin{equation}\label{proof:low-frequency-maps}
 \begin{aligned}
 A_z,\frac{A_z-I}{z}&:\C^m\longrightarrow\C^m,
 &zT_z&:\C^m\longrightarrow X_2,\\
 \mathscr L_z&:\C^m\longrightarrow X_4,
 &z\nu^*H_zV&:X_2\longrightarrow\C^m,\\
 G_zV&:X_2\longrightarrow X_2,
 &zH_zV&:X_2\longrightarrow X_4.
 \end{aligned}
\end{equation}
We control $\mathscr L_z'$ by the bounded kernel $H_z'$,
$\nu\in L^1$, and $1\in X_4$. For the remaining maps, we use
$G_z=G_0+z^2H_z$, \eqref{proof:kernel-derivatives}, and
\begin{equation*}
 X_2\xrightarrow{\ V\ }Y_2\xrightarrow{\ I_4\ }X_2
 \hookrightarrow X_4.
\end{equation*}
We apply $P$ or pair with $\nu\in Y_2$ as needed.
Only $zT_z$ is needed as a Lipschitz map into $X_2$
(Appendix~\ref{app:zero-frequency}).

Pairing the first equation in \eqref{short:resolvent-system} with $\nu$,
we obtain
\begin{equation*}
 A_z\gamma=-\nu^*\widehat u_f-z^2\nu^*H_zVX.
\end{equation*}
We set $\gamma_1=\frac{\gamma+\widehat q}{z}$. Using
\begin{equation*}
 \widehat q=\nu^*\widehat u_f-iz\nu^*\widehat z_f,
\end{equation*}
we rewrite the system as
\begin{align}\label{short:low-matrix}
 (I-PG_zV)X-(zT_z)\gamma_1
 &=P(\widehat u_f-\widehat U),\notag\\
 (z\nu^*H_zV)X+A_z\gamma_1
 &=\frac{A_z-I}{z}\nu^*\widehat u_f-iA_z\nu^*\widehat z_f,\notag\\
 Z&=\widehat z_f-\widehat W-i(zH_zV)X-i\mathscr L_z\gamma_1.
\end{align}
At zero, the coefficient matrix of the first two equations is
$\operatorname{diag}(I-PG_0V,I)$ on $X_2\oplus\C^m$.
The perturbation estimate and the inverse identity bound its inverse
and its Lipschitz norm near zero. By \eqref{proof:low-frequency-maps},
we obtain a bounded Lipschitz family
\[
 \mathbb M_{\rm lo}(z):\mathbb X\longrightarrow X_2\oplus X_4 .
\]
In Fourier variables, this map sends the free observations to $(X,Z)$.
On the full space $\mathbb X$, we treat $U,W$ as independent inputs
without imposing the charge identity.

For $z\ne0$, we substitute $X+\frac{K\gamma}{z^2}$ in
\eqref{short:resolvent-system}. By
\begin{equation*}
 \Phi X=0,\qquad VK=\nu,\qquad
 \Phi\widehat u_f=\frac{i}{z}\Phi D^2\widehat z_f,
\end{equation*}
we obtain
\begin{equation*}
 \begin{aligned}
 X&=P(I-G_zV)^{-1}\widehat u_f,\\
 Z&=(I-H_0\nu\Phi D^2)\widehat z_f\\
  &\quad-iz(H_zV-H_0\nu\Phi G_zV)(I-G_zV)^{-1}\widehat u_f,
 \end{aligned}
\end{equation*}
whenever the inverse exists. Here $H_0\nu:\C^m\to X_4$ and
$\Phi D^2:X_4\to\C^m$ are bounded.

For real $\lambda\ne0$, we prove compactness of $G_\lambda V$ on $X_2$
by the same tail and local argument. If $x=G_\lambda Vx$, we have
\begin{equation*}
 Vx\in Y_2\hookrightarrow L^{\frac{10}{7}}
 \hookrightarrow L^{\frac{10}{7},2},\qquad
 x\in L^{\frac{10}{3},2},
\end{equation*}
by the Lorentz-space uniform Sobolev estimate
\cite[Proposition~1.1 and Corollary~2.2]{Miz}, based on \cite{KY89,KRS87}.
By positive-energy exclusion and \eqref{proof:resolvent-conjugation},
we have $x=0$. Thus $(I-G_\lambda V)^{-1}$ is bounded on each compact
nonzero annulus.
At high frequency, we have $\|G_zV\|<\frac{1}{2}$ by
\eqref{short:high-resolvent} and use the Neumann series.
To control the coefficients and their frequency derivatives in both regions,
we use \eqref{proof:kernel-derivatives}, \eqref{short:high-resolvent}, and
\begin{equation}\label{proof:annular-resolvent-derivative}
 \begin{aligned}
 \partial_z(I-G_zV)^{-1}
 &=(I-G_zV)^{-1}(\partial_zG_z)V(I-G_zV)^{-1},\\
 zH_zV&=\frac{(G_z-G_0)V}{z}.
 \end{aligned}
\end{equation}
On the real axis away from zero, the preceding formula defines a
bounded Lipschitz map $\mathbb M_{\rm reg}$. We choose a strip on which
the inverses persist,
\begin{equation*}
 0\le\operatorname{Im}z\le\epsilon_0,\qquad
 \epsilon_0<\min_l\kappa_l
 \quad\text{if negative modes are present}.
\end{equation*}
The strip avoids the poles $i\kappa_l$. We choose $\chi\in C_c^\infty$
equal to one near zero and supported in the low-frequency region, and set
\begin{equation*}
 \mathbb M_\epsilon(\lambda)
 =\chi(\lambda)\mathbb M_{\rm lo}(\lambda+i\epsilon)
  +(1-\chi(\lambda))\mathbb M_{\rm reg}(\lambda+i\epsilon).
\end{equation*}
The two formulas agree on inputs from the Green system, although their
extensions to $\mathbb X$ may differ. Compactness of the transition
region gives a globally bounded Lipschitz multiplier $\mathbb M_0$,
with
\[
 \sup_{0<\epsilon\le\epsilon_0,\,\lambda\in\R}
       \|\mathbb M_\epsilon(\lambda)\|<\infty,\qquad
 \lim_{\epsilon\downarrow0}
       \|\mathbb M_\epsilon(\lambda)-\mathbb M_0(\lambda)\|=0
 \quad(\lambda\in\R).
\]

For the solution constructed above, we apply damped Plancherel and obtain
\begin{equation*}
 \|e^{-\epsilon t}(x,G_0y)\|_{L^2(X_2\oplus X_4)}\lesssim\mathfrak D.
\end{equation*}
By Fatou's lemma and strong convergence as $\epsilon\to0$, we have
\begin{equation}\label{short:observation-multiplier}
 (x,G_0y)=\mathbb M_0(D_t)\mathcal Z_f,\qquad
 \LE(x)+\OO(y)\lesssim\mathfrak D.
\end{equation}
Finally, we set $v=P_0y$. Since $K^*y=\nu^*G_0y$, the projected pair satisfies
\begin{equation}\label{proof:green-free-system}
x_t=v-G_0(\beta-\nu K^*y),\qquad
v_t=\Delta x+Vx+F_2+P_0F_1.
\end{equation}
By \eqref{short:observation-multiplier}, we have
\begin{equation*}
 \|\beta-\nu K^*y\|_{L^2Y_4}
 +\|Vx+F_2\|_{L^2Y_1}+\|P_0F_1\|_{L^1L^2}
 \lesssim\mathfrak D.
\end{equation*}
We conclude Theorem~\ref{i1:thm:main} by \eqref{short:free-estimate}.
For density, we preserve compatibility by applying $P^*$ to weighted forces,
subtracting $D^2\eta\Phi G_0\beta$ from Newton densities with $\Phi\eta=I$,
and applying $P(I-\Pi_-)$ to complementary data.

By Lipschitz continuity and Young's inequality, we obtain
\begin{equation}\label{short:cutoff-commutator}
 \|[\chi(\frac{t}{R}),\mathbb M_0(D_t)]\mathcal Z\|_{L^2(X_2\oplus X_4)}
 \lesssim R^{-1}\|\mathcal Z\|_{L^2\mathbb X},\qquad\chi\in C_c^\infty,
\end{equation}
since
\begin{equation*}
 \int_\R|\lambda-\mu|R|\widehat\chi(R(\lambda-\mu))|\,d\mu
 =R^{-1}\int_\R|s\widehat\chi(s)|\,ds\lesssim R^{-1}.
\end{equation*}

\subsection{Decay of earlier forcing}\label{part:old}

For the Green system \eqref{short:green-system}, we also assume $\rho^2\nabla V\in L^{\frac{7}{2}}$, which follows from
$|\nabla V|\lesssim\rho^{-5}$ for our stationary potential.
The sources vanish after $\frac{2T}{3}$, $T\ge1$, and the outer horizon is $H\ge\frac{2T}{3}$.
We assume
\begin{equation}\label{short:old-source-classes}
 \begin{gathered}
 \beta\in L^2Y_4,\quad \Phi G_0\beta=0,\qquad
F_2\in L^2Y_{\frac{13}{8}}\cap L^2W^{1,\frac{14}{11}},\quad K^*F_2=0,\\
F_1=F_{\rm p}+G_0\Gamma,\quad
F_{\rm p}\in L^1L^2\cap L^{\frac{7}{4}}W^{1,\frac{14}{11}}\cap L^4W^{1,1},\quad
 \Gamma\in L^1Y_4,\\
 (x_c,y_c)=(G_0B_c,G_0\Gamma_c),\qquad B_c,\Gamma_c\in Y_4,
 \end{gathered}
\end{equation}
with complementary initial data. We denote by $A$ the sum of the
displayed data and source norms. All Sobolev source norms are inhomogeneous.
\begin{theorem}[Decay of earlier forcing]\label{i2:thm:old}
The corresponding mixed Green solution satisfies
\begin{equation}\label{short:old-decay}
 \begin{aligned}
 \LE_{[bT,\infty)}(x)+\OO_{[bT,\infty)}(y)
 &\lesssim_b AT^{-\frac{5}{8}},\qquad b>\frac{2}{3},\\
 \|x\|_{L^2([T,\infty);X_{\frac{5}{3}})}&\lesssim AT^{-\frac{1}{2}},\\
 \|x\|_{L^q([T,\infty);L^r)}&\lesssim AT^{-\frac{1}{r}},
 \qquad (q,r)\text{ as in \eqref{proof:energy-strichartz-pairs}}.
 \end{aligned}
\end{equation}
Its global $\mathcal R$ norm is $\lesssim A$. Constants are independent of
$T,H$. Increasing $H$ beyond the source support leaves the solution unchanged.
\end{theorem}

We take $s=\frac{13}{8}<\frac{5}{3}$ as required in \eqref{short:halo-norms}.

\subsubsection{Free separated estimates}

We first take the sources supported on $[0,T/2]$ and measure the free
response on $J_T=[T,4T]$. The same proof applies to any fixed
separated fractions.
For ordinary forcing $F\in L^2Y_s$, $1<s\le2$, and Newton densities
given by time measures $\mathfrak B,\mathfrak G$ in $Y_4$, the free responses obey
\begin{equation}\label{short:free-tail-table}
\begin{array}{c|cc}
 \text{source}&\|u\|_{L^2(J_T;X_2)}&\|G_0y\|_{L^2(J_T;X_4)}\\ \hline
 (0,F)&\lesssim T^{1-s}\|F\|_{L^2Y_s}&\lesssim T^{1-s}\|F\|_{L^2Y_s}\\
 (-G_0d\mathfrak B,0)&\lesssim T^{-\frac{3}{2}}\|\mathfrak B\|&\lesssim T^{-3}\|\mathfrak B\|\\
 (0,G_0d\mathfrak G)&\lesssim T^{-1}\|\mathfrak G\|&\lesssim T^{-2}\|\mathfrak G\|
\end{array}
\end{equation}
Here the measure norms are total variation in $Y_4$, and initial Newton
data are atoms at zero. To prove \eqref{short:free-tail-table}, we split densities at $|x|=\frac{T}{8}$.
For far inputs, we use the global free estimate and the bounds
\begin{equation*}
 \begin{aligned}
 \|F_{\rm far}\|_{Y_1}&\lesssim T^{1-s}\|F\|_{Y_s},\\
 \|B_{\rm far}\|_{Y_1}&\lesssim T^{-3}\|B\|_{Y_4},\qquad
 \|\Gamma_{\rm far}\|_{Y_2}\lesssim T^{-2}\|\Gamma\|_{Y_4}.
 \end{aligned}
\end{equation*}
For near ordinary inputs, we use exterior support from the Huygens
principle and \eqref{short:newton-kernels}. By Hardy--Rellich duality
and the exterior weights, we obtain
\begin{align*}
 \|u\|_{L^2(J_T;X_2)}
 &\lesssim T^{-\frac{3}{2}}\|u\|_{L^\infty L^2}
 \lesssim T^{-1}\|F\|_{L^2Y_s},\\
 \|G_0y\|_{L^2(J_T;X_4)}
 &\lesssim T^{-\frac{7}{2}}\|G_0y\|_{L^\infty L^2}
 \lesssim T^{-1-s}\|F\|_{L^2Y_s}.
\end{align*}
For near Newton inputs, the translated spatial norms of the first and third
kernels are $O(T^{-\frac{3}{2}})$ in $X_2$ and $O(T^{-\frac{5}{2}})$ in $X_4$.
The first kernel has $X_4$ norm $O(T^{-\frac{7}{2}})$. For $C_0G_0$, we use
exterior support and
\begin{equation*}
 \|G_0B\|_2\lesssim\|B\|_{Y_2}.
\end{equation*}
Minkowski's inequality on $J_T$ contributes a factor $T^{1/2}$.
For a first density $B\in L^2Y_4$ supported on $[0,T/2]$,
Cauchy--Schwarz also gives
\[
 \|B\|_{L^1Y_4}\le (T/2)^{1/2}\|B\|_{L^2Y_4}.
\]
For the $X_{\frac{5}{3}}$ observation, we obtain in the same way
\begin{equation}\label{proof:free-weighted-tail}
 \|u\|_{L^2(J_T;X_{\frac{5}{3}})}
 \lesssim T^{-\frac{5}{8}}\|F\|_{L^2Y_{\frac{13}{8}}}
 +T^{-\frac{7}{6}}\|\mathfrak B\|
 +T^{-\frac{2}{3}}\|\mathfrak G\|.
\end{equation}
Indeed, the near ordinary and first Newton terms have the exterior factor
$T^{\frac{1}{2}-\frac{5}{3}}=T^{-\frac{7}{6}}$, and the translated
kernel $S_0(t)G_0$ has $X_{\frac{5}{3}}$ norm $O(T^{-\frac{7}{6}})$.
We use \eqref{short:free-estimate} for the far terms.

For the unweighted source, we prove
\begin{equation}\label{short:pure-tail}
 \begin{aligned}
 \|u\|_{L^2(J_T;X_{\frac{5}{3}})}
 &\lesssim T^{-\frac{1}{2}}\bigl(\|F\|_{L^1L^2}
                                  +\|F\|_{L^4W^{1,1}}\bigr),\\
 \|u\|_{L^2(J_T;X_2)}
 &\lesssim T^{-\frac{9}{14}}\bigl(\|F\|_{L^1L^2}
                                  +\|F\|_{L^4W^{1,1}}\bigr),\\
 \|G_0u_t\|_{L^2(J_T;X_4)}
 &\lesssim T^{-\frac{5}{4}}\|F\|_{L^4W^{1,1}}.
 \end{aligned}
\end{equation}
For smooth annular projections $P_N$ and fixed $0<a<b$, we use the
spherical Fourier formula and Plancherel in time to obtain
\begin{equation}\label{proof:annular-kernels}
 \sup_x\|S_0(t)P_N(x)\|_{L^2([aT,bT])}\lesssim N^{\frac{3}{2}}T^{-2},\qquad
 \sup_x\|C_0(t)G_0P_N(x)\|_{L^2([aT,bT])}\lesssim N^{\frac{1}{2}}T^{-2}.
\end{equation}
For $|x|\ge aT/2$, the decay in \eqref{proof:annular-kernels} follows from
\[
 |\widehat{d\sigma_{S^4}}(\zeta)|\lesssim(1+|\zeta|)^{-2}.
\]
For $|x|<aT/2$, we integrate by parts in the radial frequency when
$NT>1$, and use Fourier inversion when $NT\le1$.
We put
\begin{equation*}
 A=\|F\|_{L^1L^2}+\|F\|_{L^4W^{1,1}}.
\end{equation*}
For $N\ge1$ and $\sigma\in\{\frac{5}{3},2\}$, we use
\eqref{proof:annular-kernels}, \eqref{proof:free-position-duality}, and
the homogeneous velocity estimate in \eqref{short:free-estimate},
respectively, to obtain
\begin{equation*}
 \begin{aligned}
 \|P_{\le N}u\|_{L_x^\infty L_t^2(J_T)}
 &\lesssim AT^{-\frac{5}{4}}N^{\frac{1}{2}},\\
 \|P_{\le N}u\|_{L^2(J_T\times\R^5)}
 &\lesssim T^{\frac{1}{2}}\|F\|_{L^1L^2}^{\frac{3}{7}}
       \|F\|_{L^4L^1}^{\frac{4}{7}}
 \lesssim AT^{\frac{1}{2}},\\
 \|P_{>N}u\|_{L^2(J_T;X_\sigma)}
 &\lesssim\|D^{-1}P_{>N}F\|_{L^1L^2}\lesssim AN^{-1}.
 \end{aligned}
\end{equation*}
For the first bound, we apply H\"older on the source interval and sum
with enlarged annular cutoffs,
\begin{equation*}
 \sum_{M\le N}M^{\frac{3}{2}}\|P_MF\|_{L^4L^1}
 \lesssim A\sum_{M\le N}M^{\frac{3}{2}}\min(1,M^{-1})
 \lesssim AN^{\frac{1}{2}}.
\end{equation*}
For the third bound, we use $S_0(t)=\sin(tD)D^{-1}$ and Minkowski.
We apply
\begin{equation*}
 \|\rho^{-\sigma}g\|_2
 \lesssim\|g\|_2^{1-\frac{2\sigma}{5}}\|g\|_\infty^{\frac{2\sigma}{5}}
\end{equation*}
to $g(x)=\|P_{\le N}u(\cdot,x)\|_{L^2(J_T)}$ and obtain
\begin{equation*}
 \|u\|_{L^2(J_T;X_\sigma)}
 \lesssim A\left(N^{-1}+T^{-\frac{7\sigma-5}{10}}N^{\frac{\sigma}{5}}\right).
\end{equation*}
Taking $N=T^{\frac{7\sigma-5}{2\sigma+10}}$, up to a dyadic factor,
we obtain the first two bounds in \eqref{short:pure-tail}.
For the Newton observation, we use $\rho^{-4}\in L^2$ and
$\sum_MM^{\frac{1}{2}}\min(1,M^{-1})<\infty$ in
\eqref{proof:annular-kernels}.

For the Strichartz tails, we use the standard localized wave estimate
\cite[Section~8]{KT},
\begin{equation}\label{proof:free-wave-dispersion}
 \|e^{itD}P_N\|_{L^1\to L^\infty}\lesssim N^5(1+N|t|)^{-2}.
\end{equation}
By interpolation and Bernstein, we obtain, for $t\ge1$,
\begin{align}\label{short:separated-dispersion}
 \|S_0(t)F\|_r
 &\lesssim t^{-2+\frac{4}{r}}\|F\|_{W^{1,\frac{14}{11}}},
 \qquad 4\le r\le\frac{14}{3},\notag\\
 \|C_0(t)G_0B\|_r
 &\lesssim t^{-2+\frac{4}{r}}\|B\|_{Y_4},\qquad
 \|S_0(t)G_0\Gamma\|_r
 \lesssim t^{-2+\frac{5}{r}}\|\Gamma\|_{Y_4}.
\end{align}
For the first bound, Bernstein from $L^{\frac{14}{11}}$ to $L^{r'}$
reduces the dyadic sum to
\begin{equation*}
 \sum_{N\le1}N^{\frac{13}{14}-\frac{1}{r}}
 +\sum_{N>1}N^{-\frac{1}{14}-\frac{1}{r}}<\infty.
\end{equation*}
For the second bound, we use $Y_4\subset L^1\cap L^{r'}$ and sum
\begin{equation*}
 \sum_{N\le1}N^{1-\frac{1}{r}}+\sum_{N>1}N^{1-\frac{6}{r}}<\infty.
\end{equation*}
We obtain the last bound directly from \eqref{short:newton-kernels}.

For the pairs in \eqref{proof:energy-strichartz-pairs}, we have
\begin{equation*}
 2-\frac{4}{r}-\frac{1}{q}-\frac{1}{2}=\frac{1}{r},\qquad
 2-\frac{5}{r}-\frac{1}{q}=\frac{1}{2}.
\end{equation*}
By H\"older on the source support, we obtain
\begin{equation*}
 \|u\|_{L^q([T,\infty);L^r)}\lesssim AT^{-\frac{1}{r}}.
\end{equation*}
These estimates hold on every later dyadic scale, with constants depending
only on the separation fractions.

\subsubsection{Transfer to the projected equation}

Once $H$ exceeds the source support, \eqref{short:negative-green}
gives $|a_l(0)|\lesssim A$ uniformly in $H$. Thus the full initial
pair satisfies
\[
 (x_0,y_0)=(G_0B_0,G_0\Gamma_0),\qquad
 \|B_0\|_{Y_4}+\|\Gamma_0\|_{Y_4}\lesssim A.
\]
By \eqref{short:charge} and $K^*F_2=0$, the coefficient measure satisfies
\begin{equation*}
 \supp q\subset[0,\frac{2T}{3}],\qquad \|q\|_{\rm TV}\lesssim A.
\end{equation*}
The auxiliary waves have Newton density $\nu\,dq$. Summing the free
estimates over later dyadic intervals, we obtain
\begin{equation}\label{proof:free-tail-observations}
 \|\mathcal Z_f\|_{L^2\mathbb X}\lesssim A,\qquad
 \|\mathcal Z_f\|_{L^2([cT,\infty);\mathbb X)}\lesssim_c AT^{-\frac{5}{8}},
 \qquad c>\frac{2}{3}.
\end{equation}
For $b>\frac{2}{3}$, we choose $\frac{2}{3}<c<b$ and
$\chi\in C_c^\infty$ with $0\le\chi\le1$, $\chi=1$ on $[0,c]$,
and $\chi=0$ on $[b,\infty)$. Since $\mathcal Z_f=0$ for $t<0$,
we apply \eqref{short:observation-multiplier} and
\eqref{short:cutoff-commutator} to obtain
\begin{align*}
 \LE_{[bT,\infty)}(x)+\OO_{[bT,\infty)}(y)
 &\lesssim\|(1-\chi(\frac{t}{T}))(x,G_0y)\|_{L^2(X_2\oplus X_4)}\\
 &\lesssim_b\|(1-\chi(\frac{t}{T}))\mathcal Z_f\|_{L^2\mathbb X}
        +T^{-1}\|\mathcal Z_f\|_{L^2\mathbb X}\\
 &\lesssim_b A(T^{-\frac{5}{8}}+T^{-1})\lesssim_b AT^{-\frac{5}{8}}.
\end{align*}

For the remaining norms, we use the free form \eqref{proof:green-free-system},
with $v=P_0y$, and split the feedback at $\frac{3T}{4}$.
By the global gradient and Newton observations, the preceding tail bound,
and $\rho^2\nabla V\in L^{\frac{7}{2}}$, we have
\begin{equation*}
 \begin{aligned}
 \|(Vx)_{\rm early}\|_{L^2(Y_{\frac{13}{8}}\cap W^{1,\frac{14}{11}})}
 +\|(\nu K^*y)_{\rm early}\|_{L^2Y_4}&\lesssim A,\\
 \|(Vx)_{\rm late}\|_{L^2Y_1}
 +\|(\nu K^*y)_{\rm late}\|_{L^2Y_4}&\lesssim AT^{-\frac{5}{8}}.
 \end{aligned}
\end{equation*}
We apply \eqref{proof:free-weighted-tail} and \eqref{short:separated-dispersion}
to the early part, and \eqref{short:free-estimate} and the free Strichartz
bound to the late part. Since $\frac{5}{8}>\frac{1}{2}>\frac{1}{r}$,
we obtain the last two bounds in \eqref{short:old-decay} for both contributions.
We treat the external sources by
\eqref{proof:free-weighted-tail}, \eqref{short:pure-tail}, and
\eqref{short:separated-dispersion}, summing over later dyadic intervals.
The projection correction
\begin{equation*}
 P_0F_1-F_1=-G_0\nu K^*F_1
\end{equation*}
is a second Newton source of size $O(A)$ in $L^1Y_4$.
We conclude \eqref{short:old-decay} by approximation with the compatibility
corrections of Section~\ref{part:linear}.

\section{Radiation packets and quadratic controls}\label{part:catalogue}

\subsection{The radiation family}

With $K$ and $L$ as in Section~\ref{part:linear}, we define the force
and the impulse of a wave $v$ by
\begin{equation}\label{i3:eq:force}
w_i=f''(Q_*)K_i,\qquad
F_i(v(t))=\frac{1}{2}\int_{\R^5}w_i(x)v(t,x)^2\dd x,\qquad
J_i(v)=\int_\R F_i(v(t))\dd t.
\end{equation}
Here $|w_i|\lesssim\rho^{-4}$, including at the nodal set.

\begin{theorem}[Finite family of radiation waves]\label{i3:thm:catalogue}
We can choose a conformal representative $Q_*$ of the single-node state
in Lemma~\ref{np:nondegenerate}, with the balancing property of
Proposition~\ref{presentation:geometry}, and finitely many real solutions
\begin{equation*}
v_{\nu,tt}+Lv_\nu=0,\qquad 1\le\nu\le M,
\end{equation*}
with the following properties.
\begin{enumerate}[label=(\roman*)]
\item For pairwise disjoint compact bands $B_\nu\subset\R\setminus\{0\}$,
\begin{equation*}
 (v_\nu,\partial_tv_\nu)\in C_b(\R;\cH),\qquad
 v_\nu\in C_b(\R;L^2),\qquad
 \supp\mathcal F_t(v_\nu,\partial_tv_\nu)\subset B_\nu.
\end{equation*}
Both components are orthogonal to the zero and negative modes, and
the pair scatters in energy in both time directions.
\item Each wave has finite global $\mathcal R$ and $L_t^4L_x^4$ norms.
For every integer $k\ge0$,
\begin{equation*}
 \langle t\rangle\partial_t^kv_\nu\in L_t^2X_2,\qquad
 \langle t\rangle\nabla v_\nu\in L_t^2X_2,\qquad
 \langle t\rangle G_0\partial_tv_\nu\in L_t^2X_4.
\end{equation*}
For $T\ge1$,
\begin{equation*}
 T^{\frac{5}{7}}\|v_\nu\|_{S(|t|>T)}
 +T^{\frac{3}{4}}\|v_\nu\|_{L_t^4L_x^4(|t|>T)}\lesssim1.
\end{equation*}
\item The columns $J(v_\nu)$ satisfy
\[
 0\in\operatorname{int}\operatorname{conv}\{J(v_\nu):1\le\nu\le M\}.
\]
On a neighborhood of $\mathcal C=[-3,3]^{15}$, there are smooth,
strictly positive functions $c_\nu$ such that
\begin{equation}\label{i3:eq:amplitudes}
 \sum_{\nu=1}^M c_\nu(z)^2J(v_\nu)=z,\qquad z\in\mathcal C.
\end{equation}
Their derivatives of every fixed order are bounded on the cube.
\item The projected pair $(Pv_\nu,P\partial_tv_\nu)$ is exactly the homogeneous projected Green solution of Section~\ref{part:linear} with its own complementary initial data, for every outer terminal horizon.
\end{enumerate}
We fix this family and its constants before choosing the small
parameter or the number of intervals.
\end{theorem}

In Section~\ref{i3:sec:gauge}, we relate the force of $Pv_\nu$ to
that of $v_\nu$.

\subsection{Construction and estimates for the radiation waves}\label{presentation:packets}

We first construct a free packet whose impulse approximates a given
line integral. We fix $\omega\in S^4$, $b\perp\omega$, and a real
Schwartz function $a$ with $\|a\|_2=1$ and $\widehat a\in C_c^\infty$.
For $\delta>0$ we set $a_\delta(x)=\delta^{-5/2}a(x/\delta)$ and
\begin{equation*}
a_h=e^{\frac{i\omega\cdot x}{h}}a_\delta(x-b),\qquad
u_h=\operatorname{Re}(he^{-itD}a_h).
\end{equation*}
With $b,\delta$ fixed, we have, uniformly as $h\downarrow0$,
\begin{equation}\label{short:free-packet}
 \|\langle t\rangle u_h\|_{L^2X_2}\lesssim_{b,\delta}h,\qquad
 \|u_h\|_{L^2X_2(|t|>T)}\lesssim_{b,\delta}hT^{-1}.
\end{equation}
By the free velocity observation, we have
\begin{equation*}
 \|e^{-itD}g\|_{L^2X_1}\lesssim\|g\|_2.
\end{equation*}
Writing $D_j=-i\partial_{x_j}$, we obtain by Fourier differentiation
\begin{equation}\label{proof:packet-commutator}
t he^{-itD}a_h=\sum_j\frac{D_j}{D}
 \bigl[x_jhe^{-itD}a_h-he^{-itD}(x_ja_h)\bigr].
\end{equation}
Since $\rho^{-4}\in A_2(\R^5)$, we apply the weighted Riesz-transform
estimate \cite[Theorem~1.3]{Hyt12} to obtain
\begin{equation*}
 \left\|\frac{D_j}{D}g\right\|_{X_2}\lesssim\|g\|_{X_2},
 \qquad |x|\rho^{-2}\le\rho^{-1}.
\end{equation*}
By \eqref{proof:packet-commutator} and the free observation,
we conclude \eqref{short:free-packet}.

For a continuous weight $w$ with $|w|\lesssim\rho^{-4}$, we write
\[
 J_w(u)=\frac12\int_{\R}\int_{\R^5}w(x)u(t,x)^2\,dx\,dt .
\]
We extend $Xw(b,\omega)$ to all $b\in\R^5$ by projecting $b$ onto
$\omega^\perp$. We then obtain
\begin{equation}\label{i3:packet-ray-limit}
h^{-2}J_w(u_h)\longrightarrow
 \frac{1}{4}\int |a_\delta(y)|^2Xw(b+y,\omega)\,dy
 \longrightarrow_{\delta\downarrow0}\frac{1}{4}Xw(b,\omega).
\end{equation}
For the first limit, we expand
\begin{equation*}
 |\frac{\omega}{h}+\zeta|=h^{-1}+\omega\cdot\zeta+O_\delta(h).
\end{equation*}
For each fixed $T$, we then have
\begin{equation*}
 \sup_{|t|\le T}\left\|h^{-1}u_h(t)
 -a_\delta(x-b-t\omega)\cos(\frac{\omega\cdot x-t}{h})\right\|_2
 \lesssim_{T,b,\delta}h.
\end{equation*}

We average $\cos^2$ by Riemann--Lebesgue and let $T\to\infty$,
using the tail bound $O(h^2T^{-2})$ from \eqref{short:free-packet}.
We obtain the second limit from the approximate identity $|a_\delta|^2$.

\subsubsection{Linearized realization}\label{i3:sec:R-realization}

To pass from the free packet to a solution of the linearized equation,
we use the temporal Fourier transform
\begin{equation*}
 \widehat v(\tau)=\int e^{it\tau}v(t)\,dt.
\end{equation*}
We denote by $B_h$ the two fixed-width temporal bands about $\pm h^{-1}$
containing the support of $\widehat u_h$.
On a neighborhood of $B_h$, we use the free high-frequency estimate of
Section~\ref{part:linear} to obtain
\begin{equation*}
 \|G_{\tau+i0}V\|_{X_2\to X_2}
 +\|\partial_\tau G_{\tau+i0}V\|_{X_2\to X_2}\lesssim h.
\end{equation*}
We set
\begin{equation*}
 A=(I-G_{\tau+i0}V)^{-1},\qquad
 \|A\|\le2,\qquad \|A-I\|+\|A'\|\lesssim h,
\end{equation*}
and define $\widehat v_h$ by
\[
 \widehat v_h(\tau)=
 \begin{cases}
 A(\tau)\widehat u_h(\tau),&\tau\in B_h,\\
 0,&\tau\notin B_h .
 \end{cases}
\]
After a smooth cutoff near $B_h$, we obtain by the $H^1_\tau(X_2)$ product rule
\begin{equation}\label{short:packet-realization}
 \|\langle t\rangle v_h\|_{L^2X_2}\lesssim h,\qquad
 \|\langle t\rangle(v_h-u_h)\|_{L^2X_2}\lesssim h^2.
\end{equation}
By conjugate symmetry, $v_h$ is real. By Hardy duality and
Cauchy--Schwarz in time, we also have
\begin{equation*}
 \|\langle t\rangle Vv_h\|_{L^2Y_2}+\|Vv_h\|_{L^1L^2}
 +\|D^{-1}(Vv_h)\|_{L^1L^2}\lesssim h.
\end{equation*}
The weighted construction therefore has the energy representative
\begin{equation}\label{short:packet-duhamel}
v_h(t)=u_h(t)+\int_{-\infty}^tS_0(t-s)Vv_h(s)\,ds.
\end{equation}
Indeed, by free observation, Minkowski, and limiting absorption on $B_h$,
the damped integral converges in $L^2X_2$ and its Fourier transform satisfies
\begin{equation*}
 G_{\tau+i\epsilon}\widehat{Vv_h}
 \longrightarrow G_{\tau+i0}V\widehat v_h
 =\widehat v_h-\widehat u_h.
\end{equation*}
We obtain the equation, energy continuity, free observations, and scattering
from \eqref{short:packet-duhamel}. The energy error is $O(h)$, whereas
\begin{equation*}
 \|(u_h,u_{h,t})\|_{\mathcal H}^2=1+h^2\|\nabla a_\delta\|_2^2.
\end{equation*}

We apply \eqref{short:packet-duhamel} in $L^2\times\dot H^{-1}$ and obtain
\begin{equation}\label{i3:R:minusone}
 \|v_h\|_{C_bL^2}+\|D^{-1}v_{h,t}\|_{C_bL^2}
 \lesssim h+\|D^{-1}(Vv_h)\|_{L^1L^2}\lesssim h.
\end{equation}
By temporal Bernstein on $B_h$ and $D^2v_h=Vv_h-v_{h,tt}$,
we have $v_h\in C_bD(L)$. By the equation and $0\notin B_h$, we have
\begin{equation*}
 \begin{gathered}
 \tau^2K^*\widehat v_h=0,\qquad
 (\tau^2+\kappa_l^2)\Psi_l^*\widehat v_h=0,\\
 K^*v_h=\Psi_l^*v_h=0.
 \end{gathered}
\end{equation*}

\subsubsection{Weighted observations and tails}\label{presentation:packet-observations}

On each fixed band, temporal multipliers give
$\langle t\rangle\partial_t^kv\in L^2X_2$.
Integrating $D^2v=Vv-v_{tt}$ by parts, we obtain
\begin{equation*}
 \|\rho^{-2}\nabla v\|_2^2\lesssim
\|\rho^{-2}v\|_2(\|\rho^{-2}v_{tt}\|_2+\|\rho^{-2}v\|_2),
\end{equation*}
hence $\langle t\rangle\nabla v\in L^2X_2$. For the Newton observation,
we use Rellich duality and
\begin{equation*}
 \widehat{G_0v_t}=-i\tau^{-1}(I-G_0V)\widehat v,\qquad
I-G_0V:X_2\longrightarrow X_4,
\end{equation*}
and apply the same temporal multipliers.

By H\"older and $|\nabla V|\lesssim\rho^{-5}$, we have
\begin{equation*}
 \langle t\rangle Vv\in L^2(Y_1\cap W^{1,\frac{14}{11}}).
\end{equation*}
For $t>T\ge1$, we split \eqref{short:packet-duhamel} at $s=\frac{T}{2}$,
writing $v=v_{\rm old}+v_{\rm new}$ and including $u_h$ in $v_{\rm old}$.
For the pairs in \eqref{proof:energy-strichartz-pairs}, we use
\eqref{short:separated-dispersion}, the fixed-band Schwartz data, and
the free Strichartz bound to obtain
\begin{equation}\label{proof:packet-tails}
 \begin{aligned}
 \|v_{\rm old}\|_{L^q([T,\infty);L^r)}
 &\lesssim T^{-\frac{1}{2}-\frac{1}{r}}
       (1+\|Vv\|_{L^1W^{1,\frac{14}{11}}})
 \lesssim T^{-\frac{1}{2}-\frac{1}{r}},\\
 \|v_{\rm new}\|_{L^qL^r}
 &\lesssim\|\mathbf1_{s>\frac{T}{2}}Vv\|_{L^2Y_1}\lesssim T^{-1}.
 \end{aligned}
\end{equation}
Taking $(q,r)=(7/3,14/3)$ and $(4,4)$, we obtain the two tails in Theorem~\ref{i3:thm:catalogue}.
For $t<-T$, we use the free-packet bound and the free Strichartz estimate
for the source restricted to $s<-T$.

\subsection{Controls and the two-moment map}\label{presentation:packet-controls}
\label{i3:sec:gauge}

By \eqref{short:packet-realization}, we have
\begin{equation*}
 |J_w(v_h)-J_w(u_h)|\lesssim h^3.
\end{equation*}
By \eqref{i3:packet-ray-limit} and the strict margin \eqref{i3:cone:margin},
we choose the envelopes, then take the frequencies successively large.
Thus the bands $B_\nu$ are pairwise disjoint and the columns
$J(v_\nu)$ span the full positive cone. Positive rescaling preserves
this cone. Writing $\mathsf J=(J(v_1)\ \cdots\ J(v_M))$, we choose
\begin{equation*}
 \mathsf J\mathsf S=I,\qquad \beta\in(0,\infty)^M,
 \qquad \mathsf J\beta=0.
\end{equation*}
For $\beta$, we express $-\sum_\nu J(v_\nu)$ as a nonnegative combination
and add one to each coefficient. We rescale $\beta$ so that
\begin{equation*}
 c_\nu(z)=\sqrt{\beta_\nu+(\mathsf Sz)_\nu}
\end{equation*}
is smooth and positive near the control cube. These amplitudes satisfy
\eqref{i3:eq:amplitudes}.
By Plancherel applied to $\rho^{-2}v_\nu$ and the bounded spatial
multiplier $\rho^4w_i$, we have
\begin{equation*}
 \int_\R\int_{\R^5}w_i v_\nu v_\mu\,dx\,dt=0,\qquad \nu\ne\mu.
\end{equation*}
We therefore set
\begin{equation}\label{proof:group-impulse}
 v_z=\sum_{\nu=1}^M c_\nu(z)v_\nu,\qquad J(v_z)=z,\qquad z\in\mathcal C.
\end{equation}
The norm and tail bounds are uniform in $z$.

For any finite train of these waves, we put
\begin{equation*}
 \alpha=\Phi v,\qquad r=Pv.
\end{equation*}
We have $r_{tt}+PLr=0$ and zero negative coordinates.
By Theorem~\ref{i1:thm:main}, $(r,r_t)$ is the homogeneous Green solution
at every horizon. To retain the observations, we use
\begin{equation*}
 G_0r_t=G_0v_t-G_0^2\nu\,\alpha'.
\end{equation*}
For a fixed packet, we obtain by compact testing
\begin{equation*}
 \langle t\rangle(\alpha,\alpha')\in H^1(\R).
\end{equation*}
By Sobolev embedding, we have
\begin{equation*}
 \|\alpha\|_{L^q(|t|>T)}+\|\alpha'\|_{L^q(|t|>T)}\lesssim T^{-1},
 \qquad 2\le q\le\infty.
\end{equation*}
Since $Pv_\nu=v_\nu-K\alpha$, the projection preserves the Strichartz
tails. This proves Theorem~\ref{i3:thm:catalogue}.

For $b_{ij}=\langle R_{ij}(0),v_t\rangle$, we have, by
\eqref{short:chart-identities},
\begin{equation*}
 b_{ij}'=-\int w_iK_jv,\qquad \int w_iK_jK_k=0.
\end{equation*}
Expanding $F_i(v-K\alpha)$, we obtain
\begin{equation}\label{i3:cons:gaugeforce}
F_i(Pv)=F_i(v)+\alpha_jb_{ij}'
       =F_i(v)+(b_{ij}\alpha_j)'-b_{ij}\alpha_j'.
\end{equation}
We sum over repeated indices. Since $R_{ij},v\in D(L)$, the pairings
are well defined. For each packet, $\alpha,b\in H^1(\R)$, so
$\int_\R(b\alpha)'=0$. For the weighted moment, we integrate this
derivative by parts.

\subsubsection{The two-moment map}\label{i3:sec:blocks}

We place two copies of the family \eqref{proof:group-impulse} on each
time interval. For $\varepsilon>0$ and $T_0\ge1$, we set
\begin{equation*}
T_j=4^jT_0,\qquad \eta_j=\varepsilon2^{-j},\qquad
d=\varepsilon^2T_0=T_j\eta_j^2.
\end{equation*}
The intervals and their right endpoints are
\begin{equation*}
 I_0=[0,4T_0],\qquad I_j=[T_j,4T_j]\quad(j\ge1),\qquad e_j=4T_j.
\end{equation*}
For $N\ge1$ and controls
\begin{equation*}
 a_j=(x_j,y_j)\in\mathcal C^2,
\end{equation*}
we define
\begin{equation}\label{i3:eq:train}
v_a(t)=\sum_{j=0}^{N-1}\eta_j
 \bigl[v_{x_j}(t-2T_j)+v_{y_j}(t-3T_j)\bigr].
\end{equation}
The train solves $v_{tt}+Lv=0$, and its projected wave is a
homogeneous Green solution. All estimates below are uniform in $N$ and
in the controls. We use the maximum component norm on each finite control
space and the corresponding $\ell^\infty$ norm over the blocks.

We define the two-moment map and its diagonal approximation by
\begin{align}
 \mathcal M_N(a)_j
 &=\left(\eta_j^{-2}\int_{I_j}F(v_a(t))\dd t,\quad
d^{-1}\int_{I_j}(e_j-t)F(v_a(t))\dd t\right),\label{i3:eq:momentmap}\\
 \mathcal L(x,y)&=(x+y,2x+y),\qquad
 \mathcal L^{-1}(u,v)=(v-u,2u-v).\label{i3:eq:blocklinear}
\end{align}
The map $\mathcal L$ acts on each block.

\begin{theorem}[Uniform reference map]\label{i3:thm:controls}
The moment map satisfies
\begin{equation}\label{i3:eq:uniform-reference}
 \sup_{a\in(\mathcal C^2)^N}
 \|\mathcal M_N(a)-\mathcal L a\|_{\ell^\infty}
 \lesssim T_0^{-1}
\end{equation}
for every $N$, and $\mathcal M_N$ is continuous. The train and its projected wave have global $\mathcal R$ and
$L_t^4L_x^4$ norms $\lesssim\varepsilon$. On each block,
\begin{align}
 &\|v_a\|_{L_t^2X_2(I_j)}+\|v_a\|_{S(I_j)}
 +\|v_a\|_{L_t^4L_x^4(I_j)}
 +\|G_0\partial_tv_a\|_{L_t^2X_4(I_j)}\notag\\
 &\hspace{1cm}
 +\|\Phi v_a\|_{L^2(I_j)}+\|\Phi\partial_tv_a\|_{L^2(I_j)}
 +\|b_a\|_{L^2(I_j)}\lesssim \eta_j,
 \qquad (b_a)_{ik}=\langle R_{ik}(0),\partial_tv_a\rangle.
 \label{i3:eq:blocknorms}
\end{align}
The corresponding local-energy, $S$, $L_t^4L_x^4$, and weighted potential
bounds hold for $Pv_a$ as well.

\end{theorem}

\begin{proof}
By Minkowski and $\sum_{j<N}\eta_j\le2\varepsilon$, we obtain the
global bounds, also after projection. For the block bounds, we use
\begin{equation}\label{proof:off-block-tail}
 T_j^{-\sigma}\sum_{k<j}\eta_k+
 \sum_{k>j}\eta_kT_k^{-\sigma}
 \lesssim T_0^{-\sigma}\eta_j,\qquad \sigma\ge\frac{1}{2}.
\end{equation}
The first-moment observations have tail exponent one, and the
Strichartz tails have exponents $\frac{5}{7}$ and $\frac{3}{4}$.
Together with the two current groups, \eqref{proof:off-block-tail}
proves \eqref{i3:eq:blocknorms}, including the projected bounds.

For the current square terms, we use $J(v_z)=z$ and
\begin{equation*}
 \int_{|t|>T}\|v_z(t)\|_{X_2}^2dt\lesssim T^{-2},\qquad
 \int_{|t|>T}|t|\|v_z(t)\|_{X_2}^2dt\lesssim T^{-1}.
\end{equation*}
Extending the block integrals to $\R$, we obtain the normalized moments
\begin{equation*}
 x_j+y_j+O(T_j^{-2}),\qquad 2x_j+y_j+O(T_j^{-1}).
\end{equation*}
By $|s_1-s_2|\le|t-s_1|+|t-s_2|$ and Cauchy--Schwarz, we have
\begin{equation}\label{proof:packet-cross-terms}
 \int_\R\|v_z(t-s_1)\|_{X_2}\|v_{z'}(t-s_2)\|_{X_2}\,dt
 \lesssim\frac{1}{1+|s_1-s_2|}.
\end{equation}
The cross term of the two current groups is $O(\eta_j^2T_j^{-1})$ in
$L^1(I_j)$. The groups from the other intervals contribute
$O(\eta_j^2T_0^{-1})$ by \eqref{proof:off-block-tail} with $\sigma=1$.
Since $|e_j-t|\le4T_j$ and $d=T_j\eta_j^2$, both normalized errors
are $O(T_0^{-1})$, uniformly in $N$. This proves
\eqref{i3:eq:uniform-reference}. Continuity follows from the finite sum
in \eqref{i3:eq:train} and the continuity of the amplitudes.
\end{proof}
We fix the family first, then take $\varepsilon=d^3$ and $T_0=d^{-5}$
with $d$ small. The block bound $O(\eta_j)$ applies to the observations,
but not to the energy supremum.

\section{The nonlinear correction along a prescribed path}\label{part:nonlinear}

\subsection{Setting and the prescribed-path theorem}

We fix the background and the radiation family of the preceding
sections. We write
\begin{equation*}
P^*=I-\Phi^*K^*,\qquad \cD=P\Phi^*=\Phi^*-K\Phi\Phi^*.
\end{equation*}
For the chart \eqref{i3:eq:chart}, we write
\begin{equation*}
 \mathcal Q=(Q_i)=K+DR(\xi),\qquad V_\xi=f'(q_\xi).
\end{equation*}

In addition to \eqref{short:chart-identities}, we use
$LR_i=(V_\xi-V)Q_i$. For an interval $J$, we supplement the norms of
Section~\ref{part:linear} with
\begin{equation*}
E_J(u)=\|\nabla u\|_{L_t^\infty L_x^2(J)},\quad
M_J(u)=\|u\|_{L_t^\infty L_x^2(J)}.
\end{equation*}

For the blocks of Section~\ref{i3:sec:blocks} and the radiation train
\eqref{i3:eq:train}, we take
\[
 0<d\le1,\qquad \varepsilon=d^3,\qquad T_0=d^{-5},\qquad
 \alpha=\Phi v_a,\qquad r=Pv_a .
\]
For fixed $D_*,C_2,C_\infty$, we require the path $Y$ to satisfy
\begin{equation}\label{i4:eq:path4}
 \begin{gathered}
 Y\in W^{1,\infty}([0,T_N];\R^{15}),\qquad Y(0)=0,\\
 \|Y\|_\infty\le D_*d,\qquad
 \|Y'\|_{L^2(I_j)}\le C_2\sqrt d\,\eta_j,\qquad
 \|Y'\|_\infty\le C_\infty d\eps.
 \end{gathered}
\end{equation}
We use the Lipschitz representative of $Y$, with derivative bounds
understood almost everywhere.
We set $\xi=Y+\alpha$ and
\begin{equation*}
M_\xi=I+\Phi DR(\xi),\quad
p=-M_\xi^{-*}DR(\xi)^*\bar v,\quad
 \widetilde v=\bar v+\cD p.
\end{equation*}
The matrix $M_\xi$ is invertible for small $d$.
We define
\begin{gather}
N_q(w)=f(q+w)-f(q)-f'(q)w,\quad B_q(w)=N_q(w)-f(w),\notag\\
 \cA=(V_\xi-V)w+N_{q_\xi}(w),\qquad
H_i=\xi_j'\langle R_{ij}(\xi),\widetilde v\rangle
       +\langle Q_i(\xi),N_{q_\xi}(w)\rangle.
 \label{i4:eq:nonlinearity4}
\end{gather}
We sum repeated chart indices.

\begin{theorem}[Uniform nonlinear prescribed-path field]\label{i4:thm:field4}
For all sufficiently small $d$, every $N$, every admissible control $a$,
and every path satisfying \eqref{i4:eq:path4}, there is a unique small field in the
contraction class solving
\begin{align}
w_t&=\bar v+\cD p-PDR(\xi)\xi',\notag\\
 \bar v_t&=-PLw+PP^*\cA+\cD H,\qquad
 \Phi w=\Phi\bar v=0.\label{i4:eq:field4}
\end{align}
Its complementary initial data and mixed negative-mode conditions have the form
\begin{align}
w(0)&=r(0)+\sum_l z_l\Psi_l,\notag\\
 \bar v(0)&=P\partial_tv_a(0)+\sum_l\kappa_l z_l\Psi_l-\cD p(0),
 \label{i4:eq:data4}
\end{align}
where the Green operator determines $z_l$ by the stable initial and
unstable terminal conditions of Section~\ref{part:linear},
\[
 b_l(0)=0,\qquad a_l(T_N)=0 .
\]
The field satisfies, uniformly in $N,a,Y$,
\begin{equation}\label{i4:eq:mainfield4}
 \|(w,\bar v)\|_{\cR([0,T_N])}\lesssim \eps,\qquad
 \|(w-r,\bar v-P\partial_tv_a)\|_{\cR([0,T_N])}\lesssim d\eps,
 \qquad |z|\lesssim d\eps.
\end{equation}
The same field bounds hold with $\widetilde v$ in place of $\bar v$.
Moreover,
\begin{equation}\label{i4:eq:charges4}
 |p(0)|+\|H\|_{L^1(0,T_N)}\lesssim \eps^2,
 \qquad p'+H-K^*\cA=0
\end{equation}
in distributions. The field depends continuously on $(a,Y)$ for each
fixed finite horizon, with the admissible paths given the $C^0$ topology.
\end{theorem}

\subsection{A priori estimates and construction of the field}

To construct the field by contraction, we use the path and radiation
bounds
\begin{equation*}
 \|\xi\|_\infty\lesssim d,\qquad \|\xi'\|_2\lesssim\varepsilon.
\end{equation*}
Every solution is the Green response with complementary data
$(r(0),Pv_{a,t}(0)-\mathcal Dp(0))$ and sources
\begin{equation}\label{short:basic-sources}
 \begin{split}
 \beta&=D^2(PDR\xi'-\mathcal Dp),\qquad
 F_{2,\mathrm b}=P^*((V_\xi-V)w+B_{q_\xi}(w)),\\
 F_{1,\mathrm b}&=P^*f(w)+\Phi^*H.
 \end{split}
\end{equation}
We have the compatibility identities
\begin{equation*}
 \Phi G_0\beta=K^*F_{2,\mathrm b}
 =\Phi\mathcal D=\Psi_l^*\mathcal D=0,\qquad K^*F_{1,\mathrm b}=H.
\end{equation*}
For $I\subset[0,T_N]$, we write
\begin{equation*}
 \mathfrak S_I(w,\bar v)
 =\|\beta\|_{L^2Y_4(I)}+\|F_{2,\mathrm b}\|_{L^2Y_1(I)}
                       +\|F_{1,\mathrm b}\|_{L^1L^2(I)},
\end{equation*}
and omit $I$ on the full horizon. We set
$\mathfrak X=\|(w,\bar v)\|_{\mathcal R}$.
By Newton duality and the chart bounds, we have
\begin{equation*}
 |\langle R_i,y\rangle|\lesssim d\|G_0y\|_{X_4},\qquad
 |\langle R_{ij},y\rangle|\lesssim\|G_0y\|_{X_4},\qquad
 \|p\|_2\lesssim d\mathfrak X.
\end{equation*}
For the nonlinear terms, we use
\begin{equation}\label{short:nonlinear-pointwise}
 |B_q(w)|\lesssim\min(|q|^{\frac{1}{3}}|w|^2,|q||w|^{\frac{4}{3}}),\qquad
 |\partial_wB_q(w)|\lesssim\min(|q|^{\frac{1}{3}}|w|,|q||w|^{\frac{1}{3}}),
\end{equation}
which follow from the $\frac{1}{3}$-H\"older continuity of $f''$.
By \eqref{proof:critical-L4-interpolation} and $\rho|q|^{\frac{1}{3}}\lesssim1$,
we have
\begin{equation*}
 \|B_q(w)\|_{L^2Y_1}\lesssim\|w\|_{L^4L^4}^2\lesssim\mathfrak X^2,\qquad
 \|f(w)\|_{L^1L^2}=S(w)^{\frac{7}{3}}\le\mathfrak X^{\frac{7}{3}}.
\end{equation*}
Since $|Q_i||q|^{\frac{1}{3}}\lesssim\rho^{-4}$, we obtain
\begin{equation}\label{proof:field-source-bounds}
 \begin{aligned}
 |p(0)|+\mathfrak S(w,\bar v)
 &\lesssim d(\varepsilon+\mathfrak X)+\mathfrak X^2+\mathfrak X^{\frac{7}{3}},\\
 |p(0)|+\|H\|_1
 &\lesssim\varepsilon\mathfrak X+\mathfrak X^2+\mathfrak X^{\frac{7}{3}}.
 \end{aligned}
\end{equation}
For the initial trace, we use $\xi(0)=\alpha(0)=O(\varepsilon)$.
For the geometric term in $H$, we use Cauchy--Schwarz.
The source estimates apply to arbitrary constrained fields.
We apply the Green map to the ball
\begin{equation*}
 \|(w-r,\bar v-Pv_{a,t})\|_{\mathcal R}\le C_1d\varepsilon.
\end{equation*}
By the homogeneous bound, every field in this ball satisfies
$\mathfrak X\lesssim\varepsilon$ for $C_1d\le1$, with a constant
independent of $C_1,N$. By \eqref{proof:field-source-bounds}, we have
\begin{equation*}
 |p(0)|+\mathfrak S(w,\bar v)\lesssim d\varepsilon,\qquad
 |p(0)|+\|H\|_1\lesssim\varepsilon^2.
\end{equation*}
We first fix $C_1$ large, then choose $d$ small. With these choices,
\eqref{short:green-bound} maps the ball into itself for every $N$.

For two constrained fields of $\mathcal R$ norm $O(\varepsilon)$ on the same path,
we write $\delta$ for their difference. On every $I\subset[0,T_N]$,
we have, by the mean value formula and the preceding estimates,
\begin{equation}\label{read:nonlinear-differences}
 \|\delta\beta\|_{L^2Y_4(I)}
 +\|\delta F_{2,\mathrm b}\|_{L^2Y_1(I)}
 +\|\delta F_{1,\mathrm b}\|_{L^1L^2(I)}
 \lesssim d\|(\delta w,\delta\bar v)\|_{\mathcal R(I)}.
\end{equation}
At $t=0$, we also have
\begin{equation*}
 |\delta p(0)|\lesssim\varepsilon\|\delta\bar v(0)\|_2.
\end{equation*}
By \eqref{read:nonlinear-differences}, the Green map has Lipschitz constant
$Cd<\frac{1}{2}$. We obtain the unique field in this ball by contraction,
together with \eqref{i4:eq:mainfield4}.
We use \eqref{short:negative-green} for the bound on $z$.
Since $\Psi_l^*\mathcal D=0$, its initial data satisfy
\begin{equation}\label{short:field-initial}
 w(0)=r(0)+\sum_lz_l\Psi_l,\qquad
 \widetilde v(0)=Pv_{a,t}(0)+\sum_l\kappa_lz_l\Psi_l.
\end{equation}
By the identities
\begin{equation*}
 DR^*\mathcal D=M_\xi^*-I,\qquad PP^*=P-\mathcal DK^*,
\end{equation*}
we have $p=-DR^*\widetilde v$. Differentiating this constraint and using
$LR_i=(V_\xi-V)Q_i$, we compute
\begin{equation*}
 \begin{aligned}
 \bigl[M_\xi^*(p'+H-K^*\mathcal A)\bigr]_i
 &=\langle LR_i,w\rangle+\langle Q_i,N_{q_\xi}(w)-\mathcal A\rangle\\
 &=0.
 \end{aligned}
\end{equation*}
The inclusions
\[
 \bar v_t\in L^1_{\rm loc}H^{-1},\qquad
 R_i(\xi)\in W^{1,\infty}H^1
\]
justify the product rule by $H^{-1}$--$H^1$ duality. Since $M_\xi$ is invertible,
we obtain $p\in W^{1,1}_{\rm loc}$ and
\begin{equation}\label{short:field-identities}
p'+H-K^*\mathcal A=0,\qquad
w_t=\widetilde v-PDR\xi',\qquad
 \widetilde v_t=-PLw+P\mathcal A.
\end{equation}

\begin{lemma}[Continuity in the uniform topology]\label{i4:lem:uniform-path}
For a fixed horizon and fixed path bounds, the field in $\mathcal R$,
its selected negative coefficients, and $p$ in $C^0$ depend continuously
on the controls and on the path in $C^0$.
\end{lemma}
\begin{proof}
We set
\begin{equation*}
 W=w+PR(\xi),\qquad
 U=Q_*+W+K(\xi+\Phi R(\xi))=q_\xi+w.
\end{equation*}
The transformed field satisfies
\begin{equation}\label{proof:path-free-equation}
W_t=\widetilde v,\qquad W_{tt}=P[\Delta U+f(U)].
\end{equation}
For convergent inputs, we pass to a subsequence with $z_n\to z$.
Then \eqref{short:field-initial} gives convergence of $\vec W_n(0)$.
In \eqref{proof:path-free-equation}, the path enters only through $\xi$. The chart is Lipschitz in $\xi$, $P$ is
bounded on $L^2$, and
\begin{equation*}
 P\Delta W-\Delta W=K\Phi D^2W,\qquad
 K\Phi D^2:\dot H^1\to L^2.
\end{equation*}
For small $\varepsilon$, we choose a common subdivision with
\begin{equation*}
 |J|\ll1,\qquad
 \sup_n S_J(U_n)\le C|J|^{\frac{3}{7}}+C\varepsilon\ll1.
\end{equation*}
By H\"older and the energy--Strichartz estimate \cite[(2.2)]{DKM16},
we absorb the nonlinear and finite-rank terms on each $J$ and iterate
to obtain
\begin{equation}\label{proof:path-continuity}
 \|\vec W_n-\vec W_m\|_{C\mathcal H}+S(W_n-W_m)
 \lesssim_{T_N}\bigl(\|\vec W_n(0)-\vec W_m(0)\|_{\mathcal H}
                         +\|\xi_n-\xi_m\|_\infty\bigr).
\end{equation}
The continuity constant may depend on $T_N$, whereas the bounds in
Theorem~\ref{i4:thm:field4} remain uniform in $N$.
By Hardy's inequality and
\begin{equation*}
 \|G_0y\|_{X_4}\lesssim\|y\|_2,\qquad
 w=W-PR(\xi),\qquad p=-DR^*\widetilde v,
\end{equation*}
we recover convergence of the original fields and observations.
The limit remains in the closed contraction class, satisfies the mixed
conditions and \eqref{short:field-identities}, and hence solves
\eqref{i4:eq:field4}. By uniqueness, the full sequence converges.
\end{proof}

\subsection{The Newton monopole and stronger source bounds}

To bound the remainder in $L^2$ after subtracting the Newton monopole,
we use two weighted estimates
\begin{align}\label{short:mixed-estimates}
 \|B_q(w)\|_{L^2Y_s}&\lesssim
 \LE(w)^{\frac{2}{3}}\|w\|_{L^4L^4}^{\frac{2}{3}},
 \quad s\le\frac{5}{3},\notag\\
 \|B_q(w)\|_{L^2W^{1,\frac{14}{11}}}&\lesssim E(w)^{\frac{1}{3}}(\LE(w)+\LG(w)).
\end{align}
For the first estimate, we apply H\"older to
\begin{equation*}
 \rho^s|B_q(w)|\lesssim
 \rho^{s-\frac{5}{3}}(\rho^{-2}|w|)^{\frac{2}{3}}|w|^{\frac{2}{3}}.
\end{equation*}
For the second estimate, we use
\begin{equation*}
 |\nabla B_q(w)|\lesssim|\nabla q||w|^{\frac{4}{3}}+|q||w|^{\frac{1}{3}}|\nabla w|
\end{equation*}
and H\"older with
\begin{equation*}
 \rho^2q,\rho^2\nabla q,\rho q\in L^{\frac{70}{13}},\qquad
 \frac{1}{2}+\frac{1}{10}+\frac{13}{70}=\frac{11}{14}.
\end{equation*}
Neither bound uses $\|w\|_2$.

Projecting only the velocity, we have
\begin{equation}\label{short:mass-system}
 \begin{aligned}
 w_t&=P_0\bar v-G_0(\beta-\nu K^*\bar v),\\
 (P_0\bar v)_t&=\Delta w+Vw+P^*\mathcal A+(\Phi^*-K)H.
 \end{aligned}
\end{equation}
Here
\begin{equation*}
 K^*\bar v=\nu^*G_0\bar v,\qquad
 \|\beta-\nu K^*\bar v\|_{L^2Y_4}\lesssim\varepsilon,
\end{equation*}
and
\begin{equation*}
 P_0\bar v(0)=v_{a,t}(0)+\sum_l\kappa_lz_l\Psi_l
             -\Phi^*p(0)+G_0\nu p(0).
\end{equation*}

With the same kernel $\kappa_0$, we put
\begin{equation}\label{short:halo}
h(t)=m(\nu p(0))\kappa_0(t)-\int_0^tm(\nu H(s))\kappa_0(t-s)\,ds,
 \qquad G=w-h.
\end{equation}
By \eqref{i4:eq:charges4}, the driving measure has total variation $\lesssim\varepsilon^2$.
Using the kernel in \eqref{short:newton-kernels}, we obtain, for large $t$,
\begin{equation*}
 \begin{aligned}
 |\kappa_0(t,x)|&\lesssim t|x|^{-3}\mathbf1_{|x|\ge t-R},\\
 |\nabla\kappa_0(t,x)|&\lesssim
 t|x|^{-4}\mathbf1_{|x|\ge t-R}
 +t^{-2}\mathbf1_{||x|-t|\le R}.
 \end{aligned}
\end{equation*}
By radial integration, including bounded times by convolution with
$\psi$, we have
\begin{equation*}
 \begin{aligned}
 \|\rho^\sigma\kappa_0(t)\|_p
 &\lesssim\langle t\rangle^{\sigma+\frac{5}{p}-2},
 &&0\le\sigma<3-\frac{5}{p},\\
 \|\nabla\kappa_0(t)\|_p
 &\lesssim\langle t\rangle^{\frac{5}{p}-3}
          +\langle t\rangle^{\frac{4}{p}-2},
 &&p>\frac{5}{4}.
 \end{aligned}
\end{equation*}
We apply the free estimate and Minkowski before taking powers to obtain
\begin{align}\label{short:halo-norms}
 &\|(h,h_t)\|_{\mathcal R}+\|\rho^{\frac{13}{8}} h\|_{L^4L^\infty}
   +\|\nabla h\|_{L^4L^{\frac{7}{2}}}+\|\rho h\|_{L^\infty L^{\frac{70}{13}}}
   \lesssim \varepsilon^2,\notag\\
 &\|\rho^{\frac{13}{8}}|h|^{\frac{7}{3}}\|_{L^2L^2}\lesssim \varepsilon^{\frac{14}{3}},\qquad
 \|\rho|h|^{\frac{4}{3}}\|_{L^\infty L^{\frac{7}{2}}}\lesssim \varepsilon^{\frac{8}{3}}.
\end{align}

To prove $M(G)\lesssim\varepsilon$, we decompose
\begin{equation*}
 \mathcal A=(V_\xi-V)w+B_q(w)+[f(w)-f(G)]+f(G).
\end{equation*}
The Newton density and the remaining finite-rank force are, respectively,
\begin{equation*}
 -\nu H,\qquad \Phi^*(H-K^*f(G)).
\end{equation*} 
After subtracting $h$, we control each density $g-m(g)\psi$ by
\eqref{proof:mean-zero-newton}, with contribution $O(\varepsilon^2)$.
The other data are $O(\varepsilon)$ in $L^2\times\dot H^{-1}$ by
\eqref{i3:R:minusone}. For the sources, we use
\eqref{proof:free-position-duality}, Minkowski, and
$\|D^{-1}\Phi^*c\|_2\lesssim|c|$.

On a fixed horizon, we have $M(G)<\infty$ because
\[
 w(0)\in L^2,\qquad
 w_t=\bar v-G_0\beta\in L^1L^2,\qquad h\in C([0,T_N];L^2).
\]
By \eqref{short:mixed-estimates}, the weighted terms other than
$f(w)-f(G)$ are $O(\varepsilon)$. By interpolation, we have
\begin{equation}\label{proof:mass-interpolation}
 \|G\|_{L^{\frac{16}{3}}L^{\frac{8}{3}}}
 \le M(G)^{\frac{9}{16}}S(G)^{\frac{7}{16}}.
\end{equation}
By H\"older, \eqref{short:halo-norms}, and Young's inequality, we obtain
\begin{align*}
 \|f(G+h)-f(G)\|_{L^2Y_{\frac{13}{8}}}
 &\lesssim\varepsilon^2M(G)^{\frac{3}{4}}S(G)^{\frac{7}{12}}
             +\varepsilon^{\frac{14}{3}}\\
 &\lesssim\varepsilon^{\frac{7}{3}}(M(G)+\varepsilon),\\
 \|f(G)\|_{L^{\frac{7}{4}}L^{\frac{14}{11}}}
 &\lesssim\varepsilon^{\frac{4}{3}}M(G).
\end{align*}
The compact coefficient is $O(\varepsilon^2)$ in $L^1$. Hence
\begin{equation}\label{proof:mass-closure}
 M(G)\lesssim\varepsilon+\varepsilon^{\frac{4}{3}}M(G).
\end{equation}
For small $\varepsilon$, we absorb the last term and conclude $M(G)\lesssim\varepsilon$,
uniformly in $T_N$.

\subsubsection{Stronger source bounds}

By \eqref{proof:mass-interpolation}, H\"older, the chain rule, and
$M(G)+E(G)+S(G)\lesssim\varepsilon$,
we have
\begin{equation}\label{short:strong-sources}
 \begin{aligned}
 &\|f(G)\|_{L^1L^2}
 +\|f(G)\|_{L^{\frac{7}{4}}W^{1,\frac{14}{11}}}
 +\|f(G)\|_{L^4W^{1,1}}
 \lesssim\varepsilon^{\frac{7}{3}}.
 \end{aligned}
\end{equation}
By the mass bound, \eqref{short:halo-norms}, and $\LG(G)\lesssim\varepsilon$,
we obtain
\begin{equation}\label{proof:halo-source-bounds}
 \|f(w)-f(G)\|_{L^2(Y_{\frac{13}{8}}\cap W^{1,\frac{14}{11}})}
 \lesssim\varepsilon^{\frac{10}{3}}.
\end{equation}
For the derivative, we apply H\"older to
\begin{equation*}
 \begin{aligned}
 |\nabla(f(G+h)-f(G))|
 &\lesssim(|G|^{\frac{1}{3}}|h|+|h|^{\frac{4}{3}})|\nabla G|\\
 &\quad+(|G|^{\frac{4}{3}}+|h|^{\frac{4}{3}})|\nabla h|,
 \end{aligned}
\end{equation*}
using \eqref{proof:mass-interpolation} and
\begin{equation*}
 \frac{1}{10}+\frac{13}{70}+\frac{1}{2}
 =\frac{2}{7}+\frac{1}{2}=\frac{11}{14}.
\end{equation*}
For the term without derivatives, we use the embedding
$Y_{13/8}\hookrightarrow L^{14/11}$.

At the nodal set, we have
\begin{equation*}
 |V_\xi-V|\lesssim d\rho^{-4},\qquad
 |\nabla(V_\xi-V)|\lesssim d^{\frac{1}{3}}\rho^{-5}.
\end{equation*}
By H\"older, we obtain
\begin{equation*}
 \|(V_\xi-V)w\|_{L^2(Y_{\frac{13}{8}}\cap W^{1,\frac{14}{11}})}
 \lesssim d^{\frac{1}{3}}(\LE(w)+\LG(w))\lesssim\varepsilon.
\end{equation*}

We use the enhanced source decomposition
\begin{equation}\label{short:enhanced}
 \begin{aligned}
 F_{2,\mathrm e}&=P^*[(V_\xi-V)w+B_q(w)+f(w)-f(G)],\\
 F_{1,\mathrm e}&=f(G)+\Phi^*(H-K^*f(G)).
 \end{aligned}
\end{equation}
The basic split \eqref{short:basic-sources} and this split satisfy
\begin{equation}\label{proof:source-regrouping}
 \begin{gathered}
 F_{2,\mathrm e}-F_{2,\mathrm b}
 =P^*(f(w)-f(G))=F_{1,\mathrm b}-F_{1,\mathrm e},\\
 K^*F_{2,\mathrm b}=K^*F_{2,\mathrm e}=0,\qquad
 K^*F_{1,\mathrm b}=K^*F_{1,\mathrm e}=H.
 \end{gathered}
\end{equation}
Thus the total force and charge are unchanged, also after multiplication
by a common scalar cutoff. Since $\Phi^*=G_0D^2\Phi^*$, we also have
\begin{equation*}
 \|D^2\Phi^*(H-K^*f(G))\|_{L^1Y_4}\lesssim\varepsilon^2.
\end{equation*}
By \eqref{short:mixed-estimates}, \eqref{short:strong-sources},
\eqref{proof:halo-source-bounds}, and boundedness of $P^*$ on
$Y_{\frac{13}{8}}$ and $W^{1,\frac{14}{11}}$, the source norms in
Theorem~\ref{i2:thm:old}, together with $|p(0)|$, sum to
$O(\varepsilon)$.
The correction $(w-r,\bar v-Pv_{a,t})$ has complementary initial data
\begin{equation*}
 (0,-\mathcal Dp(0))=(0,G_0[-D^2\Phi^*p(0)+\nu\Phi\Phi^*p(0)]).
\end{equation*}

\subsection{Old and current responses}

For each $j\ge1$, we split the sources in \eqref{short:enhanced} at
$T_j/2$ by setting
\begin{equation*}
 \chi_j=\mathbf1_{[0,\frac{T_j}{2}]}.
\end{equation*}
We write $o_j$ for the response to the complementary data and the
sources up to $T_j/2$. We write $c_j$ for the response to the remaining
sources with zero complementary data. We select the negative modes for
each response separately, so linearity gives
\begin{equation*}
 (w-r,\bar v-Pv_{a,t})=o_j+c_j.
\end{equation*}
The unstable modes of $c_j$ may be nonzero before its forcing begins,
as in \eqref{short:negative-green}. By Theorem~\ref{i2:thm:old}, we have
\begin{equation}\label{short:old-bounds}
 \begin{aligned}
 \LE_{I_j}(o_j)+\OO_{I_j}(y_{o_j})&\lesssim\varepsilon T_j^{-\frac{5}{8}},
 &S_{I_j}(o_j)&\lesssim\varepsilon T_j^{-\frac{3}{14}},\\
 \|o_j\|_{L^4(I_j;L^4)}&\lesssim\varepsilon T_j^{-\frac{1}{4}},
 &\|o_j\|_{L^2(I_j;X_{\frac{5}{3}})}&\lesssim\varepsilon T_j^{-\frac{1}{2}}.
 \end{aligned}
\end{equation}
For $j=0$, we put $o_0=0$ and
$c_0=(w-r,\bar v-Pv_{a,t})$.

We set
\begin{equation*}
 \mu=\max_{j<N}\eta_j^{-1}\|c_j\|_{\mathcal R}<\infty.
\end{equation*}
Independently of $\mu$, we have
\begin{equation*}
 \|o_j\|_{\mathcal R}+\|c_j\|_{\mathcal R}\lesssim\varepsilon.
\end{equation*}
Since $T_j\eta_j^2=d$, $T_j\ge d^{-5}$, and $\varepsilon=d^3$, we have
\begin{equation*}
 \varepsilon T_j^{-\frac{5}{8}}
 \le d^{\frac{25}{8}}\eta_j,\qquad
 \LE_{I_j}(w)+\OO_{I_j}(\bar v)\lesssim\eta_j(1+\mu).
\end{equation*}
For the reference pair $(r+o_j,Pv_{a,t}+y_{o_j})$, we use
\begin{equation*}
 |B_q(r+o_j)|\lesssim|q|^{\frac{1}{3}}|r|^2+|q||o_j|^{\frac{4}{3}}.
\end{equation*}
By \eqref{short:mixed-estimates} and \eqref{short:old-bounds}, we have
\begin{equation*}
 \||q||o_j|^{\frac{4}{3}}\|_{L^2Y_1(I_j)}
 \lesssim\varepsilon^{\frac{4}{3}}T_j^{-\frac{7}{12}},\qquad
 S_{I_j}(o_j)^{\frac{7}{3}}\lesssim\varepsilon^{\frac{7}{3}}T_j^{-\frac{1}{2}}.
\end{equation*}
Since $\frac{7}{12}>\frac{1}{2}$, we may use $T_j^{-\frac{1}{2}}$ for both terms.
By the chart bounds, \eqref{i3:eq:blocknorms}, \eqref{short:mixed-estimates},
and \eqref{short:old-bounds}, we obtain
\begin{equation*}
 \mathfrak S_{I_j}(r+o_j,Pv_{a,t}+y_{o_j})
 \lesssim d\eta_j+\eta_j^2+\varepsilon^{\frac{4}{3}}T_j^{-\frac{1}{2}}
 \lesssim d\eta_j.
\end{equation*}
Here we use
\begin{equation*}
 \eta_j\le\varepsilon=d^3,\qquad
 \varepsilon^{\frac{4}{3}}T_j^{-\frac{1}{2}}=d^{\frac{7}{2}}\eta_j.
\end{equation*}
The reference pair and $(w,\bar v)$ have global $\mathcal R$ norm
$O(\varepsilon)$, so we apply \eqref{read:nonlinear-differences} with
difference $c_j$. Thus
\begin{equation}\label{proof:interval-source-bound}
 \mathfrak S_{I_j}(w,\bar v)
 \lesssim d\eta_j+d\|c_j\|_{\mathcal R(I_j)}
 \lesssim d\eta_j(1+\mu).
\end{equation}
By \eqref{proof:source-regrouping}, the sources
\[
 (1-\chi_j)(\beta,F_{2,\mathrm b},F_{1,\mathrm b})
\]
give the same response $c_j$. We apply the Green estimate. On each
block meeting $[T_j/2,T_N]$, we use that block's decomposition, so
\eqref{short:old-bounds} retains its time separation. We sum by
\begin{equation*}
 \sum_{k\ge j-1}\eta_k+
 \left(\sum_{k\ge j-1}\eta_k^2\right)^{\frac{1}{2}}\lesssim\eta_j.
\end{equation*}
For $j=0$, we have $\|c_0\|_{\mathcal R}\lesssim d\eta_0$
by the a priori correction bound \eqref{i4:eq:mainfield4}.
For small $d$, we absorb $d\mu$ to obtain
\begin{equation}\label{proof:current-response}
 \mu\lesssim d(1+\mu),\qquad \mu\lesssim d.
\end{equation}
No cutoff derivative occurs, since we cut off only the external sources.

\subsection{Comparison of the forces}

To compare the two normalized moments, we set
\begin{equation*}
 b_{ik}=\langle R_{ik}(0),v_{a,t}\rangle,\qquad
 \mathcal F=H-(b\alpha)'.
\end{equation*}
By \eqref{i3:cons:gaugeforce}, we write
\begin{equation}\label{short:force-residual}
 (\mathcal F-F(v_a))_i
 =\langle Q_i,N_q(w)\rangle-F_i(r)
 +\xi'_k\langle R_{ik}(\xi),\widetilde v\rangle-\alpha'_kb_{ik}.
\end{equation}
\begin{proposition}[Force comparison]\label{i4:four-res:proposition}
On every prescribed path, we have, uniformly in $j,N$,
\begin{equation}\label{read:force-error}
 \|\mathcal F-F(v_a)\|_{L^1(I_j)}\lesssim d^{\frac{1}{3}}\eta_j^2.
\end{equation}
The normalized first-moment error is also $O(d^{\frac{1}{3}})$.
\end{proposition}
\begin{proof}
By Taylor's formula and the chart bounds, we have
\begin{equation}\label{proof:force-taylor}
 \begin{aligned}
N_q(w)&=\frac{1}{2}f''(q)w^2+O(|w|^{\frac{7}{3}}),\\
 |Q_if''(q)-K_if''(Q_*)|&\lesssim d^{\frac{1}{3}}\rho^{-4}.
 \end{aligned}
\end{equation}
By \eqref{i3:eq:blocknorms}, \eqref{short:old-bounds}, and
\eqref{proof:current-response}, we have
\begin{equation*}
 \|F(w)-F(r)\|_{L^1(I_j)}\lesssim d\eta_j^2.
\end{equation*}
For the higher power, we write $w=(r+c_j)+o_j$.
By \eqref{short:old-bounds} and H\"older, we obtain
\begin{align}
 \iint_{I_j}\rho^{-3}|r+c_j|^{\frac{7}{3}}
 &\lesssim S_{I_j}(r+c_j)^{\frac{7}{3}}\lesssim\eta_j^{\frac{7}{3}},\notag\\
 \iint_{I_j}\rho^{-3}|o_j|^{\frac{7}{3}}
 &\le\|o_j\|_{L^2(I_j;X_{\frac{5}{3}})}^{\frac{5}{3}}
       \|o_j\|_{L^4(I_j;L^4)}^{\frac{2}{3}}
 \lesssim\varepsilon^{\frac{7}{3}}T_j^{-1}.
 \label{proof:old-remainder-force}
\end{align}
Since $T_j\eta_j^2=d$, $\eta_j\le\varepsilon$, and $\varepsilon=d^3$, we conclude
\begin{equation*}
 \eta_j^{-2}\iint_{I_j}\rho^{-3}|w|^{\frac{7}{3}}
 \lesssim\eta_j^{\frac{1}{3}}+\frac{\varepsilon^{\frac{7}{3}}}{d}
 \lesssim d.
\end{equation*}

For the geometric term, we differentiate $K^*R(\xi)=0$ twice to obtain
\begin{equation*}
 K^*R_{ik}(0)=0,\qquad
 b_{ik}=\langle R_{ik}(0),Pv_{a,t}\rangle.
\end{equation*}
We use $\xi'=Y'+\alpha'$, Newton duality, and
\eqref{proof:chart-holder}. All norms below are on $I_j$.
By Cauchy--Schwarz, we have
\begin{equation}\label{proof:geometric-force-bound}
 \begin{aligned}
 &\int_{I_j}\left|\xi'_k\langle R_{ik}(\xi),\widetilde v\rangle
                         -\alpha'_kb_{ik}\right|\,dt\\
 &\quad\lesssim
  (\|Y'\|_2+d^{\frac{1}{3}}\|\alpha'\|_2)\|G_0\widetilde v\|_{L^2X_4}
  +\|\alpha'\|_2\|G_0(\widetilde v-Pv_{a,t})\|_{L^2X_4}\\
 &\quad\lesssim d^{\frac{1}{3}}\eta_j^2.
 \end{aligned}
\end{equation}
Together with \eqref{proof:force-taylor}, these estimates prove
\eqref{read:force-error}. For the first moment, we have
\begin{equation*}
 \begin{aligned}
 d^{-1}\left|\int_{I_j}(e_j-t)[\mathcal F(t)-F(v_a(t))]\,dt\right|
 &\le\frac{4T_j}{d}\|\mathcal F-F(v_a)\|_{L^1(I_j)}\\
 &\lesssim d^{\frac{1}{3}},\qquad T_j\eta_j^2=d.
 \end{aligned}
\end{equation*}
\end{proof}

\section{Construction on finite time intervals}\label{part:endpoint}

We prove Proposition~\ref{presentation:finite} by imposing the
endpoint conditions at a fixed point and then showing that the
corresponding initial data are precompact.

\subsection{A compact endpoint map}\label{sec:endpoint-selection}

We use $\xi+\Phi R(\xi)$ as the coordinate for the endpoint conditions.
With the scales and fields of Section~\ref{part:nonlinear}, we write
$I_j=[s_j,s_{j+1}]$, where $s_0=0$, $s_j=T_j$ for $1\le j\le N$, and
$e_j=s_{j+1}$ for $0\le j<N$. For the alternating endpoints $\theta_j$ of
Proposition~\ref{presentation:finite}, we choose the central controls
\begin{equation*}
a_j^*=(-1)^j(e,-e).
\end{equation*}
The controls satisfy
\begin{equation*}
 |a_j-a_j^*|_\infty\le r_0<\frac{1}{2}.
\end{equation*}
For fixed $N$, we take all paths satisfying \eqref{i4:eq:path4}.
The path domain contains the zero path and is compact and convex in the
uniform topology, by Arzel\`a--Ascoli and weak lower semicontinuity.
By Theorems~\ref{i3:thm:controls}, \ref{i4:thm:field4} and
Proposition~\ref{i4:four-res:proposition}, with $T_0^{-1}=d^5$,
we have, uniformly in $N$,
\begin{align}
 \|\mathcal F\|_{L^1(I_j)}&\lesssim\eta_j^2,\notag\\
 \left(\eta_j^{-2}\int_{I_j}\mathcal F,\
d^{-1}\int_{I_j}(e_j-t)\mathcal F\right)
 &=\mathcal L a_j+O(d^{\frac{1}{3}}).
 \label{short:endpoint-moments}
\end{align}
We set
\begin{equation}\label{short:endpoint-k}
 \begin{split}
k&=-\mathcal Q^*\widetilde v+b\alpha-\alpha'\\
  &=-DR^*v_{a,t}-\mathcal Q^*(\widetilde v-Pv_{a,t})+b\alpha.
 \end{split}
\end{equation}
Here we use $K^*v_{a,t}=0$ and $\mathcal Q^*K=I$.
By the chart density bounds and correction estimates, we have
\begin{equation}\label{short:k-bounds}
 \|k\|_{L^2(I_j)}\le C_kd\eta_j,\qquad
 \|k\|_\infty\le C_{k,\infty}d\varepsilon.
\end{equation}
We use Newton duality as in \eqref{proof:field-source-bounds}, and
\eqref{short:old-bounds}--\eqref{proof:current-response} for the correction.
When $C_2\sqrt d\le1$, the path enters the contraction estimate only
through
\[
 \|\xi\|_\infty\lesssim D_*d+\varepsilon,\qquad
 \|\xi'\|_2\lesssim\varepsilon.
\]
Thus $C_{k,\infty}$ is independent of $C_\infty$.

Since $D(\xi+\Phi R(\xi))=M_\xi$ and $\|M_\xi-I\|<\frac{1}{2}$,
the map $\xi\mapsto\xi+\Phi R(\xi)$ has a local inverse with Lipschitz
constant two. We choose smooth nondecreasing $h_{j1},h_{j2}$ from zero
to one on $I_j$, with derivatives supported near and even about
$2T_j,3T_j$, respectively. Thus
\begin{equation*}
 h_{j\nu}(s_j)=0,\qquad h_{j\nu}(e_j)=1,\qquad
 \int_{I_j}h_{j1}=2T_j,\qquad\int_{I_j}h_{j2}=T_j.
\end{equation*}
We put
\begin{equation*}
 \begin{aligned}
 \mathcal P_j(t)&=\int_{s_j}^t\mathcal F,\qquad 0\le j<N,\\
 \vartheta_j&=\theta_j+\Phi R(\theta_j+\alpha(s_j)),\qquad 0\le j\le N.
 \end{aligned}
\end{equation*}
We define the map by
\begin{gather}
 \widehat a_j=a_j-\mathcal L^{-1}
 \binom{\eta_j^{-2}\mathcal P_j(e_j)}
 {d^{-1}\left[\int_{I_j}(\mathcal P_j+k)
       -(\vartheta_{j+1}-\vartheta_j)\right]},\notag\\
 \widehat\pi(t)=\mathcal P_j(t)
       -\eta_j^2\sum_{\nu=1}^2(a_j-\widehat a_j)_\nu h_{j\nu}(t),\notag\\
 \widehat Y(t)+\Phi R(\widehat Y(t)+\alpha(t))
 =\vartheta_j+\int_{s_j}^t(\widehat\pi+k).
 \label{short:endpoint-map}
\end{gather}
The local inverse defines $\widehat Y$ for small $d$. By the two rows
in the definition of $\widehat a_j$, we have
\begin{equation}\label{short:endpoint-domain}
 \widehat\pi(s_j)=\widehat\pi(e_j)=0,\qquad
 \widehat Y(s_j)=\theta_j,\qquad\widehat Y(e_j)=\theta_{j+1}.
\end{equation}
Differentiating the last identity in \eqref{short:endpoint-map}, we have
\begin{equation*}
 \widehat Y'
 =M_{\widehat Y+\alpha}^{-1}
   [\widehat\pi+k-\Phi DR(\widehat Y+\alpha)\alpha'].
\end{equation*}
Since $\widehat\pi=0$ at each join, the derivatives agree there, so
$\widehat Y\in C^1$.

By Fubini, \eqref{short:endpoint-moments}, \eqref{short:k-bounds}, and
$|\vartheta_j-\theta_j|\lesssim d^2$, we obtain
\begin{equation}\label{proof:control-error}
 |\widehat a_j-a_j^*|
 \lesssim d^{\frac{1}{3}}+C_k\sqrt d.
\end{equation}
We provisionally bound the errors in \eqref{short:endpoint-moments}
and \eqref{proof:control-error} by one. Then we have
\begin{equation*}
 |a_j-\widehat a_j|\lesssim1,\qquad \|\widehat\pi\|_\infty\le C_\pi\eta_j^2,
\end{equation*}
with $C_\pi$ depending only on the fixed radiation family and control
cube, under these provisional error bounds.
We compare $\widehat Y+\alpha(t)$ with $\theta_j+\alpha(t)$ in the local
inverse and use $|DR(\xi)|\lesssim|\xi|$. After increasing
$C_k,C_{k,\infty}$ by chart constants, we obtain
\begin{align}\label{proof:path-preservation}
 \|\widehat Y\|_\infty
 &\le (1+8C_\pi)d+C_kd^{\frac{3}{2}},\notag\\
 \|\widehat Y'\|_{L^2(I_j)}
 &\le4C_\pi\sqrt d\,\eta_j+C_kd\eta_j,\notag\\
 \|\widehat Y'\|_\infty
 &\le2C_\pi\varepsilon^2+C_{k,\infty}d\varepsilon.
\end{align}
We choose the path constants in the order
\begin{equation*}
 \begin{gathered}
D_*>2+8C_\pi,\qquad C_2>4C_\pi+2,\\
C_\infty>2C_{k,\infty}(D_*,C_2)+1,
 \end{gathered}
\end{equation*}
then take $d$ small. We enforce the provisional error
bounds, $C_2\sqrt d\le1$, and
\begin{equation*}
 |\widehat a_j-a_j^*|_\infty\le\frac{r_0}{2}.
\end{equation*}
Since $C_k\sqrt d\to0$ and $\frac{\varepsilon}{d}=d^2\to0$,
\eqref{proof:path-preservation} preserves the path domain.
All choices are independent of $N$.

To prove continuity, we write, by energy duality and the chart identities,
\begin{equation}\label{short:force-primitive}
 \mathcal P_j(t)
 =[\mathcal Q^*\widetilde v-b\alpha]_{s_j}^t
   +\int_{s_j}^t\Phi(-Lw+\mathcal A).
\end{equation}
By Lemma~\ref{i4:lem:uniform-path} and the local inverse, the endpoint
map is continuous in the uniform topology. We apply Schauder's theorem
to obtain
\begin{equation*}
 (\widehat a,\widehat Y)=(a,Y),\qquad Y\in C^1.
\end{equation*}

\subsection{Recovery of the equation and initial data}\label{presentation:recovery}

At the fixed point, the block primitives join into an absolutely
continuous $\pi$, with
\begin{equation}\label{short:modulation}
M_\xi\xi'=\pi-\mathcal Q^*\widetilde v+b\alpha,\qquad
 \pi(s_j)=\pi(e_j)=0.
\end{equation}
By \eqref{short:force-primitive} and \eqref{short:field-identities},
we have
\begin{equation}\label{proof:residual-cancellation}
 (M_\xi\xi')'=\Phi(-Lw+\mathcal A),\qquad
 \widetilde v_t=P(-Lw+\mathcal A).
\end{equation}
By continuity at the joins, we have $M_\xi\xi'\in W^{1,1}$ on the
entire horizon. Setting $u=q_\xi+w$, we use the stationary equation and the identity
$P+K\Phi=I$ to obtain
\begin{equation}\label{short:physical-identities}
 \begin{aligned}
 u_t&=\widetilde v+KM_\xi\xi',\\
 u_{tt}&=(P+K\Phi)(-Lw+\mathcal A)=\Delta u+f(u).
 \end{aligned}
\end{equation}
By the chart bounds, we have
\begin{equation*}
 S_{[0,T]}(u)\lesssim T^{\frac{3}{7}}+S_{[0,T]}(w)<\infty.
\end{equation*}
By \cite[Lemma~2.5]{DKM16}, $u$ is the
critical energy solution.

At $t=0$, we have $Y=\pi=0$ and $\xi=\alpha$, and hence
\begin{align}
 u(0)&=Q_*+v_a(0)+R(\alpha(0))+\sum_lz_l\Psi_l,\notag\\
 u_t(0)&=v_{a,t}(0)+\sum_l\kappa_lz_l\Psi_l
                +K[p(0)+b(0)\alpha(0)],\label{short:physical-data}\\
 |z|&\lesssim d\varepsilon,\qquad
 |p(0)|+|b(0)\alpha(0)|\lesssim\varepsilon^2.\notag
\end{align}
Here we use $p(0)=-DR(\alpha(0))^*\widetilde v(0)$ and
\eqref{short:field-initial}.

We pass to a subsequence with coordinatewise convergent controls and
convergent $z_N$. Since
\begin{equation*}
 \sum_{j\ge J}\eta_j=2^{1-J}\varepsilon\longrightarrow0,
\end{equation*}
the trains converge in energy, $L^2$ position, and every reference norm,
with $\alpha_N\to\alpha$ in $C_b^1$. By \eqref{short:physical-data},
the initial data converge strongly. Together with the field and path
bounds, this proves Proposition~\ref{presentation:finite}.

\section{Global solution and scattering}\label{sec:global-radiation}

\subsection{Passage to a global solution}\label{sec:global-passage}
\label{presentation:parameter-limit}

We fix $T>0$. Whenever $T_N\ge T$, the chart bounds and
\eqref{presentation:finite-bounds} give
\begin{equation*}
 S_{[0,T]}(u_N)\lesssim T^{\frac{3}{7}}+\varepsilon.
\end{equation*}
For $T_M,T_N\ge T$, we apply Strichartz \cite[Corollary~1.3]{KT} to the
difference equation and subdivide $[0,T]$ as in \cite[Section~2]{KM08} to obtain
\begin{equation}\label{i5:lim:stability}
 \|\vec u_N-\vec u_M\|_{C([0,T];\mathcal H)}+S_{[0,T]}(u_N-u_M)
 \lesssim_T\|\vec u_N(0)-\vec u_M(0)\|_{\mathcal H}\longrightarrow0.
\end{equation}
The constant depends on $T$ and is uniform in $M,N$. The limits agree
on overlapping intervals and satisfy Duhamel's formula, so they define
a global energy solution $u$.
By the uniform chart bound and $\|M_\xi-I\|<\frac{1}{2}$, we define
$\xi$ through the local inverse and obtain
\begin{gather*}
 \xi+\Phi R(\xi)=\Phi(u-Q_*),\qquad
 \xi'=M_\xi^{-1}\Phi u_t,\\
 \xi_N\to\xi,\qquad Y_N=\xi_N-\alpha_N\to Y
 \quad\text{in }C^1_{\rm loc},\qquad Y=\xi-\alpha.
\end{gather*}
We set $w=u-q_\xi$. By the chart bounds, we have
$w_N\to w$ in $C_{\rm loc}\dot H^1\cap S_{\rm loc}$.
Using Hardy's inequality and Fatou's lemma, we pass
\eqref{presentation:finite-chart}--\eqref{presentation:finite-bounds}
to $[0,\infty)$. By the uniform convergence of $Y_N$ on compact
intervals, we retain every prescribed endpoint.

\subsection{Scattering and stationary limits}\label{sec:scattering-limits}

\label{presentation:scattering-proof}
For the solution in Proposition~\ref{presentation:scattering}, the
physical pair $(w,u_t)$ satisfies
\begin{equation*}
w_t=u_t-G_0D^2(\mathcal Q\xi'),\qquad
u_{tt}=\Delta w+V_\xi w+B_{q_\xi}(w)+f(w).
\end{equation*}
For $I=[a,b]\subset[0,\infty)$, we apply \eqref{short:free-estimate}
with zero data at $a$. By the chart bounds,
\eqref{short:nonlinear-pointwise}, and \eqref{proof:critical-L4-interpolation},
we obtain
\begin{align}
 &\|U_0(-b)(w(b),u_t(b))-U_0(-a)(w(a),u_t(a))\|_{\mathcal H}\notag\\
 &\qquad\lesssim\|D^2(\mathcal Q\xi')\|_{L^2Y_2(I)}
       +\|V_\xi w+B_{q_\xi}(w)\|_{L^2Y_1(I)}
       +\|f(w)\|_{L^1L^2(I)}\label{short:source-duality}\\
 &\qquad\lesssim\|\xi'\|_{L^2(I)}+\LE_I(w)
       +A^{\frac{5}{6}}S_I(w)^{\frac{7}{6}}+S_I(w)^{\frac{7}{3}}.
 \label{proof:scattering-sources}
\end{align}
Since $A<\infty$, the right-hand side tends to zero as $a\to\infty$,
uniformly for $b>a$. We therefore obtain
\begin{equation*}
 W_+=\lim_{t\to\infty}U_0(-t)(w(t),u_t(t))
 \quad\text{in }\mathcal H.
\end{equation*}
Writing $U_0(t)W_+=(v_L,v_{L,t})$, we have
\begin{equation}\label{short:one-radiation}
 \|(u,u_t)-(q_\xi,0)-(v_L,v_{L,t})\|_{\mathcal H}\longrightarrow0.
\end{equation}
Taking $a=0$ and $b\to\infty$ in \eqref{proof:scattering-sources}, we have
\begin{equation*}
 \|W_+\|_{\mathcal H}\lesssim A+A^2+A^{\frac{7}{3}}.
\end{equation*}
By free energy conservation, the radiation is unique. This proves
Proposition~\ref{presentation:scattering}.
For the constructed solution, we use \eqref{global:stationary-limits}
and weak decay of free waves to obtain
\begin{equation*}
 \vec u(T_{2k})\rightharpoonup(Q_*,0)\ne0.
\end{equation*}
Thus $u$ cannot scatter forward to a free wave.

\section{The fixed-profile obstruction}\label{i5:i5:obstruction}

For the spatial similarities $\mathcal T_{\lambda,a,R}$ of the
introduction, we say that the parameters are pairwise orthogonal if
\begin{equation}\label{obstruction:orthogonality}
 \frac{\lambda_j}{\lambda_k}+\frac{\lambda_k}{\lambda_j}+
 \frac{|a_j-a_k|}{\sqrt{\lambda_j\lambda_k}}\longrightarrow\infty
 \quad(j\ne k).
\end{equation}

\begin{proposition}\label{prop:abstract-obstruction}
Suppose
\begin{equation*}
 U(t)=(q(t),0)+U_0(t)F_++o_{\mathcal H}(1),
\end{equation*}
where $q$ is continuous and precompact in $\dot H^1$. If its omega-limit
set contains nonzero $q_0,q_1$ with
\begin{equation*}
 \operatorname{dist}(q_1,\Sim\cdot q_0)>0,
\end{equation*}
then there exists a sequence $t_n\to\infty$ on which
\eqref{obstruction:competing-decomposition} fails for every fixed
finite family of nonzero energy pairs and every free wave, even when
the active entries vary with $n$ and their number is uniformly bounded.
In particular, no such representation holds at all large times.
The omega-limit set contains a continuum of similarity classes.
\end{proposition}

\subsection{Uniqueness of the radiation}\label{sec:radiation-comparison}

We identify the free state in every competing decomposition by proving,
for fixed energy pairs $F,G$,
\begin{equation}\label{short:radiation-separation}
 \sup_{\lambda,a,R}
 |\langle U_0(t)F,\mathcal T_{\lambda,a,R}G\rangle_{\mathcal H}|
 \longrightarrow0.
\end{equation}
By density and isometry, we take smooth annular Fourier data.
Fourier support restricts $\lambda$ to a compact interval. By the
localized wave estimate \eqref{proof:free-wave-dispersion} and
$L^\infty$--$L^1$ duality in the energy variables, we have
\begin{equation*}
 \sup_{\lambda,a,R}
 |\langle U_0(t)F,\mathcal T_{\lambda,a,R}G\rangle_{\mathcal H}|
 \lesssim_{F,G}|t|^{-2}
\end{equation*}
for these data, proving \eqref{short:radiation-separation}.
By unitarity and a finite net, the convergence is uniform for $G$ in
a compact subset of $\cH$.

We assume a competing decomposition uses a fixed family of nonzero energy
pairs $F_1,\ldots,F_M$ and a free state $\widetilde F_+$, so that
\begin{equation}\label{obstruction:competing-decomposition}
 \begin{gathered}
 U(t)=U_0(t)\widetilde F_+
       +\sum_{j=1}^{J(t)}\mathcal T_{j,t}F_{\nu_j(t)}
       +o_{\cH}(1),\\
 J(t)\le J_*,\qquad \nu_j(t)\in\{1,\ldots,M\}.
 \end{gathered}
\end{equation}
We subtract the decompositions and pair with
$U_0(t)(F_+-\widetilde F_+)$. By compactness of the core,
\eqref{short:radiation-separation}, and free energy conservation, we have
\begin{equation}\label{obstruction:radiation-identity}
 \|F_+-\widetilde F_+\|_{\cH}^2=o(1),\qquad
 F_+=\widetilde F_+.
\end{equation}
The same argument applies if the competing decomposition is only given
along a sequence of times tending to infinity.

\subsection{A continuum of similarity classes}\label{sec:connected-limits}

We use \eqref{master:omega}, with $q(t)$ in place of $q_{\xi(t)}$.
The tail closures in \eqref{master:omega} are nonempty, compact,
connected and nested. Hence $\Omega$ is nonempty, compact and connected.
Since the similarity action is isometric, the function
\begin{equation*}
 q\longmapsto\dist_{\dot H^1}(q,\Sim\cdot q_0)
\end{equation*}
is $1$-Lipschitz and similarity-invariant. Its image on $\Omega$ contains
the interval
\begin{equation}\label{obstruction:distance-interval}
 [0,\dist_{\dot H^1}(q_1,\Sim\cdot q_0)]
 \subseteq
 \{\dist_{\dot H^1}(q,\Sim\cdot q_0)\mid q\in\Omega\}.
\end{equation}
We obtain a continuum of classes, even after removing zero, since
this interval has positive length.

\subsection{Profile collapse and the contradiction}\label{sec:profile-collapse}

\begin{lemma}[Profile collapse]\label{lem:profile-collapse}
Let $G_1,\ldots,G_J$ be nonzero energy pairs, and let
$\mathcal T_{j,n}$ have pairwise orthogonal parameters as in
\eqref{obstruction:orthogonality}. If
\begin{equation*}
 \sum_{j=1}^{J}\mathcal T_{j,n}G_j\longrightarrow C\ne0
 \quad\text{strongly in }\cH,
\end{equation*}
then
\begin{equation}\label{obstruction:profile-collapse}
 J=1,\qquad C\in\Sim\cdot G_1.
\end{equation}
\end{lemma}
\begin{proof}
We pass to a subsequence along which each parameter family converges
or escapes. Compact approximation gives weak convergence to zero for
every escaping profile. Since $C\ne0$, at least one family converges,
and orthogonality permits at most one. We relabel it so that
$C=\mathcal T_*G_1$. By norm decoupling and $G_j\ne0$, we obtain
\begin{equation}\label{obstruction:norm-decoupling}
 \sum_{j=2}^{J}\|G_j\|_{\cH}^2
 =\|C\|_{\cH}^2-\|G_1\|_{\cH}^2=0,\qquad J=1.\qedhere
\end{equation}
\end{proof}

\begin{proof}[Proof of Proposition~\ref{prop:abstract-obstruction}]
For each $k$, we choose a finite $1/k$-net $\{q_{k,i}\}_{i=1}^{m_k}$
in $\Omega$. By the definition of $\Omega$, we may choose times
successively so that
\[
 t_{k,1}>k,\qquad
 t_{k,1}<\cdots<t_{k,m_k}<t_{k+1,1},\qquad
 \|q(t_{k,i})-q_{k,i}\|_{\dot H^1}<1/k .
\]
Listed in increasing order, these times give a sequence $t_n$ for
which the cluster set of $q(t_n)$ equals $\Omega$.
We suppose \eqref{obstruction:competing-decomposition} holds on this sequence.
For each $q\in\Omega\setminus\{0\}$, we pass to a subsequence with
$q(t_n)\to q$ and fixed active indices and number of entries.
By \eqref{obstruction:radiation-identity} and
Lemma~\ref{lem:profile-collapse}, we obtain
\begin{equation}\label{obstruction:finite-orbit-cover}
 \{(q,0)\mid q\in\Omega\setminus\{0\}\}
 \subseteq\bigcup_{m=1}^{M}\Sim\cdot F_m.
\end{equation}
This contradicts the continuum in \eqref{obstruction:distance-interval}.
We have excluded the representation on one sequence, and hence on the
whole half-line.
\end{proof}

\section{Proof of Theorem~\ref*{thm:sequence} on single-core resolution}\label{sec:one-core}
We fix $Q,\ell$ as in Theorem~\ref{thm:sequence}, write $U=u-v$, and set
\[
 \mathcal M=\mathcal M_\ell(Q),\qquad
 d(t)=\dist_{\cH}(\vec U(t),\mathcal M),\qquad
 D_\ell=\ell\cdot\nabla,\qquad \gamma=(1-|\ell|^2)^{-\frac{1}{2}}.
\]
By hypothesis, $d(t_n)\to0$. We use $S(I)$ as in
Section~\ref{master:proof}, and allow constants to depend on $Q$ and $\ell$.

\subsection{The free wave and the residual action}\label{sc:radiation}
By the exterior radiation theorem of Duyckaerts, Kenig and Merle
\cite{DKM19}, a forward-global solution bounded in $\cH$ has a unique
free wave $v_{\rm can}$ such that, for every $A\in\R$,
\begin{equation}\label{sc:canonical}
 \begin{gathered}
 \|\nabla_{t,x}(u-v_{\rm can})(t)\|_{L^2(|x|>t+A)}\longrightarrow0,\\
 U_0(-t)\vec u(t)\rightharpoonup\vec v_{\rm can}(0).
 \end{gathered}
\end{equation}
We choose $p_n\in\mathcal M$ so that $\vec U(t_n)-p_n\to0$ in $\cH$.
Since the normalized family $\{\mathbf Q_\ell(q):q\in\mathcal K_Q\}$
is compact, the uniform estimate in Section~\ref{sec:radiation-comparison}
gives
\[
 U_0(-t_n)p_n\rightharpoonup0,\qquad
 U_0(-t_n)\vec u(t_n)\rightharpoonup\vec v(0).
\]
Comparing the weak limits in \eqref{sc:canonical}, we obtain $v=v_{\rm can}$.

For $F=(w,k)$, we define
\[
 \Mom(F)=\int_{\R^5} k\nabla w\dd x,\qquad
 \mathcal J_\ell(F)=E(F)+\ell\cdot\Mom(F).
\]
We write
\[
 \begin{gathered}
 m=E(Q,0)=\frac{1}{5}\|\nabla Q\|_2^2,\\
 e_c=E(\vec u)-\frac{1}{2}\|\vec v\|_{\cH}^2,\qquad
 p_c=\Mom(\vec u)-\Mom(\vec v).
 \end{gathered}
\]
Since $U_0(-t)\vec U(t)\rightharpoonup0$, the polarized free energy
and momentum of $\vec U(t),\vec v(t)$ tend to zero. By dispersion on
smooth compact data, density and the free Strichartz estimate,
\[
 \omega_v(T):=\sup_{t\ge T}\|v(t)\|_{\frac{10}{3}}\longrightarrow0,\qquad
 s_v(T):=\|v\|_{S([T,\infty))}\longrightarrow0 .
\]
The uniform $\cH$ bound and H\"older's inequality imply
\[
 \left|\int_{\R^5}\bigl(|u|^{\frac{10}{3}}-|U|^{\frac{10}{3}}\bigr)\dd x\right|
 \lesssim\bigl(\|u\|_{\frac{10}{3}}^{\frac{7}{3}}
              +\|U\|_{\frac{10}{3}}^{\frac{7}{3}}\bigr)\omega_v(t)=o(1).
\]
Thus $E(\vec U(t))\to e_c$ and $\Mom(\vec U(t))\to p_c$.

For $q\in\mathcal K_Q$, conformal invariance implies $E(q,0)=m$.
For each constant matrix $M$, stationarity and the decay in
\eqref{sc:decay} imply
\[
 0=\left.\frac{\mathrm d}{\mathrm ds}E(q((I+sM)\cdot),0)\right|_{s=0}
 =\sum_{i,j}M_{ij}\int_{\R^5}\partial_i q\,\partial_j q\dd x
      -m\,\operatorname{tr}M .
\]
Varying $M$ and changing variables by $z=A_\ell y$, we obtain
\[
 \int\nabla q\otimes\nabla q\dd x=mI,\qquad
 \int\nabla\phi_q\otimes\nabla\phi_q\dd y=\frac{m}{\gamma}A_\ell^2.
\]
For $p=\mathbf Q_\ell(q)$, these identities give
\[
 \Mom(p)=-\frac{m}{\gamma}A_\ell^2\ell,\qquad
 \mathcal J_\ell(p)=\frac{1}{\gamma}E(q,0).
\]
Scaling and translation preserve both quantities, so every $p\in\mathcal M$
satisfies
\begin{equation}\label{sc:charges}
 E(p)=\gamma m,\qquad \Mom(p)=-\gamma m\ell,\qquad
 \mathcal J_\ell(p)=\frac{m}{\gamma}.
\end{equation}
Since $E$ and $\Mom$ are uniformly continuous on bounded subsets of
$\cH$, we may pass to the limit along $t_n$ in \eqref{sc:charges}.
The convergence $\vec U(t_n)-p_n\to0$ gives
$(e_c,p_c)=(\gamma m,-\gamma m\ell)$. In particular,
\begin{equation}\label{sc:tails}
 \omega_J(T):=\sup_{t\ge T}
 \left|\mathcal J_\ell(\vec U(t))-\frac{m}{\gamma}\right|\longrightarrow0,
 \qquad
 \eta(T):=\sqrt{\omega_J(T)}+\omega_v(T)+s_v(T)\longrightarrow0 .
\end{equation}

\subsection{Modulation and the linearized equation}\label{sc:geometry}
By elliptic regularity and removability at the stereographic pole
\cite{DKM16}, we have, uniformly for $q\in\mathcal K_Q$,
\begin{equation}\label{sc:decay}
 |\nabla^j q(y)|\lesssim\langle y\rangle^{-3-j},\qquad 0\le j\le3.
\end{equation}
After $k\le2$ differentiations by rotations on $\mathbb S^5$, the same
bound holds for $j+k\le3$. In particular, the orbit is $C^2$ in $\cH$,
and the normalized spatial tangents are in $L^2$.

We translate both components by $\ell t$ and write
$\vec U_c(t)=\vec U(t,\cdot+\ell t)$. For $d(t)$ small, we choose
\begin{equation}\label{sc:gauge}
 \begin{gathered}
 p(t)=\mathcal T_{\lambda(t),b(t)}\mathbf Q_\ell(q(t)),\qquad
 b(t)=a(t)-\ell t,\qquad K_p=T_p\mathcal M,\\
 e(t)=\vec U_c(t)-p(t)\perp_{\cH}K_{p(t)},\qquad
 d(t)\le X(t):=\|e(t)\|_{\cH}\lesssim d(t).
 \end{gathered}
\end{equation}
In local orbit coordinates, the derivative of the orthogonality
conditions at $p$ is the negative Gram matrix of the tangent vectors. We apply the
implicit function theorem modulo the stabilizer. Once we normalize the scale and center, compactness of $\mathcal K_Q$
gives a uniform radius. For uniqueness, suppose that two choices
approach the same data with errors tending to zero. We normalize the
first choice and pass to a nonzero strong limit. Escaping relative scales or centers would send
the second weakly to zero. Otherwise a subsequence of the relative
parameters converges and local uniqueness applies.

For a tangent vector $(Z,-D_\ell Z)$, we test the position against
$-\Delta Z$ and the velocity against $-D_\ell Z$ in \eqref{sc:gauge}.
Since
\[
 -\Delta Z\in L^2,\qquad -D_\ell Z\in\dot H^1,\qquad
 \partial_t\vec U_c\in C_t(L^2\times\dot H^{-1}),
\]
the pairings and $p(t)$ are $C^1$.
We lift $p$ through the constant-rank map
\[
 (\lambda,b,O)\longmapsto\mathcal T_{\lambda,b}\mathbf Q_\ell(OQ),
 \qquad O\in O(6),
\]
with the metric
\[
 \frac{(\mathrm d\lambda)^2+|\mathrm db|^2}{\lambda^2}+|\mathrm dO|^2 .
\]
By compactness after normalization, the differential
has a uniformly bounded inverse on the orthogonal complement of its
kernel. The preceding argument for escaping parameters proves
properness, so the lift extends along each compact time interval.
Writing $q=OQ$, we obtain
\begin{equation}\label{sc:lift}
 \mathrm d\tau=\frac{\mathrm dt}{\lambda(t)},\qquad
 |\lambda'|+|b'|+|O_\tau|\lesssim\|\partial_\tau p\|_{\cH}.
\end{equation}

For $w=(h,k)$ the Hessian of $\mathcal J_\ell$ at $p$ and the
linearized generator in the translating frame are
\begin{equation}\label{sc:linear}
 \begin{aligned}
 \mathfrak b_p(w,w)&=\|k+D_\ell h\|_2^2+
   \langle(-\Delta+D_\ell^2-f'(p_1))h,h\rangle,\\
 \mathcal A_p w&=(k+D_\ell h,\Delta h+f'(p_1)h+D_\ell k).
 \end{aligned}
\end{equation}
The principal form controls $\|w\|_{\cH}^2$, and the potential term
is compact. At $p=\mathbf Q_\ell(q)$, we change variables by $z=A_\ell y$.
Writing $\widetilde h=h\circ A_\ell^{-1}$, we obtain
\[
 \langle(-\Delta+D_\ell^2-f'(p_1))h,h\rangle
 =\gamma^{-1}\langle L_q\widetilde h,\widetilde h\rangle.
\]
By
conformal covariance and the hypothesis on $L_Q$,
\[
 \ker\mathfrak b_p=K_p,\qquad
 \operatorname{ind}_{-}\mathfrak b_p=\operatorname{ind}_{-}L_q\ge1,
 \qquad \langle L_qq,q\rangle=-\frac{4}{3}\|\nabla q\|_2^2 .
\]
Let $Y_i$ be an orthonormal basis of the negative spectral subspace
of $L_q$, with $L_qY_i=-\kappa_i^2Y_i$. The Lorentz transform of
$e^{\kappa_i s}Y_i$ has data
\[
 V_i=(\psi_i,\sigma_i\psi_i-D_\ell\psi_i),\qquad
 \psi_i(y)=e^{-\gamma\kappa_i\ell\cdot y}Y_i(A_\ell y),\qquad
 \sigma_i=\frac{\kappa_i}{\gamma}.
\]
Thus $\mathcal A_p V_i=\sigma_iV_i$. Elliptic decay of $Y_i$ at any
rate between $|\ell|\kappa_i$ and $\kappa_i$ implies exponential
decay of $V_i$ and its first derivatives. We put
\[
 G_p=\operatorname{span}\{V_i\},\qquad
 (\beta_p(w))_i=\mathfrak b_p(V_i,w)
       =\sigma_i\int_{\R^5}\bigl(\psi_i k-(V_i)_2h\bigr)\dd y.
\]
Since $\mathcal A_p$ is skew for $\mathfrak b_p$, we obtain
\begin{equation}\label{sc:dual}
 \mathfrak b_p|_{G_p\times G_p}=0,\qquad
 \beta_p\mathcal A_p=-\operatorname{diag}(\sigma_i)\beta_p,\qquad
 \beta_p|_{K_p}=0 .
\end{equation}

\subsection{Exterior energy and coercivity}\label{sc:exterior}
At a fixed time, we normalize the scale and center. For
\[
 p\in\mathcal P=\{\mathcal T_{\alpha,c}\mathbf Q_\ell(q):
       \frac{1}{2}\le\alpha\le2,\ |c|\le1,\ q\in\mathcal K_Q\},
\]
we set
\[
 P_p(s,y)=p_1(y-\ell s),\qquad
 \chi=\mathbf1_{\{s\ge0,\,|y|>s+R\}}.
\]
By \eqref{sc:decay}, we have
\[
 \|\chi P_p\|_S\lesssim R^{-\frac{3}{2}},\qquad
 \|\chi f'(P_p)\|_{L_s^{\frac{7}{4}}L_y^{\frac{7}{2}}}\lesssim R^{-2}.
\]
We fix $R$ large and solve the exterior linear equation
\[
 \Box h_w=\chi f'(P_p)h_w,\qquad \vec h_w(0)=w,\qquad
 \rho_p(w)=\limsup_{s\to\infty}
         \|\nabla_{s,y}h_w(s)\|_{L^2(|y|>s+R)} .
\]
We absorb the potential in the free Strichartz estimate, first for
$h_w$ and then for the difference at $p,p'$. We obtain
\begin{align*}
 \|\vec h_w\|_{L^\infty_s\cH}+\|h_w\|_S&\lesssim\|w\|_{\cH},\\
 |\rho_p(w)-\rho_{p'}(w)|
 &\lesssim\|\chi(f'(P_p)-f'(P_{p'}))\|_{L_s^{\frac{7}{4}}L_y^{\frac{7}{2}}}
                   \|w\|_{\cH}.
\end{align*}
Thus $\rho_p$ is a bounded seminorm, continuous in $p$ uniformly on
the energy unit ball.

For $p=\mathbf Q_\ell(q)$ we claim that
\begin{equation}\label{sc:detection}
 \rho_p|_{K_p}=0,\qquad \rho_p(g)>0\quad(0\ne g\in G_p).
\end{equation}
By finite propagation, we may use the full linearized solution on the
cone. For $(Z,-D_\ell Z)\in K_p$, this solution is $Z(y-\ell s)$,
whose exterior energy tends to zero by \eqref{sc:decay}.
For $0\ne g\in G_p$, we group the rest-frame eigenfunctions with the
same eigenvalue. The corresponding solution has the form
\[
 H^0(s,z)=\sum_\kappa e^{\kappa s}Y^{(\kappa)}(z).
\]
By \cite[Proposition~3.9]{DKM16},
\[
 r^2 e^{\kappa r}Y^{(\kappa)}(r\omega)\longrightarrow G_\kappa(\omega)
 \quad\text{in }L^2(\mathbb S^4),\qquad
 G_\kappa\ne0\quad\text{if }Y^{(\kappa)}\ne0 .
\]
With $|z|=s+h$, polar coordinates and independence of the
exponentials on $[M,M+1]$ imply
\[
 \lim_{s\to\infty}\int_{s+M<|z|<s+M+1}|\partial_sH^0|^2\dd z
 =\int_M^{M+1}\Big\|\sum_\kappa\kappa e^{-\kappa h}G_\kappa
             \Big\|_{L^2(\mathbb S^4)}^2\dd h=:K_M>0 .
\]
Let $H$ be the Lorentz transform of $H^0$. For fixed $M$ and
$L<s<2L$ with $L$ large, the image of the shell satisfies
\[
 c_\ell L<t<C_\ell L,\qquad
 |x|-t=\frac{|z|^2-s^2}{|x|+t}\gtrsim_\ell M .
\]
We choose $M$ so that $|x|>t+R$. The transformation preserves
spacetime volume, and $|\partial_s H^0|\lesssim_\ell|\nabla_{t,x}H|$.
Averaging in $s$ and using finite propagation, we obtain
\[
 0<K_M\lesssim_\ell\limsup_{L\to\infty}
 \frac{1}{L}\int_{c_\ell L}^{C_\ell L}\int_{|x|>t+R}
            |\nabla_{t,x}H|^2\dd x\dd t
 \lesssim_\ell\rho_p(g)^2 .
\]
This proves \eqref{sc:detection}.

At $(t_0,x_0)$ and scale $\lambda_0$, we set
\[
 (t,x)=(t_0+\lambda_0s,x_0+\lambda_0y),\qquad
 (\widetilde u,\widetilde v)(s,y)=\lambda_0^{\frac{3}{2}}(u,v)(t,x).
\]
Suppose that the rescaled data satisfy
\[
 \vec{\widetilde u}(0)-\vec{\widetilde v}(0)=p+w,\qquad p\in\mathcal P.
\]
For the exterior remainder $z$, we solve
\[
 \Box z=\chi\{f(P_p+z+\widetilde v)-f(P_p)\},\qquad \vec z(0)=w .
\]
For the Taylor remainder, we have
\begin{align*}
 N_a(b)&=f(a+b)-f(a)-f'(a)b,\\
 |N_a(b)|&\lesssim |a|^{\frac{1}{3}}|b|^2+|b|^{\frac{7}{3}}.
\end{align*}
We set $X=\|w\|_{\cH}$. Since $\|\widetilde v\|_S=s_v(t_0)$, the
Strichartz estimate and contraction for small $X+s_v(t_0)$ imply
\begin{align*}
 \|\vec z\|_{L_s^\infty\cH}+\|z\|_S&\lesssim X+s_v(t_0),\\
 \|\vec z-\vec h_w\|_{L_s^\infty\cH}
 &\lesssim\|\chi(N_{P_p}(z+\widetilde v)+f'(P_p)\widetilde v)\|_{L^1_sL_y^2}
 \lesssim X^2+s_v(t_0).
\end{align*}
By finite propagation and the change of variables above,
\[
 |y|>s+R\quad\Longrightarrow\quad
 z=\widetilde u-\widetilde v-P_p,\qquad
 |x|>t+\lambda_0R-t_0-|x_0|.
\]
Thus \eqref{sc:canonical} and the exterior decay of $P_p$ give
\begin{equation}\label{sc:exterior-bound}
 \lim_{s\to\infty}\|\nabla_{s,y}z(s)\|_{L^2(|y|>s+R)}=0,\qquad
 \rho_p(w)\lesssim X^2+s_v(t_0).
\end{equation}

We fix a normalized profile $p_j=\mathbf Q_\ell(q_j)$ and use the
subscript $j$ for the corresponding form, spaces and functionals.
Let $\Pi_{K_j}$ be the energy-orthogonal projection onto $K_j$. We claim
\begin{equation}\label{sc:coerc}
 \|w\|_{\cH}\lesssim_j|\beta_j(w)|+\rho_j(w)
              +\sqrt{[\mathfrak b_j(w,w)]_+}+\|\Pi_{K_j}w\|_{\cH}.
\end{equation}
Let $\pi_j$ be the $\mathfrak b_j$-orthogonal projection onto a
maximal negative subspace. On $\ker\pi_j$ the form is nonnegative
with kernel $K_j$. By \eqref{sc:dual}, $G_j$ is isotropic. Moreover,
\[
 \dim G_j=\operatorname{rank}\pi_j,\qquad
 G_j\cap\ker\pi_j\subset G_j\cap K_j=\{0\}.
\]
For $w\in\ker\beta_j$, we can therefore choose $g\in G_j$ with
$\pi_jg=\pi_jw$. Then
\[
 \mathfrak b_j(w,w)=\mathfrak b_j(w-g,w-g)\ge0,\qquad
 \mathfrak b_j(w,w)=0\quad\Longleftrightarrow\quad w-g\in K_j .
\]
If \eqref{sc:coerc} fails, we take unit vectors $w_n$ with right-hand
side tending to zero and pass to a weak limit $w$. The principal form
and the bounded seminorm $\rho_j$ are weakly lower semicontinuous,
and the potential term is compact. Hence
\[
 \beta_j(w)=0,\qquad \mathfrak b_j(w,w)\le0,\qquad
 \rho_j(w)=0,\qquad w\perp_{\cH}K_j .
\]
The equality case puts $w$ in $G_j+K_j$, and \eqref{sc:detection}
then gives $w=0$. By compactness of the potential term, we would have
\[
 1=\|w_n\|_{\cH}^2\lesssim[\mathfrak b_j(w_n,w_n)]_++o(1)\longrightarrow0,
\]
which proves \eqref{sc:coerc}.

For $p$ close to $p_j$ and $w\perp_{\cH}K_p$, we have by continuity
\[
 |\rho_p(w)-\rho_j(w)|
 +|\mathfrak b_p(w,w)-\mathfrak b_j(w,w)|^{\frac{1}{2}}
 +\|\Pi_{K_j}w\|_{\cH}\le o_{p\to p_j}(1)X .
\]
At the normalized data $p+w$ at time $t$, we expand the action and
use \eqref{sc:charges} and \eqref{sc:exterior-bound} to obtain
\[
 \left|\mathcal J_\ell(\vec U(t))-\frac{m}{\gamma}
            -\frac{1}{2}\mathfrak b_p(w,w)\right|\lesssim X^3,\qquad
 \rho_p(w)\lesssim X^2+s_v(t).
\]
We absorb $o(1)X+O(X^{\frac{3}{2}}+X^2)$ in \eqref{sc:coerc} and
fix $C_j$ such that
\begin{equation}\label{sc:algebra}
 X\le C_j\bigl(|\beta_j(w)|+\sqrt{\omega_J(t)}+s_v(t)\bigr).
\end{equation}

\subsection{Control at all later times}\label{sc:dynamics}
We choose $T=t_n$ with $d(T)+\eta(T)$ sufficiently small and work
on the maximal interval starting at $T$ on which \eqref{sc:gauge}
holds with $X$ small. We write $e=(h,k)$ and $v_c(t,y)=v(t,y+\ell t)$.
For fixed $t$, $p(t)$ is stationary in the translating frame, so
\begin{equation}\label{sc:error}
 \begin{aligned}
 \partial_t e&=\mathcal A_{p(t)}e+(0,F)-\partial_t p,\\
 F&=N_{p_1}(h)+f(p_1+h+v_c)-f(p_1+h),\qquad
 \|F\|_{\frac{10}{7}}\lesssim X^2+\omega_v(t).
 \end{aligned}
\end{equation}
The bound on the source follows from H\"older's inequality and Sobolev
embedding. At a fixed time,
we normalize the scale and center and denote time in those units by
$s$. For a local tangent basis $Z_i(p)$, we differentiate the
orthogonality in \eqref{sc:gauge}. By \eqref{sc:decay} and integration
by parts,
\begin{align*}
 \langle\partial_s p,Z_i(p)\rangle_{\cH}
 &=\langle\mathcal A_p e+(0,F),Z_i(p)\rangle_{\cH}
      +\langle e,\mathrm d Z_i(p)[\partial_s p]\rangle_{\cH},\\
 |\langle\partial_s p,Z_i(p)\rangle_{\cH}|
 &\lesssim X+\omega_v+X\|\partial_s p\|_{\cH}.
\end{align*}
We invert the Gram matrix and absorb the last term. Since
$\partial_s=\partial_\tau$ at the normalization time, \eqref{sc:lift}
implies
\begin{equation}\label{sc:speed}
 |\lambda'|+|b'|+|O_\tau|\lesssim
 \|\partial_\tau p\|_{\cH}\lesssim X+\omega_v .
\end{equation}

At $\tau_0=\tau(t_0)$, we fix $\mathcal T_0=\mathcal T_{\lambda_0,b_0}$
at the scale and center of $p(t_0)$ and set
\[
 w=\mathcal T_0^{-1}e,\qquad p_*=\mathcal T_0^{-1}p,\qquad
 \alpha=\frac{\lambda}{\lambda_0},\qquad
 \tau(t)=\int_T^t\frac{\dd s}{\lambda(s)} .
\]
Here $w\perp_{\cH}K_{p_*}$. In these coordinates \eqref{sc:error} reads
\[
 \begin{gathered}
 \partial_\tau w=\alpha(\mathcal A_{p_*}w+(0,F_*))-\partial_\tau p_*,\\
 F_*(t,y)=\lambda_0^{\frac{7}{2}}F(t,b_0+\lambda_0y),\qquad
 \|F_*\|_{\frac{10}{7}}=\|F\|_{\frac{10}{7}} .
 \end{gathered}
\]
For $p_*\in\mathcal P$ close to $p_j$, we apply \eqref{sc:algebra}
at scale $\lambda_0$ and laboratory center $b_0+\ell t$.
With $\Gamma_j=\operatorname{diag}(\sigma_{j,i})$ and
$B_j=|\beta_j(w)|$, we have by \eqref{sc:dual}
\begin{align*}
 \partial_\tau\beta_j(w)+\alpha\Gamma_j\beta_j(w)
 &=\alpha\beta_j((\mathcal A_{p_*}-\mathcal A_{p_j})w)
   +\alpha\beta_j(0,F_*)-\beta_j(\partial_\tau p_*),\\
 \|\mathcal A_{p_*}-\mathcal A_{p_j}\|_{\cH\to\cH}
 &\lesssim\|f'((p_*)_1)-f'((p_j)_1)\|_5=o_{p_*\to p_j}(1),\\
 |\beta_j(\partial_\tau p_*)|
 &\lesssim_j\|\Pi_{K_{p_*}}-\Pi_{K_j}\|(X+\omega_v).
\end{align*}
For the last bound, we use $\partial_\tau p_*\in K_{p_*}$ and
\eqref{sc:speed}. Exponential decay of the modes gives the other
forcing bound,
\[
 |\beta_j(0,F_*)|\lesssim_j\|F_*\|_{10/7}.
\]
We put $\sigma_j=\min_i\sigma_{j,i}$. For
$\frac{1}{2}\le\alpha\le2$, we choose $r_j$ and then $X$ small
enough that on $B(p_j,r_j)$,
\[
 |\partial_\tau\beta_j(w)+\alpha\Gamma_j\beta_j(w)|
 \le\varepsilon_jX+C'_j(X^2+\omega_v),\qquad
 C_j(\varepsilon_j+C'_jX)\le\frac{\sigma_j}{4}.
\]
Finitely many $B(p_j,\frac{r_j}{2})$ cover the normalized profiles.
On the corresponding full balls, \eqref{sc:algebra} therefore implies
\[
 \partial_\tau B_j+cB_j\lesssim\eta(t)\quad\text{a.e.},\qquad
 c=\frac{1}{4}\min_j\sigma_j>0 .
\]
Writing $t(r)$ for the inverse of $\tau$, we set
\[
 \mathsf X(r)=X(t(r)),\qquad g(r)=\eta(t(r)).
\]
We integrate from $\tau_0$, use $B_j(\tau_0)\lesssim_j\mathsf X(\tau_0)$,
and apply \eqref{sc:algebra} once more to obtain
\begin{equation}\label{sc:window}
 \mathsf X(r)\lesssim e^{-c(r-\tau_0)}\mathsf X(\tau_0)+g(r)
       +\int_{\tau_0}^r e^{-c(r-s)}g(s)\dd s .
\end{equation}

We fix $L\ge1$ so that the coefficient of $\mathsf X(\tau_0)$ at
$r=\tau_0+L$ is at most $\frac{1}{2}$. For $I=[\tau_0,r]$ of length
at most $L$, we integrate \eqref{sc:speed}, first for $\log\alpha$, to
obtain
\[
 |\log\alpha|+\frac{|b-b_0|}{\lambda_0}
    +\|p_*(r)-p_*(\tau_0)\|_{\cH}
 \lesssim_L\sup_{s\in I}\bigl(\mathsf X(s)+\omega_v(t(s))\bigr).
\]
If $p_*(\tau_0)\in B(p_j,\frac{r_j}{2})$, smallness of this supremum
keeps $p_*$ in $B(p_j,r_j)\cap\mathcal P$. Thus \eqref{sc:window}
holds on $I$. Since $g$ is nonincreasing, at the endpoints we have
\[
 \mathsf X((n+1)L)\le\frac{1}{2}\mathsf X(nL)
              +C_L\int_{nL}^{(n+1)L}g(s)\dd s .
\]
We iterate this inequality, then use \eqref{sc:window} on the last
partial interval. For $0<c_1\le\min(c,\frac{\log2}{L})$, the factors
$2^{-(n-1-k)}$ on $[kL,(k+1)L]$ are bounded by
$2e^{-c_1(nL-s)}$. Returning to $t$ and using
$\mathsf X(0)\lesssim d(T)$, we obtain
\begin{equation}\label{sc:memory}
 \begin{aligned}
 X(t)\lesssim{}& e^{-c_1\tau(t)}d(T)+\eta(t)\\
 &+\int_T^t e^{-c_1(\tau(t)-\tau(s))}
                \eta(s)\frac{\dd s}{\lambda(s)} .
 \end{aligned}
\end{equation}
In particular, $X(t)\lesssim d(T)+\eta(T)$. For $T=t_n$ sufficiently
late, this improves the bound defining the maximal interval. At a finite
endpoint $t_*$, $\vec U_c(t)$ would have a strong limit. If
$(\lambda(s_n),b(s_n))$ escaped every compact subset of
$(0,\infty)\times\R^5$ along $s_n\uparrow t_*$, compactness of the
normalized shapes would give $p(s_n)\rightharpoonup0$ along a subsequence.
We would then have
\[
 \liminf_{n\to\infty}\|e(s_n)\|_{\cH}^2
 \ge c_Q^2,\qquad
 c_Q=\min_{q\in\mathcal K_Q}\|\mathbf Q_\ell(q)\|_{\cH}>0 .
\]
Choosing $T$ so that $X<c_Q/2$ excludes escape. Thus $p(t)$ has a strong
limit along a subsequence, and the local chart in \eqref{sc:gauge}
extends the modulation past $t_*$. Hence the interval is $[T,\infty)$.

By \eqref{sc:speed}, $\lambda(t)\lesssim\lambda(T)+t-T$, so
$\tau(t)\to\infty$. By \eqref{sc:tails}, $g(r)\to0$. Splitting the
integral at $\frac{r}{2}$, we have
\[
 \int_0^r e^{-c_1(r-s)}g(s)\dd s
 \lesssim e^{-\frac{c_1r}{2}}g(0)+g\!\left(\frac{r}{2}\right)\longrightarrow0 .
\]
Equations \eqref{sc:memory} and \eqref{sc:speed}, with $q=OQ$, now imply
\[
 X(t)+|\lambda'(t)|+|a'(t)-\ell|
      +\lambda(t)\|q'(t)\|_{\dot H^1\cap L^2}\longrightarrow0 .
\]
Integrating the bounds on $\lambda'$ and $a'-\ell$, we obtain
\[
 \lambda(t)=o(t),\qquad a(t)=\ell t+o(t),
\]
which completes the proof of Theorem~\ref{thm:sequence}.

The same estimate \eqref{sc:memory} holds if $v$ is the canonical
free wave, $e_c+\ell\cdot p_c=\frac{m}{\gamma}$, and
$d(T)+\eta(T)$ is sufficiently small at one time $T$. Indeed, these conditions imply \eqref{sc:tails} and suffice to start the argument
at $T$.

For the ground state
\[
 W(x)=\left(1+\frac{|x|^2}{15}\right)^{-3/2},
\]
a fixed dilation makes $W$ the stereographic pullback of a constant.
At either sign, the spherical linearization is $-\Delta_{\mathbb S^5}-5$,
whose kernel consists of the first harmonics. These correspond to translations
and dilation, so we may take the normalized shape fixed.
For a general nondegenerate $Q$ with $v=0$, conservation and
\eqref{sc:charges} imply $\eta=0$. In this case
\[
 X+|\partial_\tau\log\lambda|+|O_\tau|+\frac{|b_\tau|}{\lambda}
       \lesssim e^{-c_1\tau}.
\]
Integrating first for $\log\lambda$ and $O$, we obtain
$\lambda_\infty>0$, $O_\infty\in O(6)$ and $\tau(t)\asymp t$ at
large times. Integrating for $b$ and using the $C^1$ dependence on
the parameters, we find $a_\infty\in\R^5$ and $c_*>0$ with
\[
 \big\|\vec u(t)-\mathcal T_{\lambda_\infty,a_\infty+\ell t}
                      \mathbf Q_\ell(O_\infty Q)\big\|_{\cH}
 \lesssim e^{-c_*t}.
\]

\appendix

\section{Energy and nonzero radiation}\label{app:energy}

We write $\vec v_L(t)=U_0(t)W_+$. The decomposition
\eqref{short:one-radiation} gives
\begin{equation}\label{proof:energy-decomposition}
 \begin{gathered}
 E(u,u_t)=E(Q_*,0)+\frac{1}{2}\|W_+\|_{\mathcal H}^2,\\
 \Mom(u,u_t)=\Mom(\vec v_L),\qquad \Mom(f,g)=\int g\nabla f.
 \end{gathered}
\end{equation}
By \eqref{short:radiation-separation} and compactness of the core,
the quadratic cross terms vanish. By free dispersion, density, and
stationarity, we have
\begin{equation*}
 \|v_L(t)\|_{\frac{10}{3}}\longrightarrow0,\qquad
 E(q_\xi,0)=E(Q_*,0).
\end{equation*}
We obtain \eqref{proof:energy-decomposition} by H\"older and conservation.

We write
\begin{equation*}
 \mathcal E_L(v)=\frac{1}{2}
       \bigl(\|v_t\|_2^2+\langle Lv,v\rangle\bigr).
\end{equation*}
Let $U_z$ and $U_a$ be the incoming free waves for $v_z$ and $v_a$.
By \eqref{short:packet-duhamel}, the disjoint bands, the positive
amplitudes, and Riemann--Lebesgue, we have
\begin{gather*}
 \|\vec U_z(0)\|_{\mathcal H}\asymp1,\qquad z\in\mathcal C,\\
 \sup_{z,z'\in\mathcal C}
 |\langle\vec U_z(s),\vec U_{z'}(0)\rangle_{\mathcal H}|
 \longrightarrow0\qquad(|s|\to\infty).
\end{gather*}
The second limit is uniform because the packet family is finite and
the amplitudes are bounded. All distinct centers $2T_j,3T_j$ are
separated by at least $T_0$. Since
\begin{equation*}
 \sum_{j<N}\eta_j\le2\varepsilon,\qquad
 \varepsilon^2\le\sum_{j<N}\eta_j^2\le\frac{4}{3}\varepsilon^2,
\end{equation*}
we obtain
\begin{equation}\label{proof:train-energy-lower-bound}
 \begin{aligned}
 2\mathcal E_L(v_a)&=\|\vec U_a(0)\|_{\mathcal H}^2\\
 &=\sum_{j<N}\eta_j^2
   \bigl(\|\vec U_{x_j}(0)\|_{\mathcal H}^2
        +\|\vec U_{y_j}(0)\|_{\mathcal H}^2\bigr)
   +o_{T_0\to\infty}(1)\varepsilon^2
 \asymp\varepsilon^2.
 \end{aligned}
\end{equation}
The first equality follows from incoming scattering and local decay of
the potential pairing. The error is uniform in $N$ and the controls.
Since $E'(Q_*,0)=0$ and $\varepsilon=d^3$, we may apply Taylor's formula
to the data in \eqref{short:physical-data}. Passing to the strong
initial-data limit and using \eqref{proof:energy-decomposition}, we obtain
\begin{equation}\label{proof:radiation-size}
 \begin{aligned}
 E(\vec u_N(0))-E(Q_*,0)
 &=\mathcal E_L(v_a)+O(d\varepsilon^2)\asymp\varepsilon^2,\\
 \frac{1}{2}\|W_+\|_{\mathcal H}^2
 &=E(u,u_t)-E(Q_*,0)\asymp\varepsilon^2>0.
 \end{aligned}
\end{equation}
\begin{remark}
The radiation-subtracted curve is precompact, but it is not in general
an exact NLW trajectory. Indeed,
\begin{equation*}
 (\partial_t^2-\Delta)(u-v_L)=f(u),
 \qquad f(u)\ne f(u-v_L)\quad\text{in general}.
\end{equation*}
We also have
\begin{equation*}
 \sum_j\eta_j^2<\infty,\qquad
 \sum_jT_j\eta_j^2=\sum_jd=\infty.
\end{equation*}
Thus finite energy alone does not verify the spatial-moment
hypotheses for weighted conformal charges.
\end{remark}

\section{Tests at zero frequency}\label{app:zero-frequency}

For a kernel column with nonzero monopole, the following tests
explain the use of $P$ on energy positions and $zT_z$ in
\eqref{proof:low-frequency-maps}.

We assume $m(\nu_i)\ne0$. Since $\nu_i=O(\rho^{-7})\in L^1$
and $K_i=D^{-2}\nu_i$, we have
\begin{equation}\label{test:energy-position}
 \widehat{D^{-1}K_i}(\zeta)
 =\frac{\widehat{\nu_i}(\zeta)}{|\zeta|^3}\notin L^2(\R^5),\qquad
 \widehat{\nu_i}(0)\ne0.
\end{equation}
Indeed, $\widehat{\nu_i}$ is continuous and nonzero at the origin,
so for some $\delta>0$,
\[
 \int_{|\zeta|<\delta}\frac{|\widehat{\nu_i}(\zeta)|^2}{|\zeta|^6}\,d\zeta
 \gtrsim\int_0^\delta r^{-2}\,dr=\infty.
\]
By $\dot H^1$ duality, $K_i^*$ is unbounded on energy positions.
We therefore use $P$ on positions and reserve $P_0=I-KK^*$ for $L^2$
velocities.

For \eqref{short:resolvent-kernels}, we rescale and apply dominated
convergence with $\nu_i\in L^1$. Uniformly on a fixed annulus, we have
\begin{equation}\label{test:rescaled-newton}
 \lambda^{-1}[(H_\lambda-H_0)\nu_i]\left(\frac{X}{\lambda}\right)
 \longrightarrow
 \frac{m(\nu_i)[h(|X|)-h(0)]}{8\pi^2|X|}
 \qquad(\lambda\downarrow0).
\end{equation}
On an annulus where the limit is nonzero, the amplitude, the squared
weight and the volume contribute $\lambda^2$, $\lambda^4$ and
$\lambda^{-5}$, respectively. Since their product is $\lambda$, we obtain
\begin{equation}\label{test:newton-lipschitz}
 \|(H_\lambda-H_0)\nu_i\|_{X_2}\gtrsim\lambda^{\frac{1}{2}},
 \qquad
 \|P(H_\lambda-H_0)\nu_i\|_{X_2}\gtrsim\lambda^{\frac{1}{2}}.
\end{equation}
The second bound follows from boundedness of the finite-rank
coefficients and
\begin{equation*}
 K(\frac{X}{\lambda})=O(\lambda^3)=o(\lambda).
\end{equation*}
Thus $T_z=PH_z\nu$ need not be Lipschitz into $X_2$, which explains
the use of $zT_z$ and auxiliary free waves in \eqref{short:low-matrix}.

\section{An analytical spectral bound for the single-node state}
\label{app:single-node}

We construct the profile in \eqref{np:profile-r} by comparison.
For the upper bound, we adapt the positive comparison in \cite[Section~8]{NP26}.
For the lower bound, we evaluate the Rayleigh quotient at an explicit test function.

\begin{lemma}[Analytical profile and spectral gap]\label{np:certificate}
There is a regular profile satisfying \eqref{np:profile-r}. For this
profile, the one-sided problem
\begin{equation}\label{np:one-sided}
 h''+2\cot r\,h'
 +\left(\mathfrak V(w)-\frac{\Lambda}{4\sin^2(\frac{r}{2})}\right)h=0
\end{equation}
has exactly one spectral value $\Lambda_*>2$ in the regular endpoint
class, and
\begin{equation}\label{np:principal-bracket}
 20<\Lambda_*<30.
\end{equation}
For $\Lambda=\alpha(\alpha+1)$ with $\alpha>0$, the regular class means the
$r^\alpha$ Frobenius branch at zero and the bounded even branch at $\pi$.
\end{lemma}

\subsection{The profile and two comparisons}

On $[0,1]$, we set
\begin{equation*}
 \mu=\sqrt{1-x^2},\qquad
 \mathcal L=-\mu^2\partial_x^2+3x\partial_x+\frac{15}{16},
\end{equation*}
and define $G=\mathcal L^{-1}$ by the boundary conditions
$(Gf)(0)=0$ and regularity at $1$. The maximum principle gives
positivity. To prove that $G:C[0,1]\to C^1[0,1]$ is bounded and
hence compact as a map into $C[0,1]$, we use
\begin{equation*}
 \|Gf\|_\infty\le\frac{16}{15}\|f\|_\infty,\qquad
 \mu(x)^3(Gf)'(x)=\int_x^1\mu(t)(f(t)-\frac{15}{16}Gf(t))\,dt.
\end{equation*}

For $f,g\in C[0,1]$ positive on $(0,1]$, we also use the ratio property
\begin{equation}\label{an:ratio}
 \frac{f}{g}\text{ nondecreasing}\quad\Longrightarrow\quad
 \frac{Gf}{Gg}\text{ nondecreasing}.
\end{equation}
Indeed, for $u=Gf$ and $v=Gg$, we have
\begin{equation*}
 \left[\mu^3v^2\left(\frac{u}{v}\right)'\right]'
 =\mu vg\left(\frac{u}{v}-\frac{f}{g}\right),\qquad
 \left.\mu^3v^2\left(\frac{u}{v}\right)'\right|_{x=0,1}=0.
\end{equation*}
On a component where $(\frac{u}{v})'<0$, the difference
$\frac{u}{v}-\frac{f}{g}$ strictly decreases. The flux cannot then have
a negative interior minimum, contrary to its zero endpoint values.
Interchanging $f,g$ proves the corresponding assertion for nonincreasing ratios.

We choose
\begin{equation}\label{an:comparisons}
 \begin{gathered}
 D(x)=\left(1-\frac{x^2}{4}\right)\left(1-\frac{x^2}{8}\right),
 \qquad b_+(x)=\frac{x}{D(x)^{\frac{3}{2}}},\\
 b_-(x)=\frac{x}{\left[(1-\frac{5x^2}{16})(1-\frac{x^2}{8})\right]^{\frac{3}{2}}}.
 \end{gathered}
\end{equation}
Both denominators are positive and decreasing. On $(0,1]$, we compute
\begin{align*}
 4x^{\frac{7}{3}}\left(\frac{\mathcal Lb_+}{b_+^{\frac{7}{3}}}\right)'
 &=-3\left(1-\frac{7x^2}{4}\right)^2-\frac{x^2}{8}
   -\frac{5x^4}{32}-\frac{x^6(980-425x^2)}{1024}<0,\\
 \frac{65536}{5}x^{\frac{7}{3}}
 \left(\frac{\mathcal Lb_-}{b_-^{\frac{7}{3}}}\right)'
 &=x^2\bigl[309+(1-x^2)(38603+2635x^2-2125x^4)\bigr]>0.
\end{align*}
Since the first ratio is positive at $1$ and the second tends to $0$ at $0$, we have
$\mathcal Lb_+,\mathcal Lb_->0$ on $(0,1]$. Moreover,
\begin{equation*}
 \left(\log\frac{b_-}{b_+}\right)'
 =\frac{3x}{16(1-\frac{x^2}{4})(1-\frac{5x^2}{16})}>0.
\end{equation*}

We consider the closed convex set
\begin{equation*}
 \mathcal K=\left\{v\in C[0,1]:v\ge0,\ v(1)=1,\
 \frac{v}{b_+}\text{ nondecreasing},\
 \frac{v}{b_-}\text{ nonincreasing}\right\}.
\end{equation*}
The ratios are taken on $(0,1]$. We have
\begin{equation*}
 \frac{b_+}{b_+(1)}\in\mathcal K,\qquad
 \frac{b_-}{b_-(1)}\le v\le\frac{b_+}{b_+(1)}\quad(v\in\mathcal K).
\end{equation*}
The map
\begin{equation*}
 T(v)=\frac{G[v^{\frac{7}{3}}]}{G[v^{\frac{7}{3}}](1)}
\end{equation*}
is continuous and compact on $\mathcal K$, since
\begin{equation*}
 G[v^{\frac{7}{3}}](1)
 \ge G\left[\left(\frac{b_-}{b_-(1)}\right)^{\frac{7}{3}}\right](1)>0.
\end{equation*}
For $v\in\mathcal K$, we write
\begin{equation*}
 \frac{v^{\frac{7}{3}}}{\mathcal Lb_\pm}
 =\left(\frac{v}{b_\pm}\right)^{\frac{7}{3}}
   \frac{b_\pm^{\frac{7}{3}}}{\mathcal Lb_\pm}.
\end{equation*}
The product is nondecreasing for $b_+$ and nonincreasing for $b_-$.
It follows from \eqref{an:ratio} and $G[\mathcal Lb_\pm]=b_\pm$ that
$T(v)\in\mathcal K$. By Schauder's theorem, $T$ has a fixed point $v$. We set
\begin{equation}\label{an:W}
 W=\left(\frac{15}{16}G[v^{\frac{7}{3}}](1)\right)^{-\frac{3}{4}}v,\qquad
 \mathcal LW=\frac{15}{16}W^{\frac{7}{3}},\qquad W(0)=0.
\end{equation}
The function $W$ is regular, and the ratio bounds ensure $W'>0$ on $[0,1]$.
Its odd extension defines $w(r)=W(\cos r)$ satisfying
\eqref{np:profile-r}.

By the ratio bounds and the regular endpoint equation, we have
\begin{equation}\label{an:amplitude}
 \begin{gathered}
 \frac{W}{b_+}\text{ nondecreasing},\qquad
 \frac{W}{b_-}\text{ nonincreasing},\\
 W(1)^{\frac{4}{3}}
 =1+\frac{16}{5}\frac{W'(1)}{W(1)}
 \le1+\frac{16}{5}\frac{b_-'(1)}{b_-(1)}=\frac{765}{77}<10,\\
 W(x)\le W(1)x\left(\frac{\frac{21}{32}}{D(x)}\right)^{\frac{3}{2}}.
 \end{gathered}
\end{equation}

\subsection{The radial kernel}

We use a Wronskian comparison as in \cite[Section~4]{NP26}.

\begin{lemma}\label{np:radial-kernel}
The only regular solution of \eqref{np:angular-equation} with
$(j,\ell)=(0,0)$ is zero.
\end{lemma}
\begin{proof}
We work on $0\le x\le1$ and write
\begin{equation*}
 \mathcal Dh=(1-x^2)h''-3xh'+\mathfrak V(W)h.
\end{equation*}
We normalize the regular solution of $\mathcal Dh=0$ by $h(1)=1$.
We show that $h(0)<0<h'(0)$.

For $0<x\le1$, we set
\begin{equation*}
 a=\frac{W'}{W},\qquad
 B=(1-x^2)a^2-2xa+\frac{15}{16}(W^{\frac{4}{3}}-1).
\end{equation*}
By \eqref{an:W}, we have
\begin{equation*}
 (1-x^2)a'=xa-B,\qquad B(1)=a(1)>0.
\end{equation*}
At any zero of $B$ in $(0,1)$, we compute
\begin{equation*}
 B'=-\frac{a\left(16[(1-x^2)a-x]^2+9-x^2\right)}{12(1-x^2)}<0.
\end{equation*}
Thus $B>0$ on $(0,1]$. We also have $xa>1$ by
\eqref{an:amplitude}. We set
\begin{equation*}
 z=(1-x^2)W'-\frac{xW}{2},\qquad
 \Phi=x\left(1-\frac{3}{10}W^{-\frac{4}{3}}\right).
\end{equation*}
We compute
\begin{equation*}
 \begin{gathered}
 \mathcal Dz=\frac{5}{4}W^{\frac{7}{3}}\Phi,\qquad
 \left(\frac{z}{W}\right)'=-B-xa-\frac{1}{2}<0,\\
 \Phi'=1+\frac{3}{10}W^{-\frac{4}{3}}
 \left(\frac{4xa}{3}-1\right)>0,\qquad
 z(0)=W'(0),\quad z'(0)=0.
 \end{gathered}
\end{equation*}

We define
\begin{equation*}
 J=\mu^3(Wh'-W'h),\qquad H=\mu^3(zh'-z'h).
\end{equation*}
By the equations for $W,h,z$, we have
\begin{equation}\label{np:radial-wronskians}
 J'=-\frac{5}{4}\mu W^{\frac{7}{3}}h,\qquad
 H'=\Phi J',\qquad (H-\Phi J)'=-\Phi'J,
 \qquad J(1)=H(1)=0.
\end{equation}
We have $J>0$ near $1$. If $J$ had a last zero $s\in(0,1)$,
then $J'(s)\ge0$ and $h(s)\le0$. If $h(s)=0$, then $J(s)=0$
implies $h'(s)=0$, contrary to uniqueness for the linear equation.
Thus $h(s)<0$. We integrate the third identity to obtain
\begin{equation*}
 H(s)=\int_s^1\Phi'J\,dx>0,
 \qquad
 H(s)=-\mu(s)^3W(s)h(s)\left(\frac{z}{W}\right)'(s)<0,
\end{equation*}
a contradiction. Thus $J>0$ on $(0,1)$, and $h/W$ is strictly increasing.

Since $W(1)>1$ by \eqref{an:amplitude}, we choose $x_0$ with $W(x_0)=1$.
If $h\ge0$ on $[x_0,1]$, we integrate the Wronskian identity for
$h$ and $W^{\frac{7}{3}}-W$ to obtain
\begin{equation*}
 -\mu(x_0)^3W'(x_0)h(x_0)
 =\frac{7}{3}\int_{x_0}^1\mu^3W^{\frac{1}{3}}W'^2h\,dx>0,
\end{equation*}
which is impossible. Hence $h$ has exactly one zero
$\sigma\in(x_0,1)$, with $h<0$ on $(0,\sigma)$.
By \eqref{np:radial-wronskians}, we have
\begin{equation*}
 W'(0)h'(0)=H(0)
 =\Phi(\sigma)J(0)-\int_0^1(\Phi-\Phi(\sigma))J'\,dx>0.
\end{equation*}
Indeed, $J(0)\ge0$, $\Phi(\sigma)>0$, and the integrand is strictly
negative on both sides of $\sigma$. At zero, $\Phi J'=O(x^2)$,
so the integral converges. Since $h<0$ near zero and $h'(0)>0$,
we also have $h(0)<0$.

Since the potential is even and the regular solution space at each
endpoint is one-dimensional, a solution regular at both endpoints
must be even or odd. The inequalities $h(0)<0<h'(0)$ exclude both
possibilities for a nonzero solution.
\end{proof}

\subsection{The principal spectral level}

We define the principal spectral level by
\begin{equation}\label{an:rayleigh}
 \Lambda_*=-\inf_{h\ne0}
 \frac{\int_0^\pi\sin^2r\,(h'^2-\mathfrak V(w)h^2)\,dr}{\int_0^\pi\cos^2(\frac{r}{2})h^2\,dr},
\end{equation}
where the infimum is taken over the Friedrichs form domain, the
completion of $C_c^\infty(0,\pi)$ under the norm
\begin{equation*}
 \left(\int_0^\pi\left[\sin^2r\,h'^2
 +\cos^2\frac{r}{2}\,h^2\right]dr\right)^{\frac{1}{2}}.
\end{equation*}

\begin{lemma}[Principal level]\label{np:principal-level}
For this profile, \eqref{np:one-sided} has exactly one regular spectral
value $\Lambda>2$.
\end{lemma}
\begin{proof}
The degree-one mode
\begin{equation*}
 \zeta(r)=\frac{3}{2}w(r)\sin\frac{r}{2}-2w'(r)\cos\frac{r}{2}
\end{equation*}
solves \eqref{np:one-sided} at $\Lambda=2$ and is positive on
$(0,\frac{\pi}{2}]$. For $x=-\cos r\in(0,1)$, we have
\begin{equation*}
 \zeta(r)=\sin\frac{r}{2}
 \left(2(1-x)W'(x)-\frac{3}{2}W(x)\right).
\end{equation*}
At every zero, we use \eqref{an:W} to compute
\begin{equation}\label{np:zero-crossing}
 \frac{W}{W'}=\frac{4(1-x)}{3},\qquad
 \frac{\zeta'(r)}{\sin(\frac{r}{2})\sin r\,W'(x)}
 =-\frac{2+5(1-x)W^{\frac{4}{3}}}{2(1+x)}<0.
\end{equation}
Since $\zeta(\frac{\pi}{2})>0>\zeta(\pi)$, the mode has exactly one
interior zero, which we denote by $r_0$. We set
\begin{equation*}
 \mathfrak q[h]=\int_0^\pi\left[\sin^2r\,(h'^2-\mathfrak V(w)h^2)
                  +2\cos^2\frac{r}{2}\,h^2\right]dr.
\end{equation*}
On the codimension-one subspace $h(r_0)=0$, we have, by integration
by parts on the two nodal intervals and density,
\begin{equation*}
 \mathfrak q[h]=\int_0^\pi\sin^2r\,\zeta^2
                 \left|\left(\frac{h}{\zeta}\right)'\right|^2dr\ge0.
\end{equation*}
For $\phi\in C_c^\infty(0,\pi)$ with $\phi(r_0)>0$, we compute
\begin{equation*}
 \mathfrak q[|\zeta|+\varepsilon\phi]
 =-4\varepsilon\sin^2r_0\,|\zeta'(r_0)|\phi(r_0)
       +\varepsilon^2\mathfrak q[\phi]<0
\end{equation*}
when $\varepsilon>0$ is small. In $L^2(\cos^2(r/2)\,dr)$ the potential
term is a bounded form, since $\sin^2 r\le4\cos^2(r/2)$. The Friedrichs
operator is therefore bounded below. The codimension-one bound and the
negative test above give it a one-dimensional negative spectral subspace.
We choose
its eigenfunction nonnegative. ODE uniqueness gives strict positivity,
and this eigenfunction realizes the infimum in \eqref{an:rayleigh}. The finite form norm selects the regular
Frobenius branches at both endpoints, so this eigenvalue corresponds
to the unique regular spectral value $\Lambda>2$.
\end{proof}

\subsection{The upper spectral bound}

We use a positive comparison function, following \cite[Lemma~8.5]{NP26}.
We set $D_0=1-\frac{x^2}{3}$. Since
\begin{equation*}
 D-\frac{47}{48}D_0
 =\frac{(3x^2-\frac{7}{3})^2+\frac{5}{9}}{288}>0,\qquad
 10\left(\frac{63}{94}\right)^2<\frac{9}{2},
\end{equation*}
we obtain from \eqref{an:amplitude}, by oddness,
\begin{equation*}
 |W(x)|^{\frac{4}{3}}\le\frac{9}{2}\frac{|x|^{\frac{4}{3}}}{D_0(x)^2}.
\end{equation*}
We use the concavity of $z^{\frac{2}{3}}+\frac{z^2}{9}$ on $[0,1]$ to obtain
\begin{equation*}
 |x|^{\frac{4}{3}}\le\frac{2+8x^2-x^4}{9},\qquad
 \mathfrak V(W)\le V_+
 :=\frac{35(2+8x^2-x^4)}{32D_0(x)^2}-\frac{15}{16}.
\end{equation*}
We set
\begin{equation*}
 \Psi=e^{-x}\left(\frac{1-x}{D_0}\right)^{\frac{5}{2}},
\end{equation*}
and compute
\begin{equation}\label{an:upper-comparison}
 \frac{-(1-x^2)\Psi''+3x\Psi'-V_+\Psi+\frac{15\Psi}{1-x}}{\Psi}
 =\frac{Q(x)}{288D_0(x)^2}.
\end{equation}
For $0\le x\le1$, we compute
\begin{align*}
 Q(-x)&=\frac{179}{5}+830\left(\frac{4}{5}-x\right)^2
        +40(1-x)^2+548x(1-x)^2\\
 &\quad+(1-x)^4\bigl[32(1+x)^2+33\bigr]>0,\\
 Q(x)-Q(-x)&=64x(3-x^2)(11-2x^2)\ge0.
\end{align*}
Thus $Q>0$ on $[-1,1]$.

In the variable $x=\cos r$, we apply the ground-state identity with
$\Psi$ to compactly supported test functions. Since $\mathfrak V(W)\le V_+$,
we may use \eqref{an:upper-comparison} and pass by density to the
principal eigenfunction $h_*$ in the Friedrichs form domain. We obtain
\begin{equation*}
 (30-\Lambda_*)\int_{-1}^1\frac{\mu h_*^2}{2(1-x)}\,dx
 \ge\int_{-1}^1\frac{\mu Qh_*^2}{288D_0^2}\,dx>0.
\end{equation*}
Hence $\Lambda_*<30$.

\subsection{The lower spectral bound}

For the lower bound at $\Lambda=20$, we set
\begin{equation*}
 g=x^{\frac{3}{2}}\left(1+\frac{3x^2}{4}+\frac{5x^4}{8}\right).
\end{equation*}
Using \eqref{an:amplitude} and the identity
\begin{equation*}
 \left(\log\frac{b_+^{\frac{3}{2}}}{g}\right)'
 =\frac{x[42(3x^2-1)^2+6+8x^2+x^4+65x^4(1-x^2)]}{(4-x^2)(8-x^2)(8+6x^2+5x^4)}>0,
\end{equation*}
we see that $\frac{W^{\frac{3}{2}}}{g}$ is nondecreasing. We compute
\begin{equation*}
 U:=\frac{-\mu^2g''+3xg'}{g}
 =-\frac{3}{4x^2}-\frac{3}{4}
   +\frac{210x^4}{8+6x^2+5x^4}.
\end{equation*}
The last term is nonnegative and increasing, so we have
\begin{equation*}
 \begin{gathered}
 U+(1+x)U'\ge\frac{3(2+x)}{4x^3}-\frac{3}{4}\ge\frac{3}{2},\\
 \left[2(1+x)\left(\frac{29}{27}U+\frac{5}{4}\right)\right]'>0.
 \end{gathered}
\end{equation*}
In \eqref{an:rayleigh}, we take the test function
\[
 h(r)=
 \begin{cases}
 0,&0\le r<\pi/2,\\
 W(-\cos r)^{3/2},&\pi/2\le r\le\pi .
 \end{cases}
\]
Let $H=W^{3/2}$ on $[0,1]$. Testing \eqref{an:W} against $\mu W^2$ gives
\[
 \int_0^1\mu\bigl(\mathfrak V(W)H^2-\mu^2H'^2\bigr)\,dx
 =\int_0^1\mu\left(\frac{29}{27}\mu^2H'^2+\frac54H^2\right)dx .
\]
The ground-state identity gives $\int_0^1\mu^3H'^2\,dx\ge
\int_0^1\mu U H^2\,dx$, with equality when $H$ is replaced by $g$.
Set $d\nu=\frac{\mu g^2}{2(1+x)}\,dx$ and
$F=2(1+x)(\frac{29}{27}U+\frac54)$. Since $F$ and $(H/g)^2$ are
nondecreasing, the weighted Chebyshev inequality gives
\begin{equation*}
 \Lambda_*\ge\frac{\int F(H/g)^2\,d\nu}{\int(H/g)^2\,d\nu}
 \ge\frac{\int F\,d\nu}{\int d\nu}=
 \frac{\displaystyle\int_0^1\mu
       \left(\frac{29}{27}\mu^2g'^2+\frac{5}{4}g^2\right)dx}{\displaystyle\int_0^1\frac{\mu g^2}{2(1+x)}\,dx}.
\end{equation*}
Both test functions are $O(x^{\frac{3}{2}})$ at zero, and the flux weight $\mu^3$
vanishes at one. Thus the boundary terms vanish, and the increasing
function is integrable against both measures in the covariance identity.

Using $\frac{\mu}{1+x}=\frac{1-x}{\mu}$ and the elementary moments
\begin{equation*}
 \int_0^1x^{2j+1}\mu\,dx=\frac{2^{2j+1}j!(j+1)!}{(2j+3)!},\qquad
 \int_0^1\frac{x^{2j}}{\mu}\,dx
 =\frac{\pi}{2}\frac{(2j)!}{4^j(j!)^2},\qquad j\ge0,
\end{equation*}
we obtain, since $\pi>\frac{157}{50}$,
\begin{align*}
 \int_0^1\mu\left[\frac{29}{27}\mu^2g'^2
       +\left(\frac{5}{4}-\frac{10}{1+x}\right)g^2\right]dx
 &=\frac{543355\pi}{65536}-\frac{2546483}{98280}
 >\frac{1}{10}.
\end{align*}
Thus $\Lambda_*>20$, completing the proof of Lemma~\ref{np:certificate}.

\section*{AI disclosure}
In 2014, J.H.\ and X.C.~organized a seminar on soliton resolution,
initially focused on the nonlinear Schr\"odinger equation. Inspired by the
work discussed in Section~\ref{intro:related-results}, the authors asked
whether radiation could drive a stationary core along the conformal orbit
of a sign-changing state, producing recurrent changes of shape modulo
spatial similarities. We constructed a counterexample to the fixed-profile formulation for
the nonlinear Schr\"odinger equation, which is available upon request. We
then turned to the nonlinear wave equation, where we expected radiation to
have a weaker effect on the core.

The main ideas of this paper are due to the authors, including, in
particular, lengthening the intervals between paired impulses so that
progressively weaker radiation packets could move the core by a fixed
amount. A further insight is that the two limiting shapes belong to the
same conformal class. This observation led the authors to propose a
reformulation of soliton resolution in which each core may change shape
within a fixed conformal class. We also made progress toward establishing
this formulation in the nonradial one-core setting. We directed ChatGPT
Astra~6.2 to assist in deriving and writing out specific technical
estimates. We also used ChatGPT to help polish the language and check the proof.

% Omit the bibliography heading from the contents.
\begingroup
\makeatletter
\let\@tocwrite\@gobbletwo
\makeatother

\endgroup

\end{document}